\documentclass[10pt, a4paper, twoside]{article}

\usepackage{exscale, amsmath, amssymb, amsthm, dsfont, amsfonts, latexsym, url,color}
\usepackage[curve,matrix,arrow,ps]{xy}
\usepackage[hypertex]{hyperref}
\newtheorem{definition}{Definition}[section]
\newtheorem{theorem}[definition]{Theorem}
\newtheorem{proposition}[definition]{Proposition}
\newtheorem{lemma}[definition]{Lemma}
\newtheorem{corollary}[definition]{Corollary}

\newtheorem{remark}[definition]{Remark}

\newtheorem{deflemma}[definition]{Definition-Lemma}

\newtheorem{examples}[definition]{Examples}
\numberwithin{equation}{section}

\newcommand{\nd}{\noindent}

\newcommand{\dE}{{\mathds E}}
\newcommand{\dC}{{\mathds C}}
\newcommand{\dQ}{{\mathds Q}}

\newcommand{\dN}{{\mathds N}}

\newcommand{\dR}{{\mathds R}}

\newcommand{\dZ}{{\mathds Z}}
\newcommand{\dP}{{\mathds P}}

\newcommand{\dH}{{\mathds H}}

\newcommand{\de}{{\mathds 1}}

\newcommand{\cA}{\mathcal{A}}

\newcommand{\cC}{\mathcal{C}}
\newcommand{\cD}{\mathcal{D}}
\newcommand{\cE}{\mathcal{E}}
\newcommand{\cF}{\mathcal{F}}
\newcommand{\cG}{\mathcal{G}}
\newcommand{\cH}{\mathcal{H}}

\newcommand{\cK}{\mathcal{K}}
\newcommand{\cL}{\mathcal{L}}
\newcommand{\cM}{\mathcal{M}}
\newcommand{\cN}{\mathcal{N}}
\newcommand{\cO}{\mathcal{O}}

\newcommand{\cQ}{\mathcal{Q}}

\newcommand{\cT}{\mathcal{T}}
\newcommand{\cU}{\mathcal{U}}
\newcommand{\cV}{\mathcal{V}}
\newcommand{\cW}{\mathcal{W}}
\newcommand{\cX}{\mathcal{X}}
\newcommand{\cY}{\mathcal{Y}}

\newcommand{\SSC}{\scriptscriptstyle}

\newcommand{\cbl}{{\check D_{BL}}}
\newcommand{\cblp}{{\check D_{BL}^{pp}}}
\newcommand{\TERP}{(H,H'_\dR,\nabla,P,w)}
\newcommand{\VB}{\mathit{VB}}
\newcommand{\VBan}{\mathit{VB}^{an}}
\newcommand{\MBL}{M_{BL}}
\newcommand{\MBLfun}{\cM_{BL}}
\newcommand{\MBLr}{|M_{BL}|}
\newcommand{\MBLred}{M_{BL}^{red}}
\newcommand{\MBLp}{M^{pp}_{BL}}
\newcommand{\MBLpr}{|M^{pp}_{BL}|}

\newcommand{\LL}{\Lambda_b(W^\omega)}
\newcommand{\LLr}{|\Lambda_b(W^\omega)|}
\newcommand{\LLp}{\Lambda^{pp}_b(W^\omega)}
\newcommand{\LLpr}{|\Lambda^{pp}_b(W^\omega)|}

\newcommand{\cAn}{\cC^{an}}
\newcommand{\uul}{{\underline{l}}}
\newcommand{\uun}{{\underline{n}}}
\newcommand{\nnn}{\nabla}
\newcommand{\www}{\widetilde}
\newcommand{\whh}{\widehat}

\newcommand{\omm}{{\{0\}\times M}}

\DeclareMathOperator{\Spec}{\textup{Spec}\,}
\DeclareMathOperator{\Specan}{\textup{Specan}\,}
\DeclareMathOperator{\Proj}{\textup{Proj}\,}
\DeclareMathOperator{\Sp}{\textup{Sp}}
\DeclareMathOperator{\Spp}{\textup{Spp}}
\DeclareMathOperator{\KMSS}{{KMSS}}
\DeclareMathOperator{\DcPMHS}{\check{\mathit{D}}_{\mathit{PMHS}}}
\DeclareMathOperator{\DPMHS}{\mathit{D}_{\mathit{PMHS}}}
\DeclareMathOperator{\DcPHS}{\check{\mathit{D}}_{\mathit{PHS}}}
\DeclareMathOperator{\DPHS}{\mathit{D}_{\mathit{PHS}}}

\DeclareMathOperator{\Id}{\mathit{Id}}

\DeclareMathOperator{\Aut}{\textup{Aut}}
\DeclareMathOperator{\diag}{\textup{diag}}

\DeclareMathOperator{\codim}{codim}

\DeclareMathOperator{\Gl}{Gl}
\DeclareMathOperator{\Gr}{Gr}

\DeclareMathOperator{\id}{id}
\DeclareMathOperator{\Imm}{Im}

\DeclareMathOperator{\End}{\mathit{End}}

\begin{document}

\title{Limits of families of Brieskorn lattices \\ and compactified classifying spaces}

\author{Claus Hertling\and Christian Sevenheck}

\maketitle

\begin{abstract}
We investigate variations of Brieskorn lattices over
non-compact parameter spaces, and discuss the corresponding
limit objects on the boundary divisor. We study the associated
variation of twistors and the corresponding
limit mixed twistor structures. We construct a compact classifying space
for regular singular Brieskorn lattices and prove that its pure polarized part
carries a natural hermitian structure and that the induced distance
makes it into a complete metric space.
\end{abstract}

\tableofcontents

\renewcommand{\thefootnote}{}
\footnote{
This research is partially supported by ANR grant ANR-08-BLAN-0317-01 (SEDIGA).\\
C.H. acknowledges partial support by the ESF research grant MISGAM.}

\setcounter{section}{0}

\section{Introduction}
\label{secIntro}

This paper deals with a basic object attached to an isolated hypersurface singularity:
the Brieskorn lattice. It was introduced in \cite{Brie} in order
to understand the monodromy of the cohomology bundle of the Milnor fibration
of such a singularity, but it turned out that it contains much more information,
it is a highly transcendental invariant of the singularity.
Since \cite{SM3}, \cite{SM} and \cite{He2} it is evident that the Brieskorn
lattice is a very well suited object to study the Torelli problem for hypersurface singularities.

In various applications, one is interested not only in
local singularities but also in regular functions on affine manifolds
with isolated critical points. In this
case, one can also define a Brieskorn lattice, which contains
more information than the sum of the local Brieskorn lattices at the critical
points, in particular, its structure depends very much on the behavior of
the function at infinity.  In \cite{Sa2}, a precise condition, called cohomological
tameness for these functions is given which ensures that this
algebraic Brieskorn lattice is a free module over the ring of polynomial functions
on the base. However, as the dimension of the cohomology of the Milnor fibre of such a function
need not to be equal to the sum of the Milnor numbers of the critical points, it might
happen that the Brieskorn lattice has the ``wrong'' rank. In order to overcome this,
and for various other reasons, it is convenient to work with a twisted version of the Brieskorn lattice,
called Fourier-Laplace transformation. This transformation can also be done
in the local case, i.e. for the Brieskorn lattice of an isolated hypersurface singularity.
Alternatively, there is a direct description of this twisted object using Lefschetz thimbles and
oscillating integrals. This description makes the definition of a polarizing
form, given by the intersection form of Lefschetz thimbles
in opposite fibres, very transparent. It goes back to the work of Pham (\cite{Ph3,Ph4}),
a short version of it, which is also valid for families of Brieskorn lattices
can be found in \cite[section 8]{He4}.

The interest in studying the global situation of tame functions on affine manifolds
comes from the mirror symmetry phenomenon,
which relates in an intricate way data defined by such a polynomial function
(called B-model in physics) to data from symplectic geometry
(called A-model), namely, the quantum cohomology of some particular
symplectic manifolds. This correspondence
can be stated as an isomorphism of \emph{Frobenius manifolds} defined
by the two geometric inputs. On the B-side, the key tool to the construction of these
Frobenius structures is exactly the Fourier-Laplace transformation
of the Brieskorn lattice of the tame functions.

In both cases (local or global), the outcome of this construction (or of the direct approach via Lefschetz thimbles and
oscillating integrals) is an object which consists of a holomorphic
vector bundle on $\dC$, a flat connection on it having a pole of order at most two at zero and a pairing
between opposite fibres of that bundle with a prescribed pole order at zero.
The same kind of object exists for the A-model, called Dubrovin- or Givental connection, where
the flatness of the connection expresses all the properties of the quantum multiplication, in particular,
its associativity, which is equivalent to the WDVV-equation of the Gromov-Witten potential.

The B-model comes canonically equipped with some extra ingredient, namely, a real or even
integer structure of the flat bundle on $\dC^*$. It is the bundle generated by Lefschetz thimble
over $\dR$ resp. $\dZ$. The flat structure and the real structure make it possible to construct
a canonical extension of the bundle to $\dP^1$ such that connection and the pairing extend appropriately.
The result of this construction is what was called an \emph{integrable twistor structure} in
\cite[chapter 7]{Sa6}, and by generalizing it to the case
of families of Brieskorn lattices, Simpson's theory of harmonic bundles and variations of twistor structures
comes into play (see, e.g., \cite{Si2,Si4,Si5}). Notice however that the
harmonic bundles defined by this construction starting from a (Fourier-Laplace transform of a) Brieskorn
lattice carry additional structure, which were called $tt^*$-geometry by Cecotti and
Vafa (\cite{CV1, CV2}).

It is by no means evident to identify the real structure on the A-side, but in the recent
papers \cite{Ir1}, \cite{Ir2} and \cite{Ir3}, Iritani has made an important progress. He shows that
one might abstractly define real or integer structures of the A-model connection
which have a good behavior under rather mild conditions and he gives a concrete description
of the real/integer structure obtained by mirror symmetry for
toric orbifolds in terms of K-groups. As a consequence, one has, at least in favorite cases of examples,
the same structure on both sides, which is an analytic or formal object inducing a rich geometry on the parameter
spaces, which mixes in a subtle way holomorphic and anti-holomorphic data.
Hence it seems to be a good idea to formalize the setup, and study these
structures abstractly. This direction has been initiated in \cite{He4}, and
pursued in \cite{HS1}, \cite{HS2}. The geometric object sketched above
was called TERP-structure in these papers. This abbreviation
stands for ``twistor, extension, real structure and pairing''.
The main philosophy which we continue to exploit in this article is
that TERP-structures are an interesting generalization
of Hodge structures, and that one should try to generalize the known
results from Hodge theory to (variations of) TERP-structures.
In particular, the notion of pure resp. pure polarized TERP-structures
are defined in a natural way generalizing the corresponding notions
for Hodge structures.

Very recently, objects quite similar to TERP-structures
have been introduced and studied in \cite{KKP},
under the name ``non-commutative Hodge structures''. According to the main
conjecture of loc.cit., they arise as the cyclic homology
of certain categories, thought of as a ``non-commutative spaces''.
Via this construction, these structures also appear in the
homological mirror symmetry program.

%Examples include the
%derived category of quasi-coherent sheaves on a scheme or the category of maximal
%Cohen-Macaulay modules on a hypersurface germ.

There is an important difference between
the two above mentioned classes of examples: Whereas the
pole at the origin of the connection has order at most two in all cases
by definition, it defines a regular singularity in the sense of
\cite{De3} in the local case (i.e., in the case where the origin is the only
critical point of the function). We call the corresponding TERP-structures
regular singular. However, starting with a
tame polynomial function, the corresponding TERP-structure
will in general have a pole defining an irregular singularity at the origin.
The analysis of TERP-structures of this type, call
irregular, is more involved. One reason is that the connection defines,
besides the monodromy, the far more subtle Stokes structures, which
have to be taken into account. On the other hand, due to a recent, fundamental result of Sabbah
(\cite{Sa8}) we know that the TERP-structure of a tame polynomial
is always pure polarized, contrary to the local (i.e. regular singular) case.

This paper can be roughly divided into two parts. Whereas the first one
(sections \ref{secBasics} to \ref{secRigidity}) applies to arbitrary
(variations of) TERP-structures, the second one (sections \ref{secCompactClassSpace} to \ref{secExApp})
concerns mainly the regular singular case. Our main motivation
for the whole article is to develop a theory of period maps for variations of regular singular TERP-structures
(e.g., for $\mu$-constant deformations of isolated hypersurface singularities)
in a way similar to the usual study of variations of Hodge structures, as
in \cite{GSch1}, \cite{Sch}. A particularly powerful tool in this theory is the
use of hyperbolic complex analysis for horizontal maps to period domains.
In order to imitate this approach for variations of regular singular TERP-structures,
one needs appropriate targets for these period maps, i.e., classifying spaces
of such regular singular TERP-structures. These spaces have been defined
and studied for Brieskorn lattices in \cite{SM3} and \cite{He2}.
However, there is no discussion of the corresponding $tt^*$-geometry on
the classifying spaces in these papers, simply because there was no clean mathematical framework
for doing this at that time. The general theory of TERP-structures and the
relation to twistor structures and harmonic bundles is worked out in
\cite{He4} and \cite{HS1}. Moreover, we
showed in \cite{HS2} how to use the twistor construction to obtain
a hermitian metric on the pure polarized part of the classifying space.
We also calculated the holomorphic sectional curvature on horizontal
tangent directions, and proved its negativity. As in the case of Hodge structures,
this is one of the key results to apply
hyperbolic complex analysis for period maps defined by
variations of TERP-structures. However,
a crucial point was left open in that paper: this metric
we constructed on the classifying space is not complete in general.
The reason behind this fact is the following: We fixed the spectral pairs in order
to relate this classifying space to the classifying spaces of Hodge structures via
a construction modeled after Steenbrink's mixed Hodge structure on
the cohomology of the Milnor fibre of an isolated hypersurface singularity. But
it might very well happen that a special member of a variation of regular
singular TERP-structures has different spectral pairs
than the general member of this family. These ``missing limit points'' of
the classifying space prevent the metric from being complete.
In order to solve this problem, one is forced to look for a suitable
compactification of the classifying space, such that the hermitian
metric on its pure polarized part extends and yields a complete distance.
In particular, the spectral pairs must not be fixed for this larger space.
This is the basic idea behind the construction of the compact classifying space
in this paper: We fix an interval for the range of the spectrum, but not the spectral numbers themselves.
Then we can expect to capture the phenomenon of jumping spectrum.
The price we have to pay for this is that the space we obtain can be very singular.
However, we will show that the expected results hold: One can still define the pure polarized
part of this space, and the distance induced by the hermitian metric coming from the
twistor construction will be shown to be complete.

Let us give a short overview on this paper.
Following the general line of arguments in Hodge theory (see, e.g., \cite{Sch}),
we discuss in the first five sections of this paper the behavior of arbitrary (i.e., possibly irregular) families
of TERP-structures lattices at boundary points of the parameter spaces.
In section \ref{secBasics}, we briefly recall the necessary definitions from \cite{He4} and \cite{HS1}
and we give some more rather elementary properties of variations of
TERP-structures. In section \ref{secLimitTwistor}, we state the main
results of this first part. More precisely, given a variation of TERP-structures
on a complex manifold $Y$ which is the complement of a normal crossing divisor in a smooth ambient manifold
$X$, we give a precise condition (which we call tame) for the family to
have a limit object on the divisor. If we start with a pure polarized variation,
this condition corresponds to the tameness of the associated harmonic
bundle, as studied in \cite{Si2} and \cite{Mo2}. In particular,
pure polarized regular singular TERP-structures are always tame.
The first result is that the limit object is a family
of TERP-structures on the divisor. This allows us to
consider the associated twistor structure, and it turns
out that this is exactly Mochizuki's
limit polarized mixed twistor structure.
The proof of these two results is given in section \ref{secLimitProofs},
and relies on a general result concerning parabolic bundles.
This result has been proved in a slightly different context
by Borne (\cite{Borne1} and \cite{Borne2}), however, in order to make
the paper self-contained, we give in section \ref{secCompatible} an adapted version of Borne's
proof, which is also technically easier.

Section \ref{secRigidity} is an application of the results on extension of TERP-structures:
We prove a generalized version of a conjecture of Sabbah
concerning a rigidity property of integrable variations of twistor
structures on quasi-projective varieties with tame behavior at the boundary.
Although this seems to indicate that TERP-structures are not
more interesting than Hodge structures in this case, it is relevant as
it helps to understand the geometry in some examples, e.g. those coming from
quantum cohomology where a natural boundary divisor is given by the
so-called semi-classical limit. In particular, this result shows that a
variation on a quasi-projective manifold which is not of Hodge type
has necessarily boundary points which are not tame.

In the second part of the paper, namely in sections \ref{secCompactClassSpace} to \ref{secExApp} we construct
the above mentioned compact classifying spaces for regular singular TERP-structures
and we state and prove some of its crucial properties.
Section 7 gives the definition and the proofs of some basic properties and discusses
the relation between the classifying spaces from \cite{HS2} to the new one.
Section 8 is devoted to the construction of the hermitian metric
on the pure polarized part of the classifying space. We prove that the corresponding
distance is complete and study the action of a discrete group under a natural condition
satisfied in all geometric applications.
We finish the paper by discussing in section \ref{secExApp} in some detail the geometry of interesting
examples of these compact classifying spaces and we give
applications to period maps defined by variations of TERP-structures in subsection \ref{subsecApplications}.
They use both the limit objects discussed in the first part and the hyperbolicity
results from \cite{HS1} as well as the metric completeness of the pure polarized part
of the compact classifying space proved before.

\vspace*{0.5cm}

\textbf{Terminology and Notations:}

We will adopt the following convention for orderings and intervals:
We consider the natural ordering on $\dC$ given by
$c < d$ if either $\Re(c) <\Re(d)$ or $\Re(c) = \Re(d)$ and
$\Im(c) <\Im(d)$. Similarly, $c\leq d$ if $c < d$ (in the previous
sense) or if $c=d$. For any two complex numbers $c,d\in\dC$,
we define $(c,d)_\dC:=\{z\in \dC\,|\,c<z<d\}$ and similarly for closed
or half-open intervals. For any complex number $c\in\dC$, we write
$\lfloor c \rfloor $ for the largest integer $k$ such that
$k\leq c$ (i.e., such that $k\leq \Re(c)$).
For any two multi-indices $\mathbf{b},\mathbf{c} \in \dC^I$ we write
$\mathbf{b} < \mathbf{c}$ if $b_i<c_i$ for any $i\in I$, and
$\mathbf{b} \lneq \mathbf{c}$ if
$\mathbf{b}\leq\mathbf{c}$ and if there is $i\in I$
with $b_i< c_i$.

For any complex space $X$, we denote by $\VB_X$ the category
of $\cO_X$-locally free sheaves. If $X$ is a complex manifold, we write
$\VB_X^\nabla$ for the category of flat bundles
(or local systems) on $X$.
Occasionally, we work with the sheaf $\cAn_X$ of real analytic functions
on a complex space $X$ and the category of coherent $\cAn_X$-modules.
If $E$ is locally free over $\cAn_X$, then we write $E\in\VBan_X$.
We denote by $\cO_{\dP^1}\cC^{an}_X(k,l)$ the sheaf of real analytic functions on $\dC^*\times X$ which are holomorphic
in the $\dP^1$-direction (i.e., annihilated by $\partial_{\overline{z}}$, where
$z$ is the coordinate on $\dC$)
and which have at most poles of order $k$ resp. $l$ along $\{0\}\times X$ resp. $\{\infty\}\times X$.

For a complex manifold $X$, we denote by $\overline{X}$ its conjugate, which is
the same $\cC^\infty$-manifold, and where we put $\cO_{\overline{X}}:=\overline{\cO}_X$.
In particular, given a holomorphic bundle $E$ on $X$, the bundle $\overline{E}$ with the conjugate
complex structure in the fibres is holomorphic over $\overline{X}$.

For the reader's convenience, we collect here the definition of some maps that will be used
at several places in the paper.
$$
\begin{array}{c}
i:\dC^*\hookrightarrow\dC
\;\; ; \;\;\;\;
\widetilde{i}:\dC^*\hookrightarrow\dP^1\backslash\{0\}
\;\; ; \;\;\;\;
\widehat{i}:\dC^*\hookrightarrow\dP^1
\;\; ; \;\;\;\;
j:\dP^1\rightarrow \dP^1\;,\;\; j(z)=-z;
\\ \\
\gamma:\dP^1\rightarrow \dP^1\;,\;\; \gamma(z)=\overline{z^{-1}}
\;\; ; \;\;\;\;
\sigma:\dP^1\rightarrow \dP^1\;,\;\; \sigma(z)=-\overline{z^{-1}}.
\end{array}
$$

\vspace*{1cm}

\section{Definition and basic properties of TERP-structures}
\label{secBasics}

We start by recalling the definition of variations of TERP-structures,
their associated topological data, the construction of twistors from them,
and the special case of regular singular TERP-structures.
The main references are \cite{He4} and \cite{HS1}.
A small generalization is the notion of families of TERP-structures
on an arbitrary complex space (possibly non-reduced), this will be
needed in the discussion of classifying spaces in section
\ref{secCompactClassSpace}. Moreover, we give a translation of the notion
of a polarized mixed twistor structure to our frame in
lemma \ref{lemTERPWeightFilt}.

\begin{definition}\label{defFamilyTERP}
Let $X$ be a complex space. A family of TERP-structures
of weight $w\in\dZ$ on $X$ consists of the following data.
\begin{enumerate}
\item
A holomorphic vector bundle $H$ on $\dC\times X$, i.e., the linear space
associated to a locally free sheaf $\cH$ of $\cO_X$-modules.
\item
A flat structure on the restriction
$H':=H_{|\dC^*\times X}$ to $\dC^*\times X$, i.e.,
the transition functions between two local trivializations
of the map $H'\rightarrow \dC^*\times X$ are constant.
Moreover, we require that for any $t\in X$, the flat connection on $H'_{|\dC^*\times \{t\}}$
extends to a meromorphic connection on $H_{|\dC\times\{t\}}$
with a pole of order at most 2. The local system associated to $H'$ will
be denoted by $(H')^\nabla$.
\item
A flat real subsystem $H'_\dR\subset H'$ of maximal rank.
\item
A non-degenerate, $(-1)^w$-symmetric pairing
$P:\cH\otimes j^*\cH\rightarrow z^w\cO_{\dC \times X}$
which is flat on $H'$ and which takes values in $i^w\dR$ on $H'_\dR$.
Here non-degenerateness along $\{0\}\times X$ means that
the induced symmetric pairing
$[z^{-w}P]:\cH/z\cH\otimes \cH/z\cH \rightarrow \cO_X$ is non-degenerate.
\end{enumerate}
If $X$ is smooth and if the flat connection on $H'$ extends to
a meromorphic connection
$$
\nabla:\cH\longrightarrow \cH\otimes z^{-1}
\Omega^1_{\dC\times X}(\log (\{0\}\times X)),
$$
of type 1 then $\TERP$ is called a variation of TERP-structures.
A single TERP-structure is a family
on $X=\{pt\}$.
\end{definition}

We give two simple examples of variations of
TERP-structures, which will be used later (see
examples \ref{exTERP-2}, \ref{exTERP-NonTame-1} and subsection \ref{subsecRedSpaces}).
Some of the interesting phenomena that may occur for general
variations are already visible here.
Many more examples will be given in section \ref{secExApp}.
\begin{examples}\label{exTERP-1}
\begin{enumerate}
\item
Consider a trivial bundle $H'$ of rank 2 on
$\dC^*\times X=\dC^*\times \dP^1$
with two generating flat sections $A_1$ and $A_2$.
The flat real structure is defined such that $\overline{A_1}=A_2$,
the pairing is defined by
$P(A_i(z,r),A_j(-z,r))=\varepsilon\cdot\delta_{i+j,3}$
for $(z,r)\in \dC^*\times \dP^1$, with $\varepsilon=\pm 1$ fixed (so in fact
this gives two examples, one for $\varepsilon=1$, one for
$\varepsilon=-1$).
The extension $H$ to $\dC\times\dP^1$ is defined by
\begin{eqnarray*}
\cH &=&\cO_{\dC\times\dC}\cdot (z^{-1}A_1+rA_2)\oplus
\cO_{\dC\times\dC}\cdot zA_2 \quad\textup{ on }\dC\times \dC,
\quad \textup{ and by }\\
\cH &=&\cO_{\dC\times(\dP^1\,\backslash\,\{0\})}\cdot (z^{-1}r^{-1}A_1+ A_2)\oplus
\cO_{\dC\times(\dP^1\,\backslash\,\{0\})}\cdot A_1
\quad\textup{ on }\dC\times (\dP^1\,\backslash\,\{0\}).
\end{eqnarray*}
Here we write $zA_2$ for the section $(z\mapsto zA_2)$.
This gives a variation of TERP-structures of weight 0 on $\dP^1$.

\item
The bundle $H'$, the flat sections $A_1$ and $A_2$,
the bundle $H$ and the pairing $P$ are as in (a).
Here $\varepsilon=1$ is chosen.
But the real structure is changed to $\overline{A_i}=A_i$.
Again this gives a variation of TERP-structures of weight 0 on $\dP^1$.
\end{enumerate}
\end{examples}

The next definition introduces several basic linear algebra objects defined
by a TERP-structure. A more detailed discussion is contained in \cite{He4} and
\cite{HS1}.
\begin{definition}\label{defTopData}
(Topological data of a family of TERP-structures)
Let $\TERP$ be a family of TERP-structures on a complex space $X$.
We denote by $H^\infty$ the vector space of multivalued flat sections of the
local system $(H')^\nabla$,
and by $H^\infty_\dR$ the subspace of real flat multivalued sections.
Let $\pi_1(\dC^*\times X)\rightarrow \Aut(H^\infty_\dR)$
be the monodromy representation associated to $(H'_\dR)^\nabla$, and
denote its image by $\Gamma$.
We write $M_z\in\Gamma$ for the automorphism corresponding
to a (counter-clockwise) loop
around the divisor $\{0\}\times X$.
We decompose $M_z$ as $M_z=(M_z)_s\cdot (M_z)_u$
into semi-simple and unipotent part.
Let $H^\infty:=\oplus H^\infty_\lambda$
be the decomposition into eigenspaces with respect to $(M_z)_s$,
put $H^\infty_{\arg=0}:=\oplus_{\arg\lambda=0} H^\infty_\lambda$,
$H^\infty_{\arg \neq 0}:=\oplus_{\arg\lambda \neq 0} H^\infty_\lambda$,
and let $N_z:=\log((M_z)_u)$
be the nilpotent part of $M$.

Finally, we denote by $S$ the non-degenerate and monodromy invariant form
on $H^\infty$ which is defined in \cite[formula (5.1)]{HS1}. It is
$(-1)^w$-symmetric on $H^\infty_{\arg=0}$ and
$(-1)^{w-1}$-symmetric on $H^\infty_{\arg\neq 0}$.
We also point the reader to the formulas \cite[(5.4) and (5.5)]{HS1}
which connect $S$ and $P$ and which will be used in the examples
in section \ref{secExApp}.
We call the tuple $(H^\infty, H^\infty_\dR, \Gamma, M_z, S, w)$
the topological data of the family
$\TERP$.
\end{definition}

As we already pointed out in the introduction, TERP-structures
are closely related to twistor structures, i.e. holomorphic
bundles over $\dP^1$. This is shown in the following definition.
\begin{definition}[Extension to infinity]\label{defExtInfinity}
Consider a family of TERP-structures $(H,H'_\dR,\nabla,P,w)$ over
a complex space $X$.
Let $\gamma:\dP^1\times X\rightarrow\dP^1\times X;
(z,t)\mapsto (\overline{z}^{-1},t)$.
\begin{enumerate}
\item
Define for any $(z,t)\in\dC^*\times X$
the following two anti-linear involutions.
$$
\begin{array}{rcl}
\tau_{real}: H_{z,t} & \longrightarrow & H_{{\gamma(z)},t}\\
s & \longmapsto & \nabla\textup{-parallel transport of }\overline{s},\\ \\
\tau: H_{z,t} & \longrightarrow & H_{{\gamma(z)},t}\\
s & \longmapsto & \nabla\textup{-parallel transport  of }\overline{z^{-w}s}.
\end{array}
$$
$\tau_{real}$ is flat.
The induced maps on sections by putting
$s\mapsto\left(z\mapsto\tau s(\overline{z}^{-1})\right)$ resp.
$s\mapsto\left(z\mapsto\tau_{real}s(\overline{z}^{-1})\right)$
will be denoted by the same letter.
They can either be seen as morphisms
$\tau,\tau_{real} :\cH' \rightarrow \overline{\gamma^*\cH'}$
which fix the base,
or as morphisms $\tau,\tau_{real}:\cH'\rightarrow\cH'$
which map sections in $U\subset \dC^*\times X$ to
sections in $\gamma(U)\subset \dC^*\times X$.
Note that for each fixed $t\in X$,
due to the two-fold conjugation (in the base and in the fibres),
$\tau$ and $\tau_{real}$ are morphisms of holomorphic bundles over $\dC^*$,
but that with respect to $X$ they are only real analytic morphisms.
Denote by $\widehat{H}$ the bundle obtained by patching
$\cH$ and $\overline{\gamma^*\cH}$ via the identification $\tau$.
It is a real analytic bundle whose restriction to $\dP^1\times\{t\}$
has a holomorphic structure for each $t\in X$.

\item
Define a sesquilinear pairing $\whh S: \cH'\otimes \overline{\sigma^*\cH'} \to \cO_{\dC^*}\cC^{an}_X$ by
\begin{eqnarray*}
\whh S:H_{z,t}\times H_{\sigma(z),t}&\to&\dC
\quad\quad\quad\textup{ for } (z,t)\in \dC^*\times X,\\
(a(z,t),b(\sigma(z),t))&\mapsto& z^{-w}P(a,\tau(b))
=(-1)^wP(a,\tau_{real}(b)).
\end{eqnarray*}
It is non-degenerate, flat and holomorphic with respect to $z$.
\end{enumerate}
\end{definition}

\begin{lemma}\label{lemGeneralitiesTERP}
Consider a single TERP-structure $\TERP$.
Let $\mu$ be the rank of $H$.
\begin{enumerate}
\item
The bundle $\whh H$ has degree zero.
The flat connection has
a pole of order at most 2 at $\infty$.
The pairing $P$ extends to a non-degenerate pairing
$P:\whh\cH\otimes j^*\whh\cH\to z^w\cO_{\dP^1}.$
By definition of $\whh H$,
$\tau(\whh \cH(U)) = \whh \cH(\gamma(U))$
for any subset $U\subset\dP^1$.
\item
The pairing $\whh S$ extends to a non-degenerate hermitian pairing
$\whh S:\whh\cH\otimes \overline{\sigma^*\whh\cH}\to \cO_{\dP^1}.$
It satisfies $\whh S(za,\sigma(z)b)=-\whh S(a,b)$ for
$a\in H_z,b\in H_{\sigma(z)}$, $z\in \dC^*$.
\item
The morphism $\tau$ acts on the space $H^0(\dP^1,\whh\cH)$
as an antilinear involution. The pairing $z^{-w}P$ has constant values on this space
and is symmetric, the pairing $h:=\whh S$ has also constant values on it and is
hermitian.
\item
Choose sections $v_1,\ldots,v_\mu$ of $\widehat{\cH}_{|\dC}$ such that $\widehat{\cH}
=\oplus_{i=1}^\mu \cO_{\dP^1}(0,k_i)\cdot v_i$ and such that $k_1\leq\ldots\leq k_\mu$.
Then $k_i=-k_{\mu+1-i}$. If $k_i+k_j > 0$, then $z^{-w}P(v_i,v_j)=0$ and
$\widehat{S}(v_i,v_j)=0$. The radicals of $z^{-w}P$ and of $h$ in $H^0(\dP^1,\widehat{\cH})$
are both equal to $H^0(\dP^1,\oplus_{i:k_i>0}\cO_{\dP^1}(0,k_i)\cdot v_i)$.
\item
If $H^0(\dP^1,\whh\cH)$ contains $\mu$ global sections
$v_1,\ldots,v_\mu$ such that $(h(v_i,v_j))\in \textup{GL}(\mu,\dC)$,
then the bundle $\whh H$ is trivial. Conversely,
if the bundle is trivial, then $h$ is non-degenerate.
\end{enumerate}
\end{lemma}

\begin{proof}
\begin{enumerate}
\item
\cite[Lemma 3.3]{HS1}.
\item
That $\whh S$ extends to a non-degenerate pairing on $\whh H$,
follows from 1. and from the definition of $\whh S$.
That it is hermitian, follows from the
calculation \cite[(3.4)]{HS1}. And
$$
\widehat{S}(za,\sigma(z)b)=z\cdot \overline{\sigma(z)}\cdot \widehat{S}(a,b)=-\widehat{S}(a,b).
$$
\item
The statement on $\tau$ follows from 1. and from $\tau^2=\id$.
The pairings $z^{-w}P$ and $\whh S$ both take values in $\cO_{\dP^1}$
and thus take constant values on global sections.
$\whh S$ is hermitian by 2., $z^{-w}P$ is symmetric because
$P$ is $(-1)^w$-symmetric on $\cH\otimes j^*\cH$.
\item
Consider the dual bundle $\widehat{\cH}^*:={\cH}\!om_{\cO_{\dP^1}}(\widehat{\cH},\cO_{\dP^1})$.
It is isomorphic to $j^*\widehat{\cH}$ via the non-degenerate pairing $z^{-w}P$. On the other
hand, $j^*\widehat{\cH}$ is non-canonically isomorphic to $\widehat{\cH}$, thus to
the direct sum $\oplus_{i=1}^\mu \cO_{\dP^1}(k_i)$. Obviously, $\widehat{\cH}^*$ is
isomorphic to $\oplus_{i=1}^\mu \cO_{\dP^1}(-k_i)$, so that $k_i=-k_{\mu+1-i}$.

Suppose that $k_i+k_j>0$. The function $z \mapsto z^{-w}P(v_i,v_j)$ is holomorphic on $\dP^1$.
From
$$
z^{-w}P(v_i,v_j) = z^{-k_i-k_j}z^{-w}P(z^{k_i}v_i,z^{k_j}v_j)
$$
we conclude that it vanishes at infinity and hence globally on $\dP^1$. Using this,
the non-degenerateness of $z^{-w}P$ and the symmetry $k_i=-k_{\mu+1-i}$ just shown we
obtain the following: For any $k\in\dZ$, the $z^{-w}P$-orthogonal complement of
$\oplus_{i:k_i\geq k} \cO_{\dP^1}(0,k_i)\cdot v_i$ is $\oplus_{j:k_j> -k}\cO_{\dP^1}(0,k_j)\cdot v_j$.
In particular, the radical of $z^{-w}P$ in $H^0(\dP^1,\widehat{\cH})$ is
$H^0(\dP^1, \oplus_{i:k_i> 0}\cO_{\dP^1}(0,k_i)\cdot v_i)$. The statements on $\widehat{S}$
and $h$ follow from those on $z^{-w}P$ and the following fact: For any $k\in \dZ$, the subbundle
$\oplus_{i:k_i\geq k}\cO_{\dP^1}(0,k_i)\cdot v_i$ is independent of the choice of the sections
$v_i$, and hence mapped to itself by the morphism $\tau$.
\item
If $\whh H$ is not trivial then by 4. we have that $\codim\textup{Rad}(h)<\mu$,
so $\mu$ global sections $v_i$ with
$(h(v_i,v_j))\in \textup{GL}(\mu,\dC)$
cannot exist. If $\whh H$ is trivial then $h$ is non-degenerate
because $\whh S$ is non-degenerate.
\end{enumerate}
\end{proof}

\begin{definition}\label{defPurePolarized}
A TERP-structure is called pure iff the bundle $\widehat{H}$ is trivial.
A pure TERP-structure is called polarized iff the hermitian
form $h$ is positive definite.
\end{definition}
Notice that lemma \ref{lemGeneralitiesTERP}, 5., gives an efficient criterion
to detect whether a given TERP-structure is pure.

For the discussion of regular singular TERP-structures
we will need {\it elementary sections} and the $V$-filtration
(also called Malgrange-Kashiwara filtration or Deligne extensions).
Let $\TERP$ be a family of TERP-structures on a complex
space $X$. For $U\subset X$ open and simply connected, denote
by $H^\infty(U)$ the space of global multivalued flat sections
on $H'_{|\dC^*\times U}$. Then for any $A\in H^\infty(U)_{e^{-2\pi i\alpha}}$,
the section
$$
es(A,\alpha):=z^{\alpha-\frac{N_z}{2\pi i}}A
$$
is holomorphic on $\dC^*\times U$ and is called an elementary section of order $\alpha\in \dC$. It satisfies
\begin{eqnarray}\label{eqElementary-z-1}
(z\nabla_z -\alpha)es(A,\alpha) &=& es(\frac{-N_z}{2\pi i}A,\alpha),\\
\label{eqElementary-z-2}
\tau(es(A,\alpha))& = & es(\overline{A},w-\overline{\alpha}).
\end{eqnarray}
The extension of $\cH'$ to $\{0\}\times X$ which is generated
by such sections of order at least $\alpha$
is called $\cV^\alpha$. Similarly we have $\cV^{>\alpha}$, which is
generated by elementary sections of order bigger than $\alpha$, and
$\cV^{>-\infty}$, which is generated by elementary sections of arbitrary order.
$\cV^{\alpha}$ and $\cV^{>\alpha}$ are locally free $\cO_{\dC\times X}$-modules,
$\cV^{>-\infty}$ is a locally free
$\cO_{\dC\times X}(*\{0\}\times X)$-module.
\begin{deflemma}
\label{defHodgeSpectrum}
\begin{enumerate}
\item
A single TERP-structure is called regular singular
if $\cH\subset \cV^{>-\infty}$.
\item
Let $\TERP$ be a single regular singular TERP-structure.
The $V$-filtration on $\cH$ induces a filtration in $H^\infty$:
We put for any $\alpha\in(0,1]_\dC$
$$
F^pH^\infty_{e^{-2\pi i \alpha}}:=
z^{p+1-w-\alpha+\frac{N}{2\pi i}}Gr_\cV^{\alpha+w-1-p}\cH.
$$
A twisted version of this filtration,
which was considered in \cite{He4,HS1,HS2} is defined
as $\widetilde{F}^\bullet:=G^{-1}F^\bullet$.
Here $G\in\Aut(H^\infty)$ is a certain automorphism of $H^\infty$,
defined by \cite[section 5]{HS1}. Its definition
is motivated by the  Fourier-Laplace transformation, and it is useful while
comparing $S$ on $H^\infty$ (see \ref{defTopData}) and $P$ on $\cH'$.
$G$ induces the identity on $\Gr^W_\bullet(H^\infty)$
where $W_\bullet$ is the weight filtration of $N_z$,
centered around 0.
\item
A regular singular TERP-structure $\TERP$ of weight $w$ is called mixed if
the tuple
$$
(H^\infty_{\arg\neq 0}, (H^\infty_{\arg\neq 0})_\dR, -N, S, \widetilde{F}^\bullet)
\mbox{ resp. }
(H^\infty_{\arg=0}, (H^\infty_{\arg=0})_\dR, -N, S, \widetilde{F}^\bullet)
$$
is a polarized mixed Hodge structure of weight $w-1$ resp. of weight $w$.
We refer to \cite{He4} or \cite{HS1} for the notion of
a polarized mixed Hodge structure (PMHS for short) used here. It is $(M_z)_s$-invariant,
and the eigenvalues of $(M_z)_s$ are automatically elements in $S^1$, so that
$H^\infty_{\arg=0}=H^\infty_1$ and $H^\infty_{\arg\neq 0}=H^\infty_{\neq 1}$ in this case
(see \cite[lemma 5.9]{HS1}).
\item
The spectrum $\Sp(H,\nabla)$ of the regular singular TERP-structure is defined by
$\Sp(H,\nabla) = \sum_{\alpha \in\dQ}d(\alpha)\cdot\alpha\in\dZ[\dC]$
where
$$
d(\alpha):=\dim_\dC\left(\frac{Gr^\alpha_\cV \cH}{Gr^\alpha_\cV z\cH}\right)
=\dim_\dC\Gr_F^{\lfloor w-\alpha\rfloor} H_{e^{-2\pi i\alpha}}^\infty.
$$
It is a tuple of $\mu$ complex numbers
$\alpha_1 \leq \ldots \leq \alpha_\mu$
with the symmetry property
$\alpha_i+\alpha_{\mu+1-i}=w$ (see, e.g., \cite[lemma 6.3]{HS1}).
By definition, $d(\alpha)\neq 0$ only if $e^{-2\pi i \alpha}$
is an eigenvalue of $M_z$.
In most applications the eigenvalues of $M_z$
are roots of unity so that the spectrum
actually lies in $\dZ[\dQ]$.
\item
The spectral pairs $\Spp(H,\nabla)$ of a regular singular TERP-structure
are a finer invariant than the spectrum itself. They are defined as follows.
\begin{eqnarray*}
\Spp &=& \sum d(\alpha,l)\cdot (\alpha,l)\in \dZ[\dC\times \dZ],\\
d(\alpha,l)&=&\dim \Gr_F^{\lfloor w-\alpha \rfloor }\Gr^W_{l-(w-1)}H^\infty_{e^{-2\pi i \alpha}}.
\end{eqnarray*}
The second entries are symmetric around $w-1$. The shift of $W$ by $w-1$ is adapted
to a PMHS as in 3. on $H^\infty_{\arg\neq 0}=H^\infty_{\neq 1}$, but not to a
PMHS on $H^\infty_{\arg=0}$, which would require a shift by $w$.
Notice also that this definition is shifted by $+1$ in the first entry
compared to the original definition in \cite{St} for Brieskorn lattices
of hypersurface singularities. Up to this shift, the current definition is also compatible
with \cite[chapter 4]{He2} if $w-1=n$.
\end{enumerate}
\end{deflemma}

This construction of a filtration $F^\bullet$ was first considered
by Varchenko for a Brieskorn lattice of an isolated hypersurface singularity
(which becomes part of a TERP-structure only after a Fourier-Laplace
transformation).

We will continue the discussion of families of regular singular
TERP-structures in section \ref{secCompactClassSpace}.
In the sequel, we state and prove a rather elementary lemma and
return afterwards to the examples considered above.

\begin{lemma}\label{lemRegSingBounded}
Let $\TERP$ be a family of regular singular TERP-structures
on a complex space with finitely many components.
Then there exists $\beta$ with
$\cV^\beta \supset \cH \supset \cV^{>w-1-\beta}$.
\end{lemma}

\begin{proof}
Consider an open and simply connected set $\Delta_\varepsilon\times U\subset\dC\times X$
such that $\cH_{|\Delta_\varepsilon\times U}$ is free with generating
sections $\sigma_1,\ldots,\sigma_\mu$. Choose a basis $(A_k)$ of $H^\infty$,
where $A_k\in H^\infty_{e^{-2\pi i \beta_k}}$ and $\beta_k\in(0,1]+i\dR$.
Then these generating sections can be written in the following way
$$
\sigma_j = \sum_{k=1}^\mu
\left(\sum_{l\in\dZ}\kappa(j,k,l)z^l\right)
es(A_k,\beta_k)=\sum_{k=1}^\mu\sum_{l\in\dZ}\kappa(j,k,l)es(A_k,\beta_k+l),
$$
where $\kappa(j,k,l)\in\cO_X(U)$.
If there were an infinite sequence
$(j_i,k_i,l_i)$ with $l_i\to-\infty$ and $\kappa(j_i,k_i,l_i)\neq 0$,
then outside of a union of countable many hypersurfaces in $U$
all these coefficients would be non-vanishing,
and the TERP-structures on this subset of $U$ would not be regular singular.
Therefore there exists $\beta_U$ with
$\cH_{|\dC\times U}\subset \cV^{\beta_U}_{|\dC\times U}$.
This inclusion  extends to all components of $X$ which meet $U$.
As $X$ has only finitely many components, one can choose
such a set $U$ for each of them. Then $\cH\subset \cV^\beta$
for a suitable $\beta$.

The other inclusion $\cH\supset \cV^{>w-1-\beta}$
uses properties of the pairing $P$.
We write $P=\sum_{k\in\dZ}z^k\cdot P^{(k)}$
with pairings $P^{(k)}:\cV^{>-\infty}\otimes j^*\cV^{>-\infty}\to \cO_X$.
The inclusion $\cH\supset \cV^{>w-1-\beta}$
follows immediately from $P^{(w-1)}(\cV^\beta,\cV^{>w-1-\beta})=0$
and $\cV^\beta\supset \cH$ and the next claim.

\medskip
{\bf Claim:} $\cH =\{\sigma\in \cV^{>-\infty}\ |\ P^{(w-1)}(\cH,\sigma)=0\}.$

\medskip
The inclusion $\subset$ is part of the definition of a family of
TERP-structures.
For the proof of $\supset$ we consider a germ $(X,x)$ of a complex space.
Let $v_1,...,v_\mu$ be an $\cO_{\dC\times X,(0,x)}$-basis of $\cH$
and let $v_1^*,...,v_\mu^*$ be the basis of $\cH$ with
$P^{(w)}(v_i,v_j^*)=\delta_{ij}$. Then $v_1,...,v_\mu$ is an
$\cO_{\dC\times X,(0,x)}[z^{-1}]$-basis of $\cV^{>-\infty}$, so any
$\sigma\in \cV^{>-\infty}$ can be written as
$\sigma=\sum_{j=1}^\mu\sum_{k\in\dZ}z^k\cdot\kappa_{j,k}\cdot v_j$
with unique coefficients $\kappa_{j,k}\in \cO_{X,x}$.
Now
$$
P^{(w-1)}(z^kv_j^*,\sigma)=P^{(w)}(z^{k+1}v_j^*,\sigma)=
(-1)^{k+1}\cdot \kappa_{j,-k-1}.
$$
This shows the claim.
\end{proof}

We return to the examples in \ref{exTERP-1} and describe the corresponding
twistors and the associated hermitian metrics in case that they are pure.
\begin{examples}\label{exTERP-2}
The TERP-structures from \ref{exTERP-1} are regular singular. In both examples
the spectral numbers are $(\alpha_1,\alpha_2)=(-1,1)$ for
$r\in\dC$ and $(\alpha_1,\alpha_2)=(0,0)$ for $r=\infty$.
We put $v_1:=z^{-1}A_1+rA_2$ (for $r\neq \infty$).
\begin{enumerate}
\item
In example 1., equation \eqref{eqElementary-z-2} yields
\begin{eqnarray*}
&&\tau(v_1)=zA_2+\overline{r}A_1 \quad (\textup{for }r\neq\infty),
\quad \tau(zA_2)=z^{-1}A_1, \quad \tau(A_1)=A_2,\\
&&H^0(\dP^1,\whh \cH(r))=\dC \cdot v_1\oplus\dC \cdot \tau(v_1)
\qquad \textup{ for }r\neq\infty,\\
&&H^0(\dP^1,\whh \cH(r))=\dC \cdot r^{-1}v_1\oplus\dC \cdot \tau(r^{-1}v_1)
\qquad \textup{ for }r\neq 0.
\end{eqnarray*}
The metric $h$ with respect to these two bases is given by the matrices
$\varepsilon\cdot(|r|^2-1)\cdot\begin{pmatrix}1&0\\0&1\end{pmatrix}$
resp.
$\varepsilon\cdot(1-|r|^{-2})\cdot\begin{pmatrix}1&0\\0&1\end{pmatrix}$
for $r\neq\infty$ resp. $r\neq 0$.
Therefore the TERP-structures are pure for $|r|\neq 1$
and either polarized for $r\in\Delta$ or for $r\in \dP^1 \backslash \overline\Delta$,
depending on the choice of $\varepsilon=\pm 1$.
For $|r|=1$, $\whh {\cH}(r)\cong \cO_{\dP^1}(-1)\oplus\cO_{\dP^1}(1).$
\item
In example 2. we have
\begin{eqnarray*}
&&\tau(v_1)=zA_1+\overline{r}A_2 (\textup{ for }r\neq\infty),
\quad \tau(zA_2)=z^{-1}A_2, \quad \tau(A_1)=A_1,\\
&&\tau(v_2)=v_2 \quad\textup{ with }\quad
v_2=r^{-1}z^{-1}A_1+A_2+{\overline r}^{-1}zA_1 \qquad
\textup{ for }r\neq 0,\\
&&H^0(\dP^1,\whh \cH(r))=\dC \cdot A_1\oplus\dC \cdot v_2
\qquad \textup{ for }r\neq 0,\\
&&H^0(\dP^1,\whh \cH(0))=\dC \cdot z^{-1}A_1\oplus\dC \cdot A_1
\oplus \dC\cdot zA_1.
\end{eqnarray*}
For $r=0$, $\whh {\cH}(r)\cong \cO_{\dP^1}(-2)\oplus\cO_{\dP^1}(2).$
For $r\neq 0$, the TERP-structure is pure, and the matrix of the
metric $h$ with respect to the basis $(A_1,v_2)$ is
$\begin{pmatrix}0&1\\1&0\end{pmatrix}$, so the signature is $(1,1)$.
\end{enumerate}
\end{examples}

The following lemma translates the notion of a polarized
mixed twistor structure into our setting of TERP-structures.

\begin{lemma}\label{lemTERPWeightFilt}
Let $\TERP$ be a variation of TERP-structures on a
manifold $M$. Additionally, let
$N:H'_\dR\to H'_\dR$ be a nilpotent flat infinitesimal isometry
of $P$, i.e. $P(Na,b)+P(a,Nb)=0$.

\begin{enumerate}
\item
(Topological part)
Let $W'_\bullet$ be the weight filtration (centered at 0) of $N$ on $H'$.
Any $W'_l$ is a flat subbundle with real structure.
Moreover, the quotients $\Gr^{W'}_l$ are also flat bundles with real structure.
The pairing $P$ has the following properties:
\begin{eqnarray*}
&& P: \cW'_{-l}\otimes j^*\cW'_{l-1}\to 0,\\
&& P:\Gr^{\cW'}_{-l}\otimes j^*\Gr^{\cW'}_l\to \cO_{\dC^*\times M}
\quad \textup{ is non-degenerate.}
\end{eqnarray*}
For $l\geq 0$, the pairing
$$P_l:= P((iN)^l.,.):\Gr^{\cW'}_l\otimes j^* \Gr^{\cW'}_l\to \cO_{\dC^*\times M}$$
is well defined, non-degenerate, $(-1)^{w-l}$-symmetric, flat,
and it takes values in $i^{w-l}\dR$ on $(\Gr^{W'}_l)_\dR$.
Let $(\Gr^{W'}_l)_{prim}=\ker N^{l+1}$ be the subbundle of primitive
subspaces.
The decomposition
$$\Gr^{W'}_l =\bigoplus_{j\geq 0} N^j(\Gr^{W'}_{l+2j})_{prim}.$$
is flat. For $l\geq 0$ it is $P_l$-orthogonal.
\item
(Induced TERP-structures)
Suppose that the map $zN:H'\to H'$ extends to a bundle endomorphism of $H$
which has the same Jordan normal form at each point of $\dC\times M$
(notice that as $N$ is flat on $\dC^*\times M$, this is a condition only at the points of $\{0\}\times M$).
Let $W_\bullet$ be the weight filtration (centered at 0) of $zN$ on $H$.

Then $W_l$ is a subbundle of $H$ which extends $W'_l$ to $\omm$,
and $\Gr^W_l$ is a quotient bundle which extends $\Gr^{W'}_l$.
The flat connections on them have poles of type 1 along $\omm$.

The decomposition in 1. extends to a decomposition
$$
\Gr^{W}_l =\bigoplus_{j\geq 0} (zN)^j(\Gr^{W}_{l+2j})_{prim}.
$$
For $l\geq 0$, the summands on the right hand side,
equipped with the pairing $P_l$, are variations of TERP-structures of weight $w-l$.

\item
For $l \geq 0$ and $0\leq j\leq l$ the map $(zN)^j$ extends
to an isomorphism from $\whh{(\Gr^{W}_{l})_{prim}}$ to
${((zN)^j(\Gr^{W}_{l})_{prim})}\mathbf{\!\!\widehat{\quad}}$.
Now fix $t\in M$. If $l-2j\geq0$, then both
$\whh{(\Gr^{W}_{l})_{prim}}(t)$ and $((zN)^j(\Gr^{W}_{l})_{prim}))^\mathbf{\!\!\widehat{\quad}}\!\!(t)$
are TERP-structures.
Using this isomorphism, the two hermitian metrics
on the spaces of the global holomorphic sections
are equal up to the factor $(-1)^j$.

It follows that for any $t\in M$ all $\Gr^W_l(t)$ for $l\in\dZ$
are pure TERP-structures if and only if all $(\Gr^{W}_l)_{prim}(t)$
for $l\geq 0$ are pure TERP-structures.
\item
(PMTS)
The map $N:H'\to H'$
extends to a map
$$\whh N: \whh\cH\to\whh\cH\otimes \cO_{\dP^1}\cC^{an}_M(1,1).$$
Remember the pairing $\whh S$ from lemma \ref{defTopData}.
For any $t\in M$ the tuple $(\whh{H}_{|\dP^1\times\{t\}},\whh S,\whh N)$
is a polarized mixed twistor
structure \cite[definition 3.48]{Mo2} iff
all the $(\Gr^{W}_l)_{prim}(t)$ are pure polarized TERP-structures.
\end{enumerate}
\end{lemma}

\begin{proof}
\begin{enumerate}
\item
The $W'_l$ are flat subbundles with real structure
because $N$ is flat and respects the real structure.
The remaining part is shown as in \cite[Lemma 6.4]{Sch} with the
exception that $P$ is a pairing between different fibers.
\item
The connections on $W_l$ and $\Gr^W_l$ have a pole of type 1 along
$\{0\}\times M$, because the same holds for $H$ and because $W'_l$
is a flat subbundle of $H'$. The decomposition follows again
as in \cite[Lemma 6.4]{Sch}.
It remains to show that $P_l$ maps $\Gr^W_l\otimes j^*\Gr^W_l$ to
$z^{w-l}\cO_{\dC\times M}$ and that it is non-degenerate.
Let $\sigma_1,\ldots,\sigma_\mu$ be a basis of the germ $\cH_{(0,t)}$
for some $(0,t)\in\dC\times M$ which is adapted to the filtration
$W_\bullet$. The matrix $(z^{-w}P(\sigma_i,\sigma_j))$ is holomorphic
and non-degenerate near $(0,t)$, and it has a block lower triangular
shape with respect to the antidiagonal.
If $a,b\in \cW_l$, then $(izN)^la\in \cW_{-l}$,
and the classes $[a],[b]\in \Gr^\cW_l$ satisfy
$$z^{-(w-l)}P_l([a],[b])=z^{-w}P((izN)^la,b)\in \cO_{\dC\times M}.$$
These observations show the properties of $P_l$
needed for a TERP-structure of weight $w-l$ on $\Gr^W_l$.
The decomposition respects real structure and pairing, thus
also the summands are TERP-structures.
\item
Fix $t\in M$. We have to compare the
TERP-structures $(\Gr^W_{l})_{prim}(t)$ of weight $w-l$
and $(zN)^j(\Gr^W_{l-2j})_{prim}(t)$ of weight $w-l+2j$.
The morphisms usually called $\tau$ differ and are called
$\tau_1$ and $\tau_2$ here. To show that $(zN)^j$ extends to an isomorphism
of vector bundles on $\dP^1$, we have to prove
$$(\gamma(z)N)^j\circ\tau_1 = \tau_2\circ (zN)^j.$$
But for any $a\in H_{z}$ where $z\in\dC^*$,
\begin{eqnarray*}
(\gamma(z)N)^j\big(\tau_1(a)(\gamma(z))\big)
&=& (\gamma(z)N)^j\big(\tau_1(a(z))\big)\\
&=& (\gamma(z)N)^j\big(\textup{flat shift of }\overline{z^{-(w-l)}a}\big)\\
&=& \textup{flat shift of }\overline{z^{-(w-l+2j)}(zN)^j(a)}\\
&=& \tau_2((zN)^j(a))(\gamma(z)).
\end{eqnarray*}
A similar calculation shows
$h_2((zN)^ja,(zN)^jb)=(-1)^jh_1(a,b)$ where $h_1$ and $h_2$ denote
the respective hermitian forms on the spaces of global sections
for some fixed $t\in M$. We leave it to the reader.
\item
We fix once and for all $t\in M$ as we do not care here about the
real analytic dependence on the parameters in $M$.
It follows from the definition of $\tau$ and the flatness of $N$ that
$$(\gamma(z)^{-1}N)\circ \tau = \tau \circ (zN).$$
Therefore the conjugate of $zN$ under $\tau$ is $z^{-1}N$.
As $zN$ is a nilpotent endomorphism of $H(t)$
which has everywhere the same Jordan normal form,
the same holds for $z^{-1}N$ as an endomorphism of $\whh H(t)_{|\dP^1-\{0\}}$.

On $H'(t)$ the two endomorphisms coincide up to the scalar $z^2$.
Therefore their weight filtrations coincide on $H'(t)$ and glue to a weight filtration
$\whh W_\bullet(t)$ on $\whh H(t)$.
Furthermore, from the above equation we get that $\whh W_l(t)$
is obtained by gluing $W_l(t)$ with $\overline{\gamma^*W_l(t)}$
via $\tau$. Also, now it is clear that $N:H'(t)\to H'(t)$
extends to a morphism
$\whh N:\whh\cH(t)\to\whh\cH(t)\otimes \cO_{\dP^1}(1,1).$
By definition, the weight filtration associated to this
morphism is $\whh W_\bullet(t)$.

The tuple $(\whh H(t),\whh N)$ is a mixed twistor
iff all $\Gr^{\whh W}_l(t)$ are pure twistors of weight $l$
\cite{Si5}\cite[definition 2.30]{Mo1}.
We first show that this is equivalent to all
$\whh{\Gr^W_l}(t)$ being pure twistors of weight 0.
Both quotients are obtained by gluing $\Gr^W_l$ with
$\overline{\gamma^*\Gr^W_l}$, the first one via $\tau$,
the second one via $\tau_l$ where
$$\tau_l: H_{(z,t)}\to H_{(\gamma(z),t)},\quad
a\mapsto \textup{ flat shift of }\overline{z^{-(w-l)}a}.$$
Comparing $\tau$ and $\tau_l$ we see that
$\tau_l(b)(z) = z^{-l}\tau(b)(z)$. This shows
$$\whh{\Gr^W_l}(t) \cong \Gr^{\whh W}_l(t)\otimes \cO_{\dP^1}(-l),$$
which proves the claim. Using 3. we obtain the statement:
The tuple $(\whh H(t),\whh N)$ is a mixed twistor
iff all $\whh{(\Gr^W_l)_{prim}}(t)$ ($l\geq 0$) are pure TERP-structures.

It remains to compare the polarization conditions.
The one for the pure TERP-structure
$(\Gr^W_l)_{prim}(t)$ reads
$$z^{-(w-l)}P_l(a,\tau_l(a)) >0 \quad \textup{ for }\quad
a\in H^0(\dP^1,\whh{(\Gr^\cW_l)_{prim}}(t))\,\backslash\,\{0\}.$$
The polarization condition for the pure twistor
$(\Gr^{\whh W}_l)_{prim}(t)$ of weight $l$ as part of the
polarized mixed twistor $(\whh\cH(t),\whh S,\whh N)$ is
that it is polarized by $\whh S(\widehat{N}^l.,.)$ \cite[definition 3.48]{Mo2}.
One has to rewrite this condition with \cite[definition 3.35]{Mo2}
as a polarization condition for a pure twistor of weight 0.
We choose the pure twistor
$(\Gr^{\whh \cW}_l)_{prim}(t)\otimes \cO_{\dP^1}(0,-l)$,
where $\cO_{\dP^1}(0,-l)$ is the sheaf of holomorphic functions
on $\dC$ with a zero of order at least $l$ at $\infty$.
Then the condition is
$$
i^{-l}\cdot \whh S(N^l(a),a)>0\quad \textup{ for }\quad
a\in H^0(\dP^1,(\Gr^{\whh \cW}_l)_{prim}(t)\otimes \cO_{\dP^1}(0,-l))\backslash\{0\}.
$$
The factor $i^{-l}$ arises from
definition 3.35 and the first line in formula (3.8)
(here applied to $(p,q)=(0,-l)$) in \cite{Mo2}.

If $a\in (\Gr^W_l)_{prim}(t)$ is glued with $\tau_l(b)$ to
a global section in $\whh{(\Gr^\cW_l)_{prim}}(t)$,
then the formula $\tau_l(b)(z) = z^{-l}\tau(b)(z)$ shows that $a$ is glued
with $z^{-l}\tau(b)$ to a global section in
$((\Gr^{\whh \cW}_l)_{prim}(t)\otimes \cO_{\dP^1}(0,-l).$
The proof of 4. is now finished by the following calculation.
\begin{eqnarray*}
z^{-(w-l)}P_l(a,\tau_l(a)) &=& z^{-w}z^lP((iN)^l(a),(-z)^{-l}\tau(a)) \\
&=& (-i)^lz^{-w}P(N^l(a),\tau(a)) = i^{-l}\whh S(N^l(a),a).
\end{eqnarray*}
\end{enumerate}
\end{proof}

\section{Limit TERP- and twistor structures}
\label{secLimitTwistor}

In this section we consider variations of TERP-structures
on the complement of a normal crossing divisor.
The fundamental result of Mochizuki \cite[theorem 12.22]{Mo2}
yields a limit mixed twistor structure starting from a tame harmonic bundle.
We will show in this section how this can be applied to
variations of pure polarized TERP-structures. We describe in detail how the
limit objects look like in this situation (theorem \ref{theoLimitTwistor}).
Moreover, we give a complementary result on limit TERP-structures
(theorem \ref{theoLimitTERP}). For a regular singular variation on a
punctured disc, this also applies if there is no harmonic
bundle associated to the variation of TERP-structures (proposition
\ref{propLimitTERPregsing}).

We start by fixing some notations and by introducing multi-elementary
sections and $V$-filtrations, which are the basic tools to
construct the limit objects.
Fix $1\leq l\leq n$ and put $\uul=\{1,\ldots,l\}$,
$\uun=\{1,\ldots,n\}$. We consider $X=\Delta^n$ with coordinates
$(r_1,\ldots,r_n)$, the open submanifold
$Y=(\Delta^*)^l\times \Delta^{n-l}$ and the normal crossing divisor
$D=X \backslash Y=\bigcup_{j\in\uul}D_j$ with irreducible components $D_j=\{r_j=0\}$.
Moreover, for $I\subset \uul$, $I\neq\emptyset$, let
$D_I=\bigcap_{j\in I}D_j$ and
$D_I^\circ = D_I \backslash \bigcup_{j\in\uul \backslash I}D_I\cap D_j$.
We will (as in \cite{Sa6} and \cite{Mo2}) denote by $\cX$ the
product $\dC\times X$, and similarly use $\cY, \cD, \cD_i, \cD_I, \cD_I^\circ$
for the corresponding products with $\dC$.
Finally, write $\pi_X$ for the canonical projection $X\to D_\uul$ and $\pi_Y:Y\to D_\uul$
for its restriction to $Y\subset X$.

Suppose that we are given a variation of TERP-structures $\TERP$ on $Y$. We will first discuss
extensions of the flat bundle $H'\in\VB^\nabla_{\dC^*\times Y}$ to $\dC^*\times D$.
In a second step, they will be used to extend $H$ to $\cX$.

We write $M_j$, $j\in\uul$, for the monodromy automorphism corresponding to
a counter-clockwise loop around $\dC^*\times D_j$.
The monodromy around $\{0\}\times Y$ is still denoted by $M_z$.
They all commute. The semisimple and unipotents parts are
denoted by $M_{j,s}$ and $M_{j,u}$, respectively. The nilpotent parts are defined by
$N_j=\log M_{j,u}$. For any $j\in\uul$, define
$$
\begin{array}{rcl}
C_j&:=&\{a_j\in\dC\ |\ e^{2\pi i a_j}\textup{ is an eigenvalue of }M_j\},\\
C_j^{b_j} &:=& C_j\cap (b_j-1,b_j]_\dC\qquad \textup{ for }b_j\in \dC,\\
C&:=&\prod_{j\in \uul}C_j,\qquad C^{\mathbf{b}}:=\prod_{j\in\uul}C_j^{b_j}
\quad \textup{ for }{\mathbf{b}}\in \dC^l.
\end{array}
$$
We have $C_j=C_j^{b_j}+\dZ$ and $C=C^{\mathbf{b}} +\dZ^l$. The existence of the flat pairing $P$ implies
that $C_j=-C_j$, similarly, the fact that $M_z, M_j\in\Aut(H^\infty_\dR)$ gives
$C_j=-\overline{C_j}$. Put ${\mathbf e^{2\pi i \mathbf{a}}}:=(e^{2\pi i a_1},\ldots,e^{2\pi i a_l})$
for ${\mathbf{a}}\in C$ and ${\mathbf{e}}_i:=(\delta_{ji})_{j\in\uul}\in \dC^l$.
Remember the relations
${\mathbf{a}}\leq {\mathbf{b}},
{\mathbf{a}}\lneq {\mathbf{b}},
{\mathbf{a}}<{\mathbf{b}}$ for
${\mathbf{a}},{\mathbf{b}}\in \dC^l$ from the introduction.

Define the flat bundle $H'(\uul)$ on $\dC^*\times D_\uul$
by iterate application of the functor of nearby cycles to the local system $(H')^\nabla$, that is $H'(\uul):=\psi_{r_1}(\psi_{r_2}(\ldots\psi_{r_l}((H')^\nabla)\ldots)$.
Its fibre over a point $(z,r)\in\dC^*\times D_\uul$ can be described concretely as
$$
H'(\uul,z,r):=\{\textup{ multivalued global flat sections in }
H'_{|\{z\}\times \pi_Y^{-1}(r)}\}.
$$
We denote its sheaf of holomorphic sections by $\cH'(\uul)$.
By definition, this bundle comes equipped with a flat connection,
the corresponding monodromy around $\{0\}\times D_\uul$, which is still denoted by $M_z\in\Aut(H^\infty)$,
a flat real subbundle $H'(\uul)_\dR$, a flat pairing $P$
(which takes values in $i^w\dR$ on $H'(\uul)_\dR$)
and with flat bundle automorphisms denoted by $M_j$ for any $j\in\uul$.
Given $\boldsymbol{\lambda}=(\lambda_1,\ldots,\lambda_l)\in(\dC^*)^l$,
we write $H'(\uul)_{\boldsymbol{\lambda}}$
for the simultaneous eigenspace with eigenvalues $(\lambda_1,\ldots,\lambda_l$)
of the semisimple parts of $(M_1,\ldots,M_l)$.

Notice that one may perform the same construction for any subset
$I$ of $\uul$, yielding a flat bundle on $\dC^*\times D^\circ_I$,
which is also equipped with a real structure and a pairing as above.

Our aim is to obtain an extension of $H'(\uul) \in\VB^\nabla_{\dC^*\times D_\uul}$
to a vector bundle on $\cD_\uul$ starting with a variation of TERP-structures which satisfies
a regularity condition near the divisor $\cD$. For that purpose,
we will need the following generalization of elementary section.

For ${\mathbf{a}}\in C$ and any holomorphic
section $A\in \cH'(\uul)_\mathbf{e^{2\pi i \mathbf{a}}}(U_1\times U_2)$
with $U_1\subset\dC^*$ and $U_2\subset D_\uul$ open, the section
$$es_\uul(A,{\mathbf{a}}):= \prod_{j\in\uul}r_j^{-a_j-\frac{N_j}{2\pi i}}A$$
is a holomorphic section in $H'$ on $U_1\times \pi_Y^{-1}(U_2)$.
Notice that contrary to the elementary sections considered
in section \ref{secBasics} we use the opposite indices here
($-a_j$ instead of $a_j$) and moreover, the section $A$ itself
is not necessarily flat.

The following identities, which will be used frequently in the sequel,
are satisfied by the $\uul$-elementary sections.
\begin{eqnarray}\label{eqElementary-r-1}
z\nabla_z es_\uul(A,{\mathbf{a}})
&=& es_\uul(z\nabla_z A,{\mathbf{a}}),\\ \label{eqElementary-r-2}
\nabla_{r_j} es_\uul(A,{\mathbf{a}}) &=& es_\uul(\nabla_{r_j} A,{\mathbf{a}})
\quad\textup{ for }j\in \uun\,\backslash\,\uul,    \\ \label{eqElementary-r-3}
(r_j\nabla_{r_j}+a_j) es_\uul(A,{\mathbf{a}})
&=& es_\uul(\frac{-N_j}{2\pi i} A,{\mathbf{a}})
\quad\textup{ for }j\in \uul, \\ \label{eqElementary-r-4}
P(es_\uul(A,{\mathbf{a}}) ,es_\uul(B,{\mathbf{b}})) &=& \left\{\begin{matrix}
0 & \textup{ if } & a+b\notin\dZ^l,\\
\prod_{j\in\uul} r_j^{-a_j-b_j}P(A,B) & \textup{ if }
& a+b\in \dZ^l\end{matrix}\right.
\end{eqnarray}
The last equation follows from the flatness of $P$,
which implies in particular that the morphisms
$N_j$ are infinitesimal isometries of $P$.

We obtain an induced increasing V-filtration by extensions
${_{\mathbf{b}}}V$ of $H'$ to vector bundles
on $\dC^*\times X$. The associated locally free sheaf is defined as
\begin{equation}
\label{eqDefDeligneExtension}
{_{\mathbf{b}}}\cV := \sum_{\mathbf{a}\leq \mathbf{b}} \cO_{\dC^*\times X} es_\uul(A,\mathbf{a})
\end{equation}
The reason for considering an increasing instead of a decreasing
V-filtration is that this notation is compatible with the one
used by Mochizuki (see section \ref{secLimitProofs}) for the parabolic filtration
defined for a harmonic bundle.
It follows from the formulas \eqref{eqElementary-r-1}, \eqref{eqElementary-r-2} and \eqref{eqElementary-r-3} that
the sheaves ${_{\mathbf{b}}}\cV$ are invariant with respect to
$\nabla_z, \nabla_{r_j}$ ($j\in \uun\,\backslash\,\uul$)
and $r_j\nabla_{r_j}$ ($j\in\uul$). The residue of $r_j\nabla_{r_j}$
($j\in \uul$) on ${_{\mathbf{b}}}V_{|\dC^*\times D_j}$ has eigenvalues in $C_j^{b_j}$.

Each ${_{\mathbf{a}}}\cV$ carries a filtration by subsheaves
${_\mathbf{b}}\cV$ indexed by $\{\mathbf{b}\in C\,|\,\mathbf{b}\leq \mathbf{a}\}$,
and we have an isomorphism
\begin{equation}
\label{eqIsoPhi}
\begin{array}{rcl}
\Phi':\cH'(\uul)_{e^{2\pi i \mathbf{a}}} & \longrightarrow & \Gr_\mathbf{a}({_\mathbf{a}\cV}) \\ \\
A & \longmapsto & [es(A,\mathbf{a})],
\end{array}
\end{equation}
in particular, the quotient sheaves $\Gr_\mathbf{a}({_\mathbf{a}\cV})$ are locally free on $\cO_{\dC^*\times D_\uul}$.
The system of locally free sheaves $({_\mathbf{a}}\cV)_{\mathbf{a}\in C}$
is a locally abelian parabolic bundle in the sense of definition \ref{defParabolicBundle},
this will be shown in lemma \ref{lemVFLocAbelian}, 1.

The $\uul$-elementary sections can be used
to describe general sections of the bundle $H'$.
More precisely, suppose that $U_1\subset \dC^*$,
$U_2\subset D_\uul$ and $U\subset \dC^*\times X$
are open subsets such that $U_1\times U_2\subset U$.
A section $\sigma \in \cH'(U\cap Y)$ is an in general infinite
sum of $\uul$-elementary sections on $U_1\times \pi_Y^{-1}(U_2)$, namely
\begin{eqnarray*}
\sigma = \sum_{{\mathbf{a}}\in C}es_\uul(A(\sigma,{\mathbf{a}}),{\mathbf{a}}),
\end{eqnarray*}
where $A(\sigma,{\mathbf{a}})$ are uniquely determined sections in
$\cH'_{\mathbf{e}^{2\pi i \mathbf{a}}}(\uul)(U_1\times U_2)$.
These pieces $A(\sigma,{\mathbf{a}})$ satisfy the following equations, which
will also be quite useful later.
\begin{eqnarray}\label{eqPiecesElementary-1}
z^2\nabla_z A(\sigma,{\mathbf{a}})
&=& A(z^2\nabla_z(\sigma),{\mathbf{a}}),\\ \label{eqPiecesElementary-2}
z\nabla_{r_j}A(\sigma,{\mathbf{a}}) &=& A(z\nabla_{r_j}(\sigma),{\mathbf{a}})
\quad\textup{ for }j\in \uun\,\backslash\,\uul,    \\ \label{eqPiecesElementary-3}
z(-a_j-\frac{N_j}{2\pi i}) A(\sigma,{\mathbf{a}})
&=& A(zr_j\nabla_{r_j}(\sigma),{\mathbf{a}})
\quad\textup{ for }j\in \uul.
\end{eqnarray}

\begin{remark}\label{remFlatBundleOnI}
As we noticed above, it is possible to define a flat bundle $H'(I)\in\VB^\nabla_{\dC^*\times D^\circ_I}$
for any nonempty subset $I\subset\uul$. In a similar way, one can define $I$-elementary sections.
Composition of these operations behaves well, more precisely,
for any $I,J\subset\uul$ such that $I\cap J=\emptyset$, we have
a canonical isomorphism $H'(H'(I),J) \cong H'(I\cup J)$
and similarly $es_J(es_I(A,{\mathbf{a}}),{\mathbf{a}})
=es_{I\cup J}(A,{\mathbf{a}})$.
\end{remark}

Up to this point, the only input data we used was the flat bundle $H'$
with its real structure and the pairing $P$. If we are given a variation of TERP-structures
$\TERP$ on $Y$, the $\uul$-elementary sections and the
increasing $V$-filtration can be used to control
the behavior of $\cH$ near the divisor $\cD$ and to
discuss possible extensions to $\cX$. More precisely, consider
the inclusions
$j^{(1)}:\dC^*\times X\hookrightarrow \cX$
and $j^{(2)}:\cY\hookrightarrow \cX$.
We define for any ${\mathbf{a}}\in C$ the sheaf
\begin{eqnarray}\label{eqDefSheavesF}
{_{\mathbf{a}}}\cF := j^{(1)}_*{_{\mathbf{a}}}\cV\cap j^{(2)}_*\cH
\end{eqnarray}
on $\cX$. It is by definition locally free on
$(\dC^*\times X) \cup \cY = \cX \backslash (\{0\}\times D)$,
but it does not even need to be coherent on the codimension two subset
$\{0\}\times D$.
\begin{definition}\label{defTERPTame}
Let $(X,Y,D)$ be as above.
A variation of TERP-structures $\TERP$ on $Y$ is called \emph{tame} (along $D$)
iff all sheaves ${_{\mathbf{a}}}\cF$ are locally free.
\end{definition}

Theorem \ref{theoLimitTERP} treats tame variations of TERP-structures.
We start with two examples which are not tame.
\begin{examples}\label{exTERP-NonTame-1}
Any of the variations of TERP-structures from example \ref{exTERP-1} is denoted
by TERP. We do not care about the real structure or the pairing here,
so any of the choices made for them in example \ref{exTERP-1} can be used here.

\begin{enumerate}
\item
Consider $X=\dC$, $Y=\dC^*$ and
$\varphi_1:Y\to \dC^*,\quad r\mapsto e^{1/r}$.
The pullback $\varphi^*_1(TERP)$ is a variation of TERP-structures
on $Y$ which is generated by the sections
$$z^{-1}A_1+e^{1/r}A_2,\quad zA_2, \quad A_1.$$
The sheaf ${_\mathbf{0}}\cF$ is locally free on $\dC\times X \backslash \{0\}$,
and $zA_2$ and $A_1$ are global sections whereas $z^{-1}A_1+e^{1/r}A_2$ is not.
This implies that ${_\mathbf{0}}\cF$ is not
coherent at $0$. The reason for this is that $zA_2$, $A_1$ do not generate
the restriction ${_\mathbf{0}}\cF$ to $\cY$ (i.e., $\cH$ itself),
as they should by the implication i) --> iii) in \cite[th\'eor\`eme 1]{Se} if ${_\mathbf{0}}\cF$
were coherent.

\item
Consider $X=\dC^2$, $Y=(\dC^*)^2$ and
$\varphi_2:X \backslash \{0\}\to \dP^1,\quad (r_1,r_2)\mapsto (r_1:r_2)$.
The pullback $\varphi^*_2(TERP)$ is a variation of TERP-structures
on $X \backslash \{0\}$ which is generated by
$$
r_2z^{-1}A_1+r_1A_2,\quad zA_2,\quad A_1,
$$
and it can be restricted to $Y$.
The sheaf ${_\mathbf{0}}\cF$ is locally free on $\dC\times X \backslash \{0\}$,
but at $0$ it is not, but only coherent with three
generators. For any fixed $(r_1:r_2)$ the restriction
of the variation on $X \backslash \{0\}$ to $\varphi_2^{-1}((r_1:r_2))$
is a constant variation. As they are all different, their limits
for $(r_1,r_2)\to 0$ are not compatible.
\end{enumerate}
\end{examples}

In the tame case, the various ingredients of the TERP-structure can be extended
to the sheaves $_\mathbf{a}\cF$. This is done in the following lemma.
\begin{lemma}
\label{lemExtensionSheaves_aF}
Let $\TERP$ be a variation of TERP-structures, tame along $D$. Then for any $\mathbf{a}\in C$, the
connection extends as
$$
\nabla: {_\mathbf{a}}\cF \longrightarrow {_\mathbf{a}}\cF \otimes z^{-1}\Omega^1_\cX\left(\log (\cD\cup(\{0\}\times X))\right).
$$
Moreover, $P$ extends to a non-degenerate pairing
$$
P: {_\mathbf{a}}\cF \otimes j^* {_\mathbf{b}}\cF \longrightarrow z^w \cO_\cX,
$$
where ${\mathbf{b}}\in C$ is the unique multi-index satisfying  $-C^{\mathbf{a}} = C^{\mathbf{b}}$,
i.e. $b_j=\max C_j\cap(-\infty,a_j+1)$.
\end{lemma}
\begin{proof}
As we have seen, the sheaves
${_\mathbf{a}}\cV$ are invariant under the connection operators $\nabla_z$, $\nabla_{r_j}$
($j\in \uun\,\backslash\,\uul$) and $r_j\nabla_{r_j}$ ($j\in \uul$). Moreover, $\cH$ is invariant
under $z^2\nabla_z$ and $z\nabla_{r_j}$ ($j\in\uun$) as it underlies a variation of TERP-structure.
It follows that ${_\mathbf{a}}\cF$ is invariant under $z^2\nabla_z$, $z\nabla_{r_j}$
($j\in \uun\,\backslash\,\uul$) and $zr_j\nabla_{r_j}$ ($j\in \uul$). This proves the first statement.
Concerning the second one, notice that formula \eqref{eqElementary-r-4} yields that
$P:{_\mathbf{a}}\cV\otimes j^*{_\mathbf{b}}\cV\rightarrow \cO_{\dC^*\times X}$ is non-degenerate.
Now choose locally on $\{0\}\times D_\uul$ arbitrary bases of ${_\mathbf{a}}\cF$ and
${_\mathbf{b}}\cF$, then the corresponding matrix of $z^{-w}P$ is holomorphic and invertible on $\dC^*\times X$ and
on $\cY$, i.e., outside of the codimension two subset $\{0\}\times D$. Therefore it is holomorphic
and invertible all over $\cX$, and the pairing $P:{_\mathbf{a}}\cF\otimes j^* {_\mathbf{b}}\cF\rightarrow z^w\cO_\cX$ is non-degenerate.
\end{proof}
For any $\mathbf{a}\in C$, the increasing filtration of ${_\mathbf{a}}\cV$ considered
above induces by definition an increasing filtration of ${_\mathbf{a}}\cF$, given by
the subsheaves ${_\mathbf{b}}\cF$ for any $\mathbf{b}\in C$ with $\mathbf{b}\leq \mathbf{a}$.
However, it is considerably less obvious that the corresponding quotients are locally free.
This is part of the next theorem, which is the first main result of this section.
\begin{theorem}\label{theoLimitTERP}
Let $X,Y,D$ be as above and let $\TERP$ be a tame variation on $Y$.
\begin{enumerate}
\item
For any ${\mathbf{a}}\in C$, the quotient
sheaf
$$
\Gr_\mathbf{a}({_\mathbf{a}}\cF):=\frac{{_\mathbf{a}}\cF}{\sum_{\mathbf{b}\lneq\mathbf{a}}{_\mathbf{b}}\cF}
$$
is locally free over $\cO_{\cD_\uul}$ and defines an extension of $Gr_\mathbf{a}({_\mathbf{a}}V)$ to a vector bundle
on $\cD_\uul$. Via the inverse of the isomorphism $\Phi'$ in formula \eqref{eqIsoPhi}, this
induces an extension of $H'(\uul)_{e^{2\pi i a}}$ to $\cD_\uul$.
This extension is independent of the choice of $\mathbf{a}$ within the set
$\mathbf{a}+\dZ^l$. It is denoted by $H(\underline{l})_{e^{2\pi i \mathbf{a}}}$. We have the following
isomorphism of the associated sheaves
\begin{equation}
\label{eqDescriptionExtensionH}
\begin{array}{rcl}
\Phi^{-1}:\Gr_\mathbf{a}({_\mathbf{a}}\cF) & \stackrel{\cong}{\longrightarrow} & \cH(\underline{l})_{e^{2\pi i \mathbf{a}}} = \sum_{\sigma\in{_\mathbf{a}}\cF} \cO_{\cD_\uul} \cdot A(\sigma,\mathbf{a}) \\ \\
{[}\sigma= \sum_{\mathbf{b}\leq \mathbf{a}} es(A(\sigma,\mathbf{b}),\mathbf{b}){]}
& \longmapsto &  A(\sigma, \mathbf{a}).
\end{array}
\end{equation}
We put $H(\uul):=\sum_{\mathbf{a}\in C} H(\uul)_{e^{2\pi i \mathbf{a}}}$, then the tuple
$(H(\uul),H'(\uul)_\dR,\nabla,P)$
is a variation of TERP-structures of weight $w$ on $D_\uul$.
Furthermore, the nilpotent endomorphisms $zN_j,j\in\uul$, on $H'(\uul)$
extend to nilpotent endomorphisms on the bundle $H(\uul)$. However,
it is unclear whether they have the same Jordan normal form at each point
in $\{0\}\times D_\uul$.

\item For any $I\subset \uul$, the same construction yields
an extension of $H'(I)_{e^{2\pi i \mathbf{a}}}$ to a vector bundle
$H(I)_{e^{2\pi i \mathbf{a}}}$ on $\cD_I^\circ$, and the sum
$H(I)$ underlies a variation of TERP-structures on $D_I^\circ$,
tame along $D_I\backslash D_I^\circ$. Moreover, for any
$J\subset \uul\backslash I$, we have $H(I\cup J)\cong H(H(I),J)$.
\end{enumerate}
\end{theorem}
The proof of this theorem will be postponed until section \ref{secLimitProofs}.
It relies on a general result concerning parabolic bundles
on $\cX$, which is given in section \ref{secCompatible}. This result applies to the
system of locally free sheaves $({_{\mathbf{a}}}\cF)_{{\mathbf{a}}\in C}$.
More precisely, we construct in theorem \ref{theoCompatibleLocallyAbelian} a compatible system of
local bases for all ${_{\mathbf{a}}}\cF,{\mathbf{a}}\in C$,
which yield the proof of theorem \ref{theoLimitTERP} essentially by using the formulas
\eqref{eqElementary-r-4} - \eqref{eqPiecesElementary-3}.

Let us turn back to the examples \ref{exTERP-NonTame-1}. These were seen to be non-tame
variations of TERP-structures, as a consequence, the limit construction of theorem
\ref{theoLimitTERP} does not work here.
\begin{examples}\label{exTERP-NonTame-2}
In example 1., $({_\mathbf{0}}\cF)_0$ at $0$ is a free $\cO_{\dC\times X,0}$ module
of rank 2, generated by $zA_2$ and $A_1$ (which do not generate ${_\mathbf{0}}\cF$ in
a neighborhood of $0$). This implies that
$A((zA_2,\mathbf{0}),\mathbf{0})=zA_2$ and
$A((A_1,\mathbf{0}),\mathbf{0})=A_1$.
As we already remarked, in example 2., the germ $({_\mathbf{0}}\cF)_0$ is not free but
generated by the 3 sections
$$
r_2z^{-1}A_1+r_1A_2,\quad zA_2,\quad A_1.
$$
We have
$A(zA_1,\boldsymbol{0})=zA_2$, $A(A_1,\boldsymbol{0})=A_1$ but
$A(r_2z^{-1}A_1+r_1A_2,\boldsymbol{0})=0$.
In both cases the construction in theorem
\ref{theoLimitTERP} gives an extension of $H'(\uul)$ to the vector bundle
generated by $zA_2$ and $A_1$.
The connection has even a logarithmic pole, and the pairing extends
holomorphically, but it is degenerate at 0. Therefore this extension
is not a TERP-structure.

Notice that the above mentioned compatibility condition is not satisfied
in example 2. More precisely, putting
$I=\{1\}$ and $J=\{2\}$ we obtain variations of limit TERP-structures
$H(I)$ on $D^\circ_{I}$ and $H(J)$ on $D^\circ_{J}$. Both are constant variations,
and they are different, so their limits on $D_{I\cup J}=D_\uul=\{0\}$,
are non-isomorphic.
\end{examples}

The second result of this section gives a much stronger result about the limit object
$H(\uul)$ under the additional hypothesis that the variation
we started with is pure polarized. It builds on
\cite[theorem 12.22]{Mo2}, which describes
a limit polarized mixed twistor structure defined by a tame harmonic bundle on $Y$.
Recall that a variation of pure polarized TERP-structures $\TERP$
on a manifold $M$ gives rise to a harmonic bundle, namely
$(H_{|\{0\}\times M},\overline{\partial},\theta,h)$ \cite[chapter 2]{He4}.
Here the operator $\overline{\partial}$ is the one defining the holomorphic structure on $H_{|\{0\}\times M}$
whereas the Higgs field $\theta$ is the pole part along $\{0\}\times M$ of the connection $\nabla$ with respect to
vector fields on $M$. The hermitian metric $h$ is obtained from $pr_*\whh\cH$ by the real analytic isomorphism
$\cH^{an}_{|\{0\}\times M}\to pr_*\whh\cH$, which exists as $H$ is pure.

For any harmonic bundle $(E,\overline\partial,\theta,h)$ on $Y=(\Delta^*)^l\times \Delta^{n-l}$ as above,
the Higgs field can be written in the coordinates $(r_1,\ldots,r_n)$ as follows:
$$\theta = \sum_{j\in \uul}\theta_j\frac{\textup{d} r_j}{r_j}
+ \sum_{j\in \uun \backslash \uul} \theta_j\textup{d}r_j.$$
$(E,\overline\partial,\theta,h)$ is called \emph{tame} along $D=X \backslash Y$
if the coefficients of the characteristic polynomials of all endomorphisms $\theta_j$
extend to holomorphic functions on $X$.
\begin{theorem}\label{theoLimitTwistor}
Let $X,Y,D$ be as above and let $\TERP$ be a variation of
pure polarized TERP-structures such that the associated
harmonic bundle is tame along $D$. Then the following holds.
\begin{enumerate}
\item
The variation of TERP-structures is tame, so that theorem \ref{theoLimitTERP} applies.
We obtain a (limit) variation of TERP-structures
$(H(\uul),H'(\uul)_\dR,\nabla,P)$ on $D_\uul$.
\item
For each ${\mathbf{a}}\in (\dR^+)^l$, the nilpotent endomorphism
$zN_{\mathbf{a}}=\sum_{j\in\uul}za_jN_j$ of the bundle $H(\uul)$
has at each point of $\cD_\uul$ the same Jordan normal form.
Therefore it induces a weight filtration $W_\bullet$ on $H(\uul)$ by subbundles.
This weight filtration does not depend on the choice of ${\mathbf{a}}\in (\dR^+)^l$.
\item
For any $r\in D_\uul$ and any ${\mathbf{a}}\in (\dR^+)^l$,
the limit TERP-structure at $r$ yields a polarized mixed twistor structure
$(\whh{\cH(\uul)}(r),\whh S,\whh {N_{\mathbf{a}}})$ in the sense of
lemma \ref{lemTERPWeightFilt}.
\item
The quotients $\Gr^W_l$ as well as the summands in the
decomposition $\Gr^{W}_l =\bigoplus_{j\geq 0} (zN)^j(\Gr^{W}_{l+2j})_{prim}$ are variations
of pure TERP-structures of weight $w-l$. Moreover, any
$(\Gr^W_l)_{prim}$ is a
variation of pure polarized TERP-structures.
\end{enumerate}
\end{theorem}
The proof of this theorem, which is essentially an
application of \cite[theorem 12.22]{Mo2}, will also be given in section \ref{secLimitProofs}.

The next result shows that theorem \ref{theoLimitTwistor} applies
in the case of a variation of regular singular pure polarized
TERP-structures.

\begin{proposition}\label{propRegSingTame}
If $\TERP$ is a variation of regular singular TERP-structures
on a manifold $M$ then all endomorphisms of the Higgs field on
$H_{|\{0\}\times M}$ are nilpotent. Therefore,
if the TERP-structures are also pure and polarized, then the associated
harmonic bundle is tame along any divisor.
\end{proposition}

\begin{proof}
The endomorphisms of the Higgs field are the endomorphisms
$[z\nabla_X]$ on $\cH/z\cH\to \cH/z\cH$, $X\in \cT_M$.
Consider the Deligne extensions $\cV^\alpha$ of $\cH'$ to $\dC\times M$.
Any $\cV^\alpha$ is stable under $\nabla_X$ by definition.
This implies
$$
[z\nabla_X]: \cV^\alpha(\cH/z\cH)\to
\cV^{\alpha+1}(\cH/z\cH).
$$
Because of this and lemma \ref{lemRegSingBounded}, $[z\nabla_X]$ is nilpotent.
The tameness is now obvious,
as the only eigenvalue of $[z\nabla_X]$ is zero.
\end{proof}

The following proposition treats the case of a regular singular
variation on $\Delta^*$ with an a priori much weaker tameness
assumption than that in definition \ref{defTERPTame}.
The proof builds on \cite{Se}. However, we do not obtain
a polarized mixed twistor in the limit.

\begin{proposition}\label{propLimitTERPregsing}
Let $\TERP$ be a variation of regular singular TERP-structures
on $Y=\Delta^*$. Suppose that at least one ${_{\mathbf{a}}}\cF$ is coherent
or that its (locally free) restriction
${_{\mathbf{a}}}\cF_{|\dC\times\Delta \backslash \{0\}}$
is generated by global sections. Then all ${_{\mathbf{b}}}\cF$ are
locally free, so theorem \ref{theoLimitTERP} applies and gives
a limit TERP-structure $\cH(\underline{1})$. Moreover, if $\cH\subset\cV^\alpha$,
then $\cH(\underline{1})\subset V^\alpha$. In particular,
in the situation of theorem \ref{theoLimitTERP}, 2., the spectrum of
the limit TERP-structure $H(I)_{|\dC\times\{x\}}$ is contained
in the interval $[\alpha,w-\alpha]_\dC$ for any $x\in D^\circ_I$.
\end{proposition}

\begin{proof}
By construction ${_{\mathbf{a}}}\cF= j^{(3)}_*{_{\mathbf{a}}}\www\cF$ where
$j^{(3)}:(\dC\times \Delta \backslash \{0\})\hookrightarrow\dC\times \Delta$
and ${_{\mathbf{a}}}\www\cF:= {_{\mathbf{a}}}\cF_{|\dC\times \Delta-\{0\}}$ is locally free.
By \cite[theorem 1]{Se}, ${_{\mathbf{a}}}\cF$ is coherent iff
there is a neighborhood $U\subset\dC\times\Delta$ of $0$ such that
at each point in $U \backslash \{0\}$ the sheaf
${_{\mathbf{a}}}\www\cF$ is generated by its global sections in $U \backslash \{0\}$.
Therefore the second assumption from above is equivalent to
the coherence of ${_{\mathbf{a}}}\cF$. As ${_{\mathbf{a}}}\www\cF$
is locally free, ${_\mathbf{a}}\cF$ is reflexive \cite[proposition 7]{Se}. As the base has
dimension 2, it is locally free.

Now consider any other ${_{\mathbf{b}}}\cF$.
We will apply \cite[theorem 1]{Se}
to show that it is coherent. Then reflexiveness and local freeness
follow as above. By lemma \ref{lemRegSingBounded} there exists $\beta$ with $\cH\supset \cV^\beta$.
The sheaf $j^{(1)}_*{_{\mathbf{b}}}\cV\cap j^{(2)}_*\cV^{\beta}$ is locally free,
as it is generated by sections which are elementary with respect to
$r$ and $z$. The sheaf
$r^{[{\mathbf{b}}-{\mathbf{a}}+1]}{_{\mathbf{a}}}\cF$ is locally free,
because ${_\mathbf{a}}\cF$ is locally free.
The union of bases of both sheaves
generates ${_{\mathbf{b}}}\www\cF$ at each point
in $U \backslash \{0\}$. Using \cite[theorem 1]{Se} again,
we obtain the coherence of ${_{\mathbf{b}}} \cF$.

As to the last statement, consider any local section $s\in\cH(\underline{1})$ and
decompose it into a sum of $z$-elementary sections.
If in this decomposition there is any $z$-elementary section with order $\beta<\alpha$,
then it necessarily also appears in some section of some $_a\cF$.
This implies that $_a\cF \not\subset j^{(1)}_*{_a}\cV\cap j^{(2)}_*\cV^\alpha$
from which we conclude that $\cH\not\subset\cV^\alpha$, which contradicts the assumption.
\end{proof}

\section{Locally abelian parabolic bundles}
\label{secCompatible}

In this section we consider parabolic bundles on an arbitrary
complex manifold $M$. The main result is
theorem \ref{theoCompatibleLocallyAbelian}, which says that
any parabolic bundle is locally abelian in the sense of \cite{SiIy1,SiIy2}.
%This condition can also be expressed as the compatibility of associated filtrations
%in the sense of \cite[chapter 4]{Mo2}. We will comment on this at the end of this section.

This theorem was proved first
by Borne %\cite[remarque 3 part 3.]{Borne2}.
\cite[th\'eor\`eme 2.4.20]{Borne2}.
However, his proof is adapted to a more general situation, it is done in the algebraic category
and moreover it is spread over the two papers \cite{Borne1} and \cite{Borne2}.
Therefore we found it useful to offer here a short proof, which is actually
a mixture of an (independent) proof we had in the first version
of this paper and of Borne's proof. We comment on the relation between the proofs
at the end of this section. We thank the referee for pointing us to Borne's work.

Theorem \ref{theoCompatibleLocallyAbelian} is applied in
the next section in the proof of theorem \ref{theoLimitTERP}. More precisely,
it shows that for a variation of TERP-structures on $Y=X\backslash D$,
tame along $D$, the system of locally free sheaves $({_\mathbf{a}}\cF)_{\mathbf{a}\in C}$
as defined by formula \eqref{eqDefSheavesF} is a locally abelian parabolic bundle on $\cX$.

We start by recalling briefly the notion of a parabolic sheaf, in order to fix the
notations. We follow \cite{SiIy1,SiIy2}, however, we consider the corresponding analytic objects,
and we also allow arbitrary complex numbers as weights of the
parabolic structure. This imposes a slight change in the definition
compared to loc.cit., on which we comment later.

\begin{definition}\label{defParabolicBundle}
Let $M$ be a complex manifold, $l\in\dN$ and $D=\coprod_{i\in\uul}D_i \subset M$ be a normal crossing
divisor with irreducible components $D_i$. Write, as before, $D_\uul$
for the intersection $\bigcap_{i\in\uul}D_i$.
\begin{enumerate}
\item
A parabolic sheaf on $(M,D)$ is a family
$({_\mathbf{a}}\cE)_{\mathbf{a}\in\dC^l}$ of
torsion free $\cO_M$-modules such that the following holds.
\begin{enumerate}
\item
For any $\mathbf{a}, \mathbf{b}\in \dC^l$ with $\mathbf{a}\leq\mathbf{b}$, ${_\mathbf{a}}\cE$ is
an $\cO_M$-submodule of $_\mathbf{b}\cE$.
\item(support condition)
For any $\mathbf{a}\in \dC^l$, we have $_{\mathbf{a}-\mathbf{e}_i}\cE= {_\mathbf{a}}\cE(-D_i)$.
\item
There exists a discrete index set
$C=\prod_{i=1}^l C_i \subset \dC^l$ such that $C_i^{a_i}:=C_i\cap (a_i-1,a_i]_\dC$ is finite
for any $a_i\in\dC$, and $C_i=C_i^{a_i}+\dZ$ such that the system of sheaves $({_\mathbf{a}}\cE)_{\mathbf{a}\in C}$
determines the parabolic sheaf $({_\mathbf{a}}\cE)_{\mathbf{a}\in \dC^l}$ in the following way:
For all $\mathbf{a}\in \dC^l$, we have ${_\mathbf{a}}\cE={_{\widetilde{\mathbf{a}}}}\cE$,
where $\widetilde{\mathbf{a}}_i := \max\left((-\infty,a_i]_\dC \cap C_i\right)$.
We will often use the notation $C^{\mathbf{a}}=\prod_{i=1}^lC_i^{a_i}$.
\end{enumerate}
\item
A parabolic sheaf $\cE$ on $(M,D)$ is called a parabolic bundle if all ${_\mathbf{a}}\cE$
are locally free.
\item
For any $a\in\dC^l$, denote by ${^\mathbf{a}}\cL:=({^\mathbf{a}_\mathbf{b}}\cL)_{\mathbf{b}\in\dC^l}$ the parabolic line bundle
on $(M,D)$ defined by ${^\mathbf{a}_\mathbf{b}}\cL:=\cO_M(\sum_{i=1}^l \lfloor a_i+b_i \rfloor D_i)$.
(Notice that $^\mathbf{a}\cL$ is called $\cO_M(\sum_{i=1}^l a_i D_i)$ in \cite{SiIy1}).
\item
A parabolic bundle $({_\mathbf{a}}\cE)_{\mathbf{a}\in \dC^l}$ is called locally abelian, if it is locally isomorphic (as a parabolic bundle)
to $\oplus_{i=1}^m({^{\mathbf{a}^i}}\cL)$ for suitable $\mathbf{a}^i\in\dC^l$.
\end{enumerate}
\end{definition}

It follows from the definition of a parabolic sheaf that for any $\mathbf{b}\in\dC^l$ and any
$\mathbf{a}\in(\mathbf{b}-\de,\mathbf{b}]_\dC$, the quotient ${_\mathbf{b}}\cE/{_\mathbf{a}}\cE$ is supported
on the divisor $D$. If $b_i=a_i$ for all
$i\in\uul\,\backslash\,\{j\}$, then
$\textup{supp}({_\mathbf{b}}\cE/{_\mathbf{a}}\cE)\subset D_j$. In
particular, the quotient
$$
\Gr_\mathbf{a}({_\mathbf{a}}\cE):=\frac{{_\mathbf{a}}\cE}{\sum_{\mathbf{c}\lneq
  \mathbf{a}}{_\mathbf{c}}\cE}
$$
is supported on the intersection $D_\uul$.
Notice that this quotient is zero if $\mathbf{a}\notin C$. Moreover,
condition 1.(c) from the above definition implies the weaker \emph{semi-continuity condition}:
For any $\mathbf{a}\in \dC^l$, there is $\varepsilon\in\dR_+$ such that for any $\mathbf{c}\in[0,\varepsilon)_\dC^l$ we have ${_{\mathbf{a}+\mathbf{c}}}\cE={_\mathbf{a}}\cE$.
If we consider a parabolic sheaf indexed by $\dR^l$, and suppose that
for any subset $I\subset \uul$, the intersection $D_I=\cap_{i\in I} D_i$ has only
finitely many components, then the semi-continuity condition and the support condition
actually imply the existence of an index set $C$ with the above properties.

\begin{theorem}\label{theoCompatibleLocallyAbelian}
(\cite[th\'eor\`eme 2.4.20]{Borne2})
Any parabolic bundle is locally abelian.
\end{theorem}

Theorem \ref{theoCompatibleLocallyAbelian} will be proved by induction over $l$,
after lemma \ref{t4.4}.
The first step, $l=1$, will be an immediate consequence of lemma \ref{t4.4}.
Lemma \ref{t4.3} rewrites the condition "locally abelian" in a more explicit way
and draws two useful conclusions.

\begin{lemma}\label{t4.3}
Let $M$ be a complex manifold and $D$ a normal crossing divisor as above.
Let $(_\mathbf{a}\cE)_{\mathbf{a}\in C}$
be a parabolic bundle on $M$.
Then the following two conditions are equivalent.
\begin{enumerate}
\item
$(_\mathbf{a}\cE)_{\mathbf{a}\in C}$ is locally abelian.
\item
For any $t\in D_\uul$ there are local coordinates $r=(r_1,...,r_n)$ on $M$ and a
neighborhood $U$ with $r:(U,t)\to (\Delta^n,0)$ an isomorphism such that for $i\in\uul$,
$D_i\cap U=\{r_i=0\}$, and there are sections
$\sigma_{j,\mathbf{a}}\in {_\mathbf{a}}\cE_{|U}$, $\mathbf{a}\in C, j\in\{1,\ldots,d(\mathbf{a})\}$,
such that the following two conditions hold.
\begin{enumerate}
\item
For any $i\in\uul$ and $j\in\{1,\ldots,d(\mathbf{c})\}$, we have $r_i\cdot \sigma_{j,\mathbf{c}}
=\sigma_{j,\mathbf{c}-\mathbf{e}_i}$.
\item
Denote by $i_{\mathbf{c}\mathbf{a}}$ the inclusion $_\mathbf{c}\cE\subset{_\mathbf{a}}\cE$ for
any $\mathbf{c}\leq\mathbf{a}$. Then
$(i_{\mathbf{c}\mathbf{a}}(\sigma_{j,\mathbf{c}}))_{c\in C^\mathbf{a},j\in\{1,\ldots,d(\mathbf{c})\}}$
is a local basis of ${_\mathbf{a}}\cE_{|U}$.
\end{enumerate}
\end{enumerate}
Suppose that 1. and 2. hold. Then $\Gr_\mathbf{a}({_\mathbf{a}}\cE)$ is a
locally free $\cO_{D_\uul}$-module of finite rank. If $U$ and $\sigma_{j,\mathbf{c}}$
are as in 2., then
\begin{eqnarray}\label{4.1}
\Gr_\mathbf{a}({_\mathbf{a}}\cE_{|U})=\bigoplus_{j\in\{1,\ldots,d(\mathbf{a})\}}
\cO_{D_\uul\cap U} \; [\sigma_{j,\mathbf{a}}].
\end{eqnarray}
Vice versa, if $U$ and $\sigma_{j,\mathbf{c}}$ are as in 2.(a) and satisfy \eqref{4.1}
then also 2.(b) holds.
\end{lemma}

\begin{proof}
The equivalence 1. $\Leftrightarrow$ 2. is clear.
Namely, a locally parabolic abelian bundle is isomorphic, for any $\mathbf{b} \in C$, to a
direct sum $\oplus_{\mathbf{a}\in C^\mathbf{b}}({^\mathbf{a}}\cL)^{d(\mathbf{a})}$. On the other hand,
given a parabolic bundle satisfying 1., then for any $\mathbf{b}\in C$, the sections $\sigma_{j,\mathbf{a}}$
for $\mathbf{a}\in C^\mathbf{b}$ correspond to the choice of such an isomorphism.

Suppose that $U$ and $\sigma_{j,\mathbf{c}}$ are as in 2.
The sheaf ${_\mathbf{a}}\cE_{|U}$ is free with basis
$i_{\mathbf{c}\mathbf{a}}(\sigma_{j,\mathbf{c}})$, $\mathbf{c}\in C^\mathbf{a},j\in\{1,...,d(\mathbf{c})\}$.
The subsheaf $\sum_{\mathbf{c}\lneq\mathbf{a}}{_\mathbf{c}}\cE_{|U}$
is coherent and is generated by
$i_{\mathbf{c}\mathbf{a}}(\sigma_{j,\mathbf{c}})$, $\mathbf{c}\in C^\mathbf{a}\backslash \{\mathbf{a}\},
j\in \{1,...,d(\mathbf{c})\}$, and by
$r_i\cdot \sigma_{j,\mathbf{a}}$, $i\in\uul,
j\in\{1,...,d(\mathbf{a})\}$.
This shows equation \eqref{4.1} and the local freeness of
$\Gr_\mathbf{a}(_\mathbf{a}\cE)$.

Finally, suppose that $U$ and $\sigma_{j,\mathbf{c}}$ are as in 2.
and that we have sections $\www\sigma_{j,\mathbf{c}} \in {_\mathbf{c}}\cE_{|U}$
satisfying 2.(a) and equation \eqref{4.1}.
Fix $\mathbf{a}\in C$. For $\mathbf{b},\mathbf{c}\in C^\mathbf{a}$,
$j\in\{1,...,d(\mathbf{b})\},k\in\{1,...,d(\mathbf{c})\}$, write
$$
i_{\mathbf{b}\mathbf{a}}(\www\sigma_{j,\mathbf{b}}) =
\sum_{k,\mathbf{c}} \kappa_{(j,\mathbf{b}),(k,\mathbf{c})}\cdot
i_{\mathbf{c}\mathbf{a}}(\sigma_{k,\mathbf{c}}).
$$
The square matrix $(\kappa_{(j,\mathbf{b}),(k,\mathbf{c})})$
is holomorphic. It is sufficient to prove that it is invertible
on $D_\uul\cap U$. Then in a neighborhood $\www U\subset U$
of $D_\uul\cap U$ the sections $\www\sigma_{j,\mathbf{b}}$ satisfy 2.(b).

But for fixed $\mathbf{b}\in C^\mathbf{a}$ the block
$(\kappa_{(j,\mathbf{b}),(k,\mathbf{b})})_{|D_\uul\cap U}$ is invertible,
and for fixed $\mathbf{b},\mathbf{c}\in C^\mathbf{a}$ with
$\mathbf{c}\not\leq \mathbf{b}$ the block
$(\kappa_{(j,\mathbf{b}),(k,\mathbf{c})})_{|D_\uul\cap U}$ vanishes.
Therefore the matrix
$(\kappa_{(j,\mathbf{b}),(k,\mathbf{c})})_{|D_\uul\cap U}$ is invertible.
\end{proof}

Given a parabolic bundle $(_\mathbf{a}\cE)_{\mathbf{a}\in C}$, the following
notations will be used:
${_\mathbf{a}}E$ denotes the vector bundle corresponding to ${_\mathbf{a}}\cE$,
$f_{\mathbf{a},\mathbf{b}}$ denotes the morphism corresponding to the inclusion
${_\mathbf{a}}\cE\hookrightarrow {_\mathbf{b}}\cE$ for $\mathbf{a}\leq \mathbf{b}$.

\begin{lemma}\label{t4.4}
Let $(_\mathbf{a}\cE)_{\mathbf{a}\in C}$ be a parabolic bundle with $l=1$
on $(M,D)$. Fix any point $t\in D=D_1$.
\begin{enumerate}
\item
For $\mathbf{a},\mathbf{b}\in C\subset \dC$ with
$\mathbf{a}\leq \mathbf{b}\leq \mathbf{a}+1$
$$
\Imm((f_{\mathbf{a},\mathbf{b}})_{|t}:{_\mathbf{a}}E_t\to{_\mathbf{b}}E_t)
=
\ker((f_{\mathbf{b},\mathbf{a}+1})_{|t}:{_\mathbf{b}}E_t\to{_{\mathbf{a}+1}}E_t).
$$
\item
Fix $\mathbf{a}\in C$. For any $\mathbf{c}\in C^\mathbf{a}$ choose germs of sections
$\sigma_{j,\mathbf{c}}\in {_\mathbf{c}}\cE_t$, $j\in \{1,...,d(\mathbf{c})\}$, such that
the vectors $(\sigma_{j,\mathbf{c}})_{|t}\in {_\mathbf{c}}E_t$ represent a basis of
${_\mathbf{c}}E_t/\sum_{\mathbf{b}<\mathbf{c}}f_{\mathbf{b},\mathbf{c}}({_\mathbf{b}}E_t)$.
Choose a neighborhood $U\subset M$ of $t$ and a coordinate function
$r_1:U\to \dC$ with $D\cap U=\{r_1=0\}$. Define for any $\mathbf{c}\in C$
germs of sections $\sigma_{j,\mathbf{c}}\in{_\mathbf{c}}\cE_t$ by
$$\sigma_{j,\mathbf{c}}=r_1^n\cdot \sigma_{j,\mathbf{c}+n}\quad \textup{ for }\quad
n\in \dZ\textup{ with }\mathbf{c}+n\in C^a.$$
Then part 2.(b) in lemma \ref{t4.3} holds, i.e. for any $\mathbf{b}\in C$
$(i_{\mathbf{c}\mathbf{b}}(\sigma_{j,\mathbf{c}}))_{\mathbf{c}\in C^\mathbf{b},
j\in \{1,...,d(\mathbf{c})\}}$ is a basis of ${_\mathbf{b}}\cE_t$.
\end{enumerate}
\end{lemma}

\begin{proof}
Choose a neighborhood $U\subset M$ of $t$ and a coordinate system
$r=(r_1,...,r_n)$ on $U$ with $r:(U,t)\to (\Delta^n,0)$ an isomorphism
such that $D\cap U=\{r_i=0\}$.

\begin{enumerate}
\item
Consider $\www M:=\{r_2=...=r_n=0\}\subset U$,
$\www D :=\www M\cap D=\{0\}=\{t\}\subset\www M$, and for $\mathbf{a}\in C$
the sheaf ${_\mathbf{a}}\www\cE$ of holomorphic sections of
${_\mathbf{a}}E_{|\www M}$.
Then $(_{\mathbf{a}}\www\cE)_{\mathbf{a}\in C}$ is a parabolic bundle
on $(\www M,\www D)$, and
${_\mathbf{a}}\www E={_\mathbf{a}}E_{|\www M}$,
$\www f_{\mathbf{a},\mathbf{b}}=(f_{\mathbf{a},\mathbf{b}})_{|\www M}$,
and in particular ${_\mathbf{a}}\www\cE ={_{\mathbf{a}+1}}\www\cE(-\www D)
=r_1\cdot {_{\mathbf{a}+1}}\www\cE$.

Consider a germ of a section $\sigma\in {_\mathbf{b}}\cE_{|t}$ with
value $\sigma_{|t}\in {_\mathbf{b}}E_t$. Then
\begin{eqnarray*}
&&\sigma_{|t}\in \ker(f_{\mathbf{b},\mathbf{a}+1})_{|t}\\
&\iff& \sigma \textup{ vanishes at }t\textup{ as a section in }
{_{\mathbf{a}+1}}\www\cE\\
&\iff& \sigma \textup{ is already a section in }{_\mathbf{a}}\www\cE\\
&\iff& \sigma_{|t}\in\Imm(f_{\mathbf{a},\mathbf{b}})_{|t}.
\end{eqnarray*}
This proves part 1.

\item

The vector space ${_\mathbf{b}}E_t$ is naturally filtered by the
subspaces $\{0\}$ and $f_{\mathbf{c},\mathbf{b}}({_\mathbf{c}}E_t)$
for $\mathbf{c}\in C^\mathbf{b}$, with quotients
$$
\Gr_\mathbf{c}({_\mathbf{b}}E_{|t})=
\frac {f_{\mathbf{c}\mathbf{b}}(_\mathbf{c}E_{|t})}
{f_{\www{\mathbf{c}}\mathbf{b}}({_{\www{\mathbf{c}}}}E_{|t})}
$$
where $\www{\mathbf{c}}:=\max C\cap(-\infty,\mathbf{c})_\dC$.
The map $(f_{\mathbf{c}\mathbf{b}})_{|t}$ induces the map
$$
\Gr(f_{\mathbf{c}\mathbf{b}})_{|t}:
\Gr_\mathbf{c}(_\mathbf{c}E_t) \to \Gr_\mathbf{c}(_\mathbf{b}E_t).
$$
This map is an isomorphism for $\mathbf{c}\in C^\mathbf{b}$
because of
$$
f_{\mathbf{c}\mathbf{b}}^{-1}(f_{\www{\mathbf{c}}\mathbf{b}}
(_{\www{\mathbf{c}}}E_t)) =
f_{\mathbf{c}\mathbf{b}}^{-1}(\ker(
(f_{\mathbf{b},\www{\mathbf{c}}+1})_{|t}))
= \ker((f_{\mathbf{c},\www{\mathbf{c}}+1})_{|t})
= f_{\www{\mathbf{c}},\mathbf{c}}(_{\www{\mathbf{c}}}E_t).
$$
Here 1. is used two times.

By construction, for all $\mathbf{c}\in C$ (not only $c\in C^\mathbf{a}$)
the vectors $(\sigma_{j,\mathbf{c}})_{|t}\in {_\mathbf{c}}E_t$
represent a basis of $\Gr_\mathbf{c}(_\mathbf{c}E_t)$.
Because of the isomorphisms
$
\Gr(f_{\mathbf{c}\mathbf{b}})_{|t}:
\Gr_\mathbf{c}(_\mathbf{c}E_t) \to \Gr_{\mathbf{c}}(_\mathbf{b}E_t)
$,
for $\mathbf{c}\in C^\mathbf{b}$ the vectors
$f_{\mathbf{c}\mathbf{b}}(\sigma_{j,\mathbf{c}})_{|t})\in
{_\mathbf{b}}E_t$ represent a basis of
$\Gr_{\mathbf{c}}(_\mathbf{b}E_t)$.
Therefore all the vectors
$(f_{\mathbf{c}\mathbf{b}}((\sigma_{j,\mathbf{c}})_{|t}))_{
\mathbf{c}\in C^\mathbf{b},j\in\{1,...,d(\mathbf{b})\}}$
are a basis of ${_\mathbf{b}}E_t$,
and all the sections
$(f_{\mathbf{c}\mathbf{b}}(\sigma_{j,\mathbf{c}}))_{
\mathbf{c}\in C^\mathbf{b},j\in\{1,...,d(\mathbf{b})\}}$
are a basis of ${_\mathbf{b}}\cE_t$.
\end{enumerate}
\end{proof}

\begin{proof}[Proof of theorem \ref{theoCompatibleLocallyAbelian}]
We will prove condition 2. in lemma \ref{t4.3} by induction on $l$.
The case $l=1$ is handled by part 2. of lemma \ref{t4.4}, which gives
exactly condition 2. in
lemma \ref{t4.3} if the divisor $D$ has only one component.

Now suppose that for some $l\geq 2$, the statement (i.e.,
condition 2. in lemma \ref{t4.3}) is true for any parabolic bundle
on a manifold with a divisor with $l-1$ components.
Let $(_\mathbf{a}\cE)_{\mathbf{a}\in C}$ be a parabolic bundle
on $(M,D)=(M,\bigcup_{i\in\uul}D_i)$. Choose $t\in D_\uul$.
Choose locally around $t\in M$ coordinates $r=(r_1,...,r_n)$
on $M$ and a neighborhood $U\subset M$ of $t$ with
$r:(U,t)\to (\Delta^n,0)$ an isomorphism such that
$D_i\cap U=\{r_i=0\}$ for $i\in\uul$.

For $I\subset \uul, I\neq\emptyset$, and $\mathbf{a}^0\in C$
define the index set
$$
C(I,\mathbf{a}^0):=\{\mathbf{a}\in C\ |\
a_j=a^0_j \textup{ for }j\in \uul-I\}.
$$
Then $(_\mathbf{a}\cE)_{\mathbf{a}\in C(I,\mathbf{a}^0)}$
is a parabolic bundle on $(M,\bigcup_{i\in I}D_i)$.

Consider the case $I=\{l\}$.
The system $(_\mathbf{a}\cE)_{\mathbf{a}\in C(\{l\},\mathbf{a}^0)}$
is locally abelian because of $|I|=1$, and by lemma \ref{t4.3}
any quotient sheaf
$$
\Gr_l(_\mathbf{a}\cE):= \frac{{_\mathbf{a}}\cE}
{\sum_{\mathbf{b}\in C(\{l\},\mathbf{a}^0),
b_l<a_l} {_\mathbf{b}}\cE}
$$
is a locally free $\cO_{D_l}$-module.
Now consider $I=\underline{l-1}=\{1,...,l-1\}$.

\medskip
{\bf Claim:} {\it The system
$(\Gr_l(_{\mathbf{a}}\cE))_{\mathbf{a}\in C(\underline{l-1},\mathbf{a}^0)}$
is a parabolic bundle on
$(D_l,\bigcup_{i\in \underline{l-1}}D_l\cap D_i)=(D_l,D_l-D_l^\circ)$.}

\medskip
\begin{proof}
All the sheaves ${_\mathbf{a}}\cE$,
$\mathbf{a}\in C(\underline{l-1},\mathbf{a}^0)$,
coincide on $M-\bigcup_{i\in\underline{l-1}}D_i$.
Therefore all the sheaves $\Gr_l(_\mathbf{a}\cE)$,
$\mathbf{a}\in C(\underline{l-1},\mathbf{a}^0)$,
coincide on $D_l^\circ$.
Because they are locally free, the natural maps
$\Gr_l(_\mathbf{a}\cE)\to \Gr_l(_\mathbf{b}\cE)$ for
$\mathbf{a}\leq \mathbf{b},
\mathbf{a},\mathbf{b}\in C(\underline{l-1},\mathbf{a}^0),$
are inclusions.
The equations
${_{\mathbf{a}-\mathbf{e}_i}}\cE = {_\mathbf{a}}\cE(-D_i)$
imply
$$
\Gr_l(_{\mathbf{a}-\mathbf{e}_i}\cE) =
\Gr_l(_\mathbf{a}\cE)(-D_l\cap D_i)\quad
\textup{ for }i\in\underline{l-1}.
$$
\end{proof}

\medskip
By induction hypothesis, the parabolic bundle
$(\Gr_l(_\mathbf{a}\cE))_{\mathbf{a}\in C(\underline{l-1},\mathbf{a}^0)}$
is locally abelian.
The equality of quotients
$$
\Gr_\mathbf{a}(\Gr_l(_\mathbf{a}\cE)) = \Gr_\mathbf{a}(_\mathbf{a}\cE)
$$
is obvious. By lemma \ref{t4.3} the quotient on
the left is a locally free $\cO_{D_\uul}$-module
of some finite rank $d(\mathbf{a})$, hence,
also $\Gr_\mathbf{a}(_\mathbf{a}\cE)$ is $\cO_{D_\uul}$-locally free.

Now choose sections $\sigma_{j,\mathbf{a}}$,
$\mathbf{a}\in C,j\in\{1,...,d(\mathbf{a})\}$ of
$_\mathbf{a}\cE_{|U}$ such that
\begin{eqnarray*}
&&r_i\cdot \sigma_{j,\mathbf{c}} = \sigma_{j,\mathbf{c}-\mathbf{e}_i}
\quad\textup{ for }i\in\uul,\mathbf{c}\in C,
j\in\{1,...,d(\mathbf{c})\}\qquad \textup{ and}\\
&&\Gr_\mathbf{c}(_\mathbf{c}\cE_{|U})= \bigoplus_{j\in\{1,...,d(\mathbf{c})\}}
\cO_{D_\uul\cap U}\cdot[\sigma_{j,\mathbf{c}}]
\quad\textup{ for }\mathbf{c}\in C.
\end{eqnarray*}
We have to show that
$(i_{\mathbf{c} \mathbf{a}}(\sigma_{j,\mathbf{c}}))_{\mathbf{c}\in C^\mathbf{a},
j\in\{1,...,d(\mathbf{c})\}}$
is a local basis of ${_\mathbf{a}}\cE_{|U}$ for any $\mathbf{a}\in C$.
Then, by lemma \ref{t4.3}, $2.\Rightarrow 1.$, $(_\mathbf{a}\cE)_{\mathbf{a}\in C}$
is locally abelian.

Denote by $[\sigma_{j,\mathbf{a}}]_l$ the class of $\sigma_{j,\mathbf{a}}$
in $\Gr_l(_\mathbf{a}\cE)$.
Lemma \ref{t4.3} applied to the locally abelian parabolic bundle
$(\Gr_l(_\mathbf{a}\cE))_{\mathbf{a}\in
  C(\underline{l-1},\mathbf{a}^0)}$ shows that
the sections $([i_{\mathbf{c} \mathbf{a}}(\sigma_{j,\mathbf{c}})]_l)_{\mathbf{c}
\in C(\underline{l-1},\mathbf{a}^0)^\mathbf{a},
j\in\{1,...,d(\mathbf{c})\}}$ form a local basis of
$\Gr_l(_\mathbf{a}\cE)$. This holds for arbitrary $\mathbf{a}^0$ and
$\mathbf{a}$ in $C$.

For any $\mathbf{b}^0\in C$, it gives the condition \eqref{4.1} for the locally abelian
parabolic bundle $(_\mathbf{b}\cE)_{\mathbf{b}\in C(\{l\},\mathbf{b}^0)}$
and the sections $i_{\mathbf{c}\mathbf{b}}(\sigma_{j,\mathbf{c}})$,
$\mathbf{b}\in C(\{l\},\mathbf{b}^0), \mathbf{c}\in
C(\underline{l-1},\mathbf{b})^{\mathbf{b}},
j\in\{1,...,d(\mathbf{c})\}$.
The last part of lemma \ref{t4.3} applies and shows that
the sections
$(i_{\mathbf{c}\mathbf{b}}(\sigma_{j,\mathbf{c}}))_{\mathbf{c}\in C^\mathbf{b},
j\in\{1,...,d(\mathbf{b})\}}$
are a local basis of ${_\mathbf{b}}\cE$.
Therefore $(_\mathbf{b}\cE)_{\mathbf{b}\in C}$ is a
locally abelian parabolic bundle.
\end{proof}

\textbf{Remark:}
%\begin{enumerate}
%\item
Borne's proof of theorem \ref{theoCompatibleLocallyAbelian} starts
with %\cite[lemme 6]{Borne2},
\cite[lemme 2.3.11]{Borne2},
which applies to the case $l=1$
and gives that the sheaves $\Gr_l(_\mathbf{a}\cE)$ are locally free
$\cO_{D_l}$-modules.
Then an induction and additional arguments in
\cite[propositions 2.3.5 and 2.3.10]{Borne2}
%\cite[propositions 1 and 2]{Borne2}
give a result which contains
that the sheaves $\Gr_\mathbf{a}(_\mathbf{a}\cE)$ are locally free
$\cO_{D_\uul}$-modules (it treats the homology of a complex
associated to a "facette", and $\Gr_\mathbf{a}(_\mathbf{a}\cE)$
is one homology group).
Then he establishes equivalence between parabolic bundles
and locally free sheaves on certain stacks. To conclude he
shows that these are sums of line bundles on these stacks.
The case $l=1$ of this is \cite[proposition 3.12]{Borne1}.

Our proof unifies the treatment of the quotient sheaves
$\Gr_l(_\mathbf{a}\cE)$ and $\Gr_\mathbf{a}(_\mathbf{a}\cE)$
with the final analysis of the locally free sheaves on the stacks
into one big induction. Therefore our step $l=1$, which is
lemma \ref{t4.4} (and the application of lemma \ref{t4.3} to
the case $l=1$) replaces
%\cite[lemme 6]{Borne2}
\cite[lemme 2.3.11]{Borne2}
and
\cite[proposition 3.12]{Borne1}.
Our inductive step is almost the same as the induction in
\cite[propositions 2.3.5 and 2.3.10]{Borne2}.
%\cite[propositions 1 and 2]{Borne2}

%\item

% In \cite[chapter 4]{Mo2}, a related notion of compatible
% filtration is considered. Given a weakly compatible system
% of vector bundles as above, one can consider for any $\mathbf{b}\in C$
% the filtrations
% $$
% (\Imm(f_{{\mathbf{b}}+(a_j-b_j){\mathbf{e}}_j,{\mathbf{b}}}))_{a_j\in (b_j-1,b_j]}
% $$
% on ${_\mathbf{b}}E$ . Then the system is strongly compatible iff
% these $l$ filtrations are compatible in the sense of
% \cite[definitions 4.2 and 4.9]{Mo2}. In that case
% we also have $\Imm(f_{{\mathbf{b}}+(a_j-b_j){\mathbf{e}}_j,{\mathbf{b}}})
% =\ker(f_{{\mathbf{b}},\mathbf{b}+(a_j+1-b_j){\mathbf{e}}_j})$, and moreover,
% the graded object
% $$
% \frac{{_\mathbf{a}}E}{\sum_{\mathbf{c}\lneq \mathbf{a}}{_\mathbf{c}}E}
% $$
% can also be obtained by taking the successive quotients (in whatever order)
% of these single index filtrations.

%\end{enumerate}

\section{Tame harmonic bundles and limit data}
\label{secLimitProofs}

This section is devoted to the proof of theorem \ref{theoLimitTERP}
and of theorem \ref{theoLimitTwistor}. The first one is essentially
an application of the results of the last section, while the second
one consists in a detailed comparison with the data occurring in
\cite[theorem 12.22]{Mo2}.

First we apply the results from the last section to
the system of sheaves considered in theorem \ref{theoLimitTERP}.
Notice that the first statement of the following lemma
is actually
shown in \cite[lemma 3.3]{SiIy1} in the algebraic context.
\begin{lemma}\label{lemVFLocAbelian}
\begin{enumerate}
\item
Let $H'\in\VB^\nabla_{\dC^*\times Y}$ be a flat bundle. Then
the system of locally free sheaves $({_\mathbf{a}}\cV)_{\mathbf{a}\in C}$
as defined by formula \eqref{eqDefDeligneExtension}
is a locally abelian parabolic bundle on $\dC^*\times X$.
\item
Let $\TERP$ be a variation of TERP-structures on $Y$, tame along $D$. Then the system
$({_\mathbf{a}}\cF)_{\mathbf{a}\in C}$ is a locally abelian parabolic bundle.
\end{enumerate}
\end{lemma}
\begin{proof}
Both system of sheaves $({_\mathbf{a}}\cV)_{\mathbf{a}\in C}$ resp.
$({_\mathbf{a}}\cF)_{\mathbf{a}\in C}$ are obviously parabolic
sheaves, and both are locally free: for
$({_\mathbf{a}}\cV)_{\mathbf{a}\in C}$
this follows from the construction using multi-elementary sections,
and for $({_\mathbf{a}}\cF)_{\mathbf{a}\in C}$ this is exactly the
condition for the variation $\TERP$ to be tame along $D$. Hence, by
theorem \ref{theoCompatibleLocallyAbelian}, both are locally abelian
parabolic bundles. Obviously, for a tame variation of TERP-structures, once we know that
$({_\mathbf{a}}\cF)_{\mathbf{a}\in C}$
is locally abelian, the same is true for
$({_\mathbf{a}}\cV)_{\mathbf{a}\in C}$ as the latter parabolic bundle
is the restriction of the former to $\dC^*\times X$.
On the other hand, an explicit basis of
the sheaves $({_\mathbf{a}}\cV)_{\mathbf{a}\in C}$ which satisfies
the conditions (a) and (b) in lemma \ref{t4.3}, 2. is defined by putting
$\sigma_{j,\mathbf{c}}=es_\uul(A_j,\mathbf{c})$, where
$(A_j)_{j=1,\ldots,d(\mathbf{c})}$ is a local basis of $\cH'(\uul)_{e^{2\pi i \mathbf{c}}}$.
\end{proof}

We can now use the adapted basis constructed in theorem \ref{theoCompatibleLocallyAbelian}
to show the first main result of section \ref{secLimitTwistor}.
\begin{proof}[Proof of theorem \ref{theoLimitTERP}]
It follows from theorem \ref{theoCompatibleLocallyAbelian} and lemma
\ref{t4.3} that the quotients
$\Gr_\mathbf{a}({_\mathbf{a}}\cF)$ are locally free over $\cO_{\cD_\uul}$. Moreover,
we have ${_\mathbf{a}\cF}_{|\dC^*\times X}={_\mathbf{a}}\cV$ by definition, so that
$\Gr_\mathbf{a}({_\mathbf{a}}\cF)_{|\dC^*\times D_\uul} = \Gr_\mathbf{a}({_\mathbf{a}}\cV)$.
Hence we obtain an extension of $\Gr_\mathbf{a}({_\mathbf{a}}V)$ to $\cD_\uul$,
and an extension of $H'(\uul)_{e^{2\pi i \mathbf{a}}}$ via $(\Phi')^{-1}$ from formula \eqref{eqIsoPhi}.
Choose a local basis $(\sigma_{j,\mathbf{c}})$ of ${_\mathbf{a}}\cF$ as in
lemma \ref{t4.3}, 2., and develop any $\sigma_{j,\mathbf{c}}$ as
a sum of $\uul$-elementary sections $\sigma_{j,\mathbf{c}}=\sum_{\mathbf{b}\leq \mathbf{c}}es_\uul(A(\sigma_{j,\mathbf{c}},\mathbf{b}),\mathbf{b})$.
Then $[\sigma_{j,\mathbf{a}}]=[es_\uul(A(\sigma_{j,\mathbf{a}},\mathbf{a}),\mathbf{a})]$ in $\Gr_\mathbf{a}({_\mathbf{a}}\cF)$,
which shows the isomorphism \eqref{eqDescriptionExtensionH}.

It follows from lemma \ref{lemExtensionSheaves_aF} and from
the equations \eqref{eqPiecesElementary-1} and \eqref{eqPiecesElementary-2} that the flat connection
on $\cH'(\uul)_{e^{2\pi i \mathbf{a}}}$ has a pole of type 1
on $\cH(\uul)_{e^{2\pi i \mathbf{a}}}$
along $\{0\}\times D_\uul$. Similarly, equation \eqref{eqPiecesElementary-3}
shows that the endomorphisms $zN_j$ extend holomorphically
to the bundle $\cH(\uul)_{e^{2\pi i \mathbf{a}}}$.

In order to show that $H(\uul)$ underlies a variation of TERP-structures
on $D_\uul$, it only remains to prove that the pairing $P$ has the correct properties
on this bundle. We have already seen that $P:{_\mathbf{a}}\cF\otimes j^*{_\mathbf{b}}\cF\rightarrow
z^{-w}\cO_\cX$ is non-degenerate, where $\mathbf{b}\in C$ such that $C^\mathbf{b}=-C^\mathbf{a}$.
Choose again local bases $(\sigma_{i,\mathbf{c}})_{\mathbf{c}\in C^\mathbf{a}}$ of ${_\mathbf{a}}\cF$
resp. $(\sigma_{j,\mathbf{d}})_{\mathbf{d}\in C^\mathbf{b}}$ of ${_\mathbf{b}}\cF$
as in theorem \ref{theoCompatibleLocallyAbelian}, 1., then
the matrix
$
z^{-w}(P(\sigma_{i,\mathbf{c}}, \sigma_{j,\mathbf{d}}))_{
\mathbf{c}\in C^\mathbf{a},\mathbf{d}\in C^\mathbf{b},
i\in\{1,\ldots d(\mathbf{c})\},j\in\{1,\ldots d(\mathbf{d})\}
}$ is holomorphic and non-degenerate on $\cX$.

These sections can be developed as sums of elementary sections
and formula \eqref{eqElementary-r-4} can be applied. The restriction of the entries
of the above matrix for $\mathbf{d}=-\mathbf{a}$
to the subvariety $\cD_\uul$ takes a particularly simple form, namely:
\begin{eqnarray*}
z^{-w}P(\sigma_{i,\mathbf{c}},\sigma_{j,\mathbf{d}})_{|\cD_\uul}=\left\{
\begin{matrix} 0 & \textup{ if }&\mathbf{c}\neq \mathbf{-a},\\
z^{-w}P(A(\sigma_{i,\mathbf{c}},\mathbf{c}),
A(\sigma_{j,-\mathbf{c}},-\mathbf{c})),
& \textup{ if }&\mathbf{c}=\mathbf{-a}.
\end{matrix}\right.
\end{eqnarray*}
Thus the matrix
$(z^{-w}P(A(\sigma_{i,\mathbf{a}},\mathbf{a}),
A(\sigma_{j,-\mathbf{a}},-\mathbf{a})))_{i,j=1,\ldots,d(\mathbf{a})}$
is holomorphic and non-degenerate on $\cO_{\cD_\uul}$. Therefore the pairing
$$
P:\cH(\uul)_{{\mathbf{e}}^{2\pi i
\mathbf{a}}}\otimes j^*\cH(\uul)_{{\mathbf{e}}^{-2\pi i \mathbf{a}}}
\to z^w\cO_{\cD_\uul}
$$
is non-degenerate. It follows that $(H(\uul),H'(\uul)_\dR,\nabla,P)$ is a variation
of TERP-structures of weight $w$. This finishes the proof of part 1. of theorem \ref{theoLimitTERP}.

The first statement of part 2. is that for any $I\subset \uul$, the same construction
yields a bundle $H(I)$ which underlies a variation of TERP-structures
on $D_I^\circ$. This is proved exactly as in part 1., it only uses the tameness of the original variation
of TERP-structures along the divisor $D \backslash \bigcup_{i\in \uul\backslash I} D_i$ in
$X \backslash \bigcup_{i\in \uul\backslash I} D_i$.

The second statement is that this limit variation of TERP-structures is tame along $D_I\backslash D_I^\circ$.
This follows from the tameness of the original variation along the ``other'' components of the divisor,
i.e., along $\bigcup_{i\in \uul\backslash I} D_i$. It uses remark \ref{remFlatBundleOnI}, the formula
$es_\uul(A,\mathbf{a})=es_{\uul\backslash I}(es_I(A,\mathbf{a}),\mathbf{a})$ and the construction
in part 1.

In a similar way the third statement, i.e., the compatibility of these constructions for
$I,J\subset \uul$, $I,J\neq \emptyset, I\cap J=\emptyset$ can be shown using remark \ref{remFlatBundleOnI}. We leave the details
of the second and the third statement to the reader.
\end{proof}

\medskip

In the second part of this section we give the proof of theorem \ref{theoLimitTwistor}.
We will identify the various objects appearing in \cite[theorem 12.22]{Mo2}
with data defined by a tame variation of TERP-structures on $Y$. Ultimately,
we show that the limit object considered in loc.cit. (which is called
 $\oplus_\mathbf{a} S^{can}_{(\mathbf{a},\boldsymbol{0})}(E)$) is the twistor $\widehat{H(\uul)}$ appearing in theorem
\ref{theoLimitTwistor}, which therefore underlies a polarized mixed twistor structure.
This will show most of the statements of this theorem.

Let us first briefly recall the main objects and results appearing in
\cite[theorem 12.22]{Mo2}.
\begin{deflemma}
\label{lemParabolicExt}
Let $(E, \overline{\partial},\theta ,h)$ be a tame harmonic bundle on
$Y$.
Consider the $\cO_\dC\cAn_Y$-module $\cE':=\cO_\dC\cAn_Y \otimes p^{-1}\cAn(E)$,
where $p:\cY\rightarrow Y$ is the projection.
Let $\cE\in\mathit{VB}_\cY$ be the kernel of
$\overline{\partial}+z\overline{\theta}:\cE' \rightarrow \cE'
\otimes\cO_\dC \cA_Y^{0,1}$. Then for any $\mathbf{a}\in\dR^I$, define
the extension
\begin{equation}\label{eqParabolicExt}
{_\mathbf{a}\cE} := \{s\in j^{(2)}_*\cE\;|\;|s|_{p^*h} \in O (\prod_{i\in I}|r_i|^{-\varepsilon-a_i})\;\;\forall \varepsilon>0\}.
\end{equation}
(recall that $j^{(2)}:\cY\hookrightarrow\cX$).
It follows that $r_i \cdot {_\mathbf{a} \cE} \cong {_{\mathbf{a}-_i}\cE}$, which endows $_\mathbf{a}\cE$ with
an $\cO_\cX$-module structure.
For any $z\in \dC^*$, the restrictions ${_\mathbf{a}\cE^z}:=j_z^*({_\mathbf{a}\cE})$ are $\cO_X$-locally free extensions
of $j_z^* \cE$ over $(z,0)$, where $j_z:\{z\}\times X \hookrightarrow
\cX$ (\cite[theorem 8.59]{Mo2}). However, ${_\mathbf{a}\cE}$ is not $\cO_\cX$-free in general.
The system $({_\mathbf{a}}\cE^z)_{\mathbf{a}\in\dR^l}$ is a locally abelian parabolic bundle
on $\{z\}\times X$ in the sense of definition \ref{defParabolicBundle}.
The operator $z\partial+\theta$ defines a $z$-connection on ${_\mathbf{a}}\cE^z$, which has a
logarithmic pole along $\{z\}\times D$ (\cite[lemma 8.88]{Mo2}). We obtain a tuple of commuting residue
endomorphisms $r_i(z\partial_{r_i}+\theta_{r_i})$ on the graded object $\Gr_\mathbf{a}({_\mathbf{a}}\cE^z)$
for any $\mathbf{a}\in \dR^l$ (which is $\cO_{\cD_\uul}$-locally free).

Denote by $\oplus_{\boldsymbol{\alpha}} \dE_{\boldsymbol{\alpha}} \Gr_\mathbf{a}
(_\mathbf{a}\cE^z)$ the generalized common eigenspace decomposition
(the eigenvalues are constant, \cite[8.8.4]{Mo2}).
Define
$$
\mathit{KMSS}(\cE^z, \uul):=\left\{
(\mathbf{a},\boldsymbol{\alpha})\in(\dR\times \dC)^l\,|\,
\dim_\dC\left(\dE_{\boldsymbol{\alpha}}\Gr_\mathbf{a}
  ({_\mathbf{a}}\cE^z)\right)\neq 0
\right\}
$$
to be the \textbf{Kashiwara-Malgrange-Sabbah-Simpson spectrum}
of $\cE^z$ and by
$\mathbf{m}(\mathbf{a},\boldsymbol{\alpha}):=\dim_\dC\left(\dE_{\boldsymbol{\alpha}}
\Gr_\mathbf{a}(\cE^z)\right)$
the multiplicity of the spectral element $(\mathbf{a},\boldsymbol{\alpha})$ (\cite[8.8.4]{Mo2}).
There is a $\dZ^l$-action on $\mathit{KMSS}$ due to
$$
\mathbf{m}(\mathbf{a},\boldsymbol{\alpha})=\mathbf{m}(\mathbf{a}+\mathbf{k},\boldsymbol{\alpha}-z\cdot\mathbf{k})
\quad\forall\,\mathbf{k}\in\dZ^l.
$$

The behavior of the KMSS-spectrum for varying $z$ is described as follows (\cite[lemma 8.108]{Mo2}):
For any $z\in\dC$ the bijective map
\begin{eqnarray}
\label{eqBehaviorKMSS}
\begin{array}{rcl}
\mathfrak{k}(z):(\dR\times \dC)^l & \longrightarrow &
(\dR \times \dC)^l \\
(\mathbf{a},\boldsymbol{\alpha}) & \longmapsto &
\left(a_j+2\Re(z\cdot\overline{\alpha_j}),\alpha_j-a_j\cdot
  z-\overline{\alpha_j}\cdot z^2\right)_{j\in \uul}
\end{array}
\end{eqnarray}
restricts to a bijection from $\KMSS(\cE^0,\uul)$ to
$\KMSS(\cE^z,\uul)$, preserving the multiplicities.
\end{deflemma}

In the following lemma, we show that for a tame harmonic bundle defined by
a variation of pure polarized TERP-structures, the objects introduced above
simplify to a large extent.
\begin{lemma}\label{lemKMSSforTERP}
Let $\TERP$ be a variation of pure polarized TERP-structures on $Y$ and
suppose that the associated harmonic bundle $E:=H_{|\{0\}\times Y}$ is tame
along $D$. Then:
\begin{enumerate}
\item
The variation of twistor structures constructed from the harmonic
bundle $E$ in \cite[11.1]{Mo2} is the bundle $\widehat{H}$, equipped with
the horizontal parts of the $z$-connection in each fibre $H_{|\{z\}\times Y}$,
and the pairing $\widehat{S}$ from definition \ref{defExtInfinity}.
In particular, the sheaf $\cE$ from above is isomorphic to $\cH$.
\item
The KMSS-spectrum satisfies
$$
\KMSS(\cE^z,\uul) = \left\{(\mathbf{a},-z \cdot \mathbf{a}) \in(\dR\times \dC)^l \,| \Gr_\mathbf{a}({_\mathbf{a}}\cE^0)\neq 0\,\right\}.
$$
In particular, the eigenvalues of the Higgs fields $\theta_{\partial_{r_i}}\in{\cE}\!nd_{\cO_Y}(\cE/z\cE)$ are
all equal to zero.
\item
The variation $H$ is tame along $D$ in the sense of definition \ref{defTERPTame}.
\end{enumerate}
\end{lemma}
\begin{proof}
\begin{enumerate}
\item
This is clear from the definitions, by comparing with the formulas
in \cite[proof of lemma 11.2]{Mo2} for the connection and \cite[lemma
11.9]{Mo2} for the pairings.
\item
On each slice $\{z\}\times Y$ with $z\neq 0$, the connection operator $\partial+z^{-1}\theta$
gives a flat structure on $\cE^z_{|\{z\}\times Y}$, and the extension ${_\mathbf{b}}\cE^z$
has a logarithmic pole along $\{z\}\times D$. The residue eigenvalues at $\{z\}\times D_j$
on $\Gr_\mathbf{b}({_\mathbf{b}}\cE^z)$ are equal to $\alpha_j\cdot
z^{-1}-a_j-\overline{\alpha_j}\cdot z$
for some $(\mathbf{a},\boldsymbol{\alpha})\in\KMSS(\cE^0,\uul)$,
due to formula \eqref{eqBehaviorKMSS}. It follows that the eigenvalues
of the corresponding monodromies around the divisors $\{z\}\times D_j$ are of the form
$\exp(-2\pi i (\alpha_j\cdot z^{-1}-a_j-\overline{\alpha_j}\cdot z))$.
However, as we have $\cE_{|\dC^*\times Y}=\cH'$, all these flat bundles
form an isomonodromic family, so that the monodromies are constant, namely, they are the endomorphisms
$M_i\in \Aut(H^\infty)$ considered at the beginning of section \ref{secLimitTwistor}. We conclude
that the eigenvalues are constant, and thus $\alpha_j=0$ for all $j\in \uul$.
Therefore, only pairs $(\mathbf{a},-z\mathbf{a})$, where $\mathbf{a}\in C$ such that
$\Gr_\mathbf{a}({_\mathbf{a}}\cE^0)\neq 0$ appear as elements of $\KMSS(\cE^z,\uul)$.

As an easy consequence, we obtain that all eigenvalues of the monodromies $M_i$ are
elements in $S^1$, as they are exponentials of the values $a_i$, where $\mathbf{a}\in C\subset \dR^l$
is a vector such that $\Gr_\mathbf{a}({_\mathbf{a}}\cE^0)\neq 0$.
\item
As we have seen in part 2., the eigenvalues of the residue endomorphism $r_i\nabla_{r_i}$ on
${_\mathbf{a}}\cE^z/r_i \cdot {_\mathbf{a}}\cE^z$ are independent of $z\in\dC^*$ and
contained in $-C^\mathbf{a}$. This yields ${_\mathbf{a}}\cE_{|\dC^*\times X} \cong {_\mathbf{a}}\cV$.
Thus ${_\mathbf{a}}\cF_{|\dC^*\times X} \cong {_\mathbf{a}}\cE_{|\dC^*\times X}$ by definition
of the sheaf ${_\mathbf{a}}\cF_{|\dC^*\times X}$ (see definition \ref{defTERPTame}). In order to show
${_\mathbf{a}}\cF \cong {_\mathbf{a}}\cE$, and the local freeness of these sheaves,
we proceed as in \cite[lemma 6.11, 4.]{HS1}.
We already know that ${_\mathbf{a}}\cF_{|\cX\backslash (\{0\}\times D)} \cong {_\mathbf{a}}\cE_{|\cX\backslash (\{0\}\times D)}$.
If $(\mathbf{a},\boldsymbol{0})\notin\KMSS(\cE^0,\uul)$, then ${_\mathbf{a}}\cE$ is $\cO_\cX$-locally free
by  \cite[proposition 1.11]{Mo2}. Suppose therefore that $(\mathbf{a},\boldsymbol{0})\in\KMSS(\cE^0,\uul)$,
then there is $\boldsymbol{\varepsilon}_0\in \dR_{>0}^l$ such that $(\mathbf{a}+\boldsymbol{\varepsilon},0)
\notin\KMSS(\cE^0,\uul)$ for all $\boldsymbol{\varepsilon}\in\dR_{>0}^l$
such that $\boldsymbol{\varepsilon}\leq \boldsymbol{\varepsilon}_0$. Then
${_{\mathbf{a}+\boldsymbol{\varepsilon}}}\cE$ is locally free and moreover, as
${_{\mathbf{a}+\boldsymbol{\varepsilon}}}\cE_{|\cX\backslash (\{0\}\times D)}\cong {_\mathbf{a}}\cF_{|\cX\backslash (\{0\}\times D)}$,
we have ${_{\mathbf{a}+\boldsymbol{\varepsilon}}}\cE \cong {_\mathbf{a}}\cF$
by \cite{Se}. It follows that ${_\mathbf{a}}\cE \subset {_{\mathbf{a}+\boldsymbol{\varepsilon}}}\cE \cong {_\mathbf{a}}\cF$.
However, this is true for all $\boldsymbol{\varepsilon}$ with
$\boldsymbol{0} < \boldsymbol{\varepsilon}\leq
\boldsymbol{\varepsilon}_0$, i.e., any section $s\in{_\mathbf{a}}\cF$
satisfies
$|s|_{p^*h} \in O (\prod_{i\in I}|r_i|^{-\varepsilon_i-\delta-a_i})$ for
all $\delta>0$ and all $\boldsymbol{\varepsilon}\in(\boldsymbol{0},\boldsymbol{\varepsilon}_0]$.
This is exactly the defining property of ${_\mathbf{a}}\cE$, so that
we obtain $s\in{_\mathbf{a}}\cE$. Hence
${_\mathbf{a}}\cE={_\mathbf{a}}\cF$, ${_\mathbf{a}}\cE$ is locally
free and the variation $\TERP$ we started with is tame along $D$.
\end{enumerate}
\end{proof}
In particular, this lemma shows that for a pure polarized variation of TERP-structures
the notion of tameness, as introduced in definition \ref{defTERPTame}, coincides with the notion of tameness
of the associated harmonic bundle, which justifies our terminology.

The next step is the discussion of the limit objects constructed from the sheaves ${_\mathbf{a}}\cE$.
In \cite[8.9.1]{Mo2}, for any $(\mathbf{a},\boldsymbol{\alpha})\in\KMSS(\cE^0,\uul)$,
the vector bundle ${^\uul}\cG_{(\mathbf{a},\boldsymbol{\alpha})}\in\VB_{\cD_\uul}$
is defined. It is characterized by the property that
for any $z\in\dC$, $({^\uul}\cG_{(\mathbf{a},\boldsymbol{\alpha})})_{|\{z\}\times D_\uul}
=\dE_{\boldsymbol{\gamma}}\Gr_\mathbf{c}({_\mathbf{c}}\cE^z)$, where $(\mathbf{c},\boldsymbol{\gamma})
=\mathfrak{k}(z)(\mathbf{a},\boldsymbol{\alpha})$.
If the harmonic bundle $E$ is defined by a variation of pure polarized TERP-structures
as above, then $\mathfrak{k}(z)(\mathbf{a},\boldsymbol{\alpha})=(\mathbf{a},-z\mathbf{a})$,
as we have just proved. It follows that on each $\Gr_\mathbf{a}({_\mathbf{a}}\cE^z)$, there
is just a single generalized common eigenspace of the operators $r_i\nabla_{r_i}$, namely
the one associated to $-z\mathbf{a}$. Therefore $\dE_{-z\mathbf{a}}\Gr_\mathbf{a}({_\mathbf{a}}\cE^z)$ is
identified with $\Gr_\mathbf{a}({_\mathbf{a}}\cE^z)$. This implies that ${^\uul}\cG_{(\mathbf{a},\boldsymbol{0})}$
is simply the quotient $\Gr_\mathbf{a}({_\mathbf{a}}\cE)=\Gr_\mathbf{a}({_\mathbf{a}}\cF)$
which is locally free over $\cO_{\cD_\uul}$ by theorem \ref{theoCompatibleLocallyAbelian} and
lemma \ref{lemVFLocAbelian}, 2.

The limit objects $S^{can}_{(\mathbf{a},\boldsymbol{0})}(E)$ in \cite{Mo2} are obtained
by gluing the bundle ${^\uul}\cG_{(\mathbf{a},\boldsymbol{0})}$ with a similar quotient bundle
on $\dP^1\backslash\{0\}$. The next lemma describes this bundle in the current situation
\begin{lemma}
Let $\TERP$ be as above, and $(E,\overline{\partial},\theta,h)$ the associated harmonic
bundle. Consider also the harmonic bundle $(E,\partial,\overline{\theta},h)$ on $\overline{Y}$.
Then:
\begin{enumerate}
\item
There is an isomorphism
$$
(\cE^\dagger,\overline{\partial}+z\overline{\theta}) \cong
\overline{\gamma^*(\cH,\nabla_Y)},
$$
here $\cE^\dagger\in\VB_{\dP^1\backslash\{0\})\times \overline{Y}}$ is the sheaf
constructed from $(E,\partial,\overline{\theta},h)$ on $\overline{Y}$ in
\cite[11.1.1]{Mo2} and
$\nabla_Y:\cH\rightarrow\cH\otimes z^{-1}\Omega^1_{\cY/\dC}$ is
the horizontal part of the connection operator on $\cH$.
\item
The bundle $j^*\overline{\gamma^*\cH}$, where
$j:\dP^1\rightarrow\dP^1$, $j(z)=-z$,
underlies a variation of pure polarized TERP-structures of weight $w$ on $\overline{Y}$,
which is tame along $\overline{D}\subset \overline{X}$.
\item
The limit object ${^\uul}\cG^\dagger_{(\mathbf{-a},\boldsymbol{0})}$ from \cite[11.2.3]{Mo2} is
equal to $j^*\Gr{_\mathbf{-a}}({_\mathbf{-a}}\cF(j^*\overline{\gamma^*\cH}))$, where ${_\mathbf{b}}\cF(j^*\overline{\gamma^*\cH})$
denotes the extension over $\dC\times\overline{X}$ of order $\mathbf{b}$ of the variation
$j^*\overline{\gamma^*\cH}$ from 2.
\end{enumerate}
\end{lemma}
\begin{proof}
\begin{enumerate}
\item
The flat real subbundle $H'_\dR$ of $H'$ induces a real structure of $E$ (denoted
by $\kappa$ in \cite[theorem 2.19]{He4}) which defines a complex conjugation on sections of $E$.
This gives a complex conjugation on $\cE'$, which interchanges $\partial$ and $\overline{\partial}$
resp. $\theta$ and $\overline{\theta}$ (due to the compatibility of the real structure with the
hermitian metric). Moreover, $\overline{\gamma^*(\cE',\overline{\partial}+z\overline{\theta})}
=(\overline{\gamma^*\cE'},\partial+z^{-1}\theta)$, which yields that
$\overline{\gamma^*\cE}=\cE^\dagger$. On the other hand, we already know that $(\cH,\nabla_Y)\cong(\cE,\partial+z^{-1}\theta)$.
Hence
$$
\overline{\gamma^*(\cH,\nabla_Y)}=\overline{\gamma^*(\cE,\partial+z^{-1}\theta)}=
(\cE^\dagger,\overline{\partial}+z\overline{\theta}).
$$
\item
That $j^*\overline{\gamma^*\cH}$ underlies a variation of TERP-structures
of weight $w$ on $\overline{Y}$ is immediately clear. Moreover,
it follows from 1. that this variation is pure polarized, namely, its corresponding
harmonic bundle on $\overline{Y}$ is $(E,\partial,\overline{\theta},h)$.
This harmonic bundle is obviously tame along $\overline{D}$, as
$(E,\overline{\partial},\theta,h)$ is tame along $D$.
Hence the variation of TERP-structures $j^*\overline{\gamma^*\cH}$ is tame along $\overline{D}$.
\item
The sheaf $j^*\Gr_{\mathbf{-a}}(_{\mathbf{-a}}\cF(j^*\overline{\gamma^*\cH}))$ is $\cO_{\dC\times \overline{D}_\uul}$-locally
free due to 2. and theorem \ref{theoLimitTERP}. For each $z\in \dP^1\backslash\{0\}$, the restriction
$({^\uul}\cG^\dagger_{(\mathbf{-a},\boldsymbol{0})})_{|\{z\}\times \overline{D}_\uul}$
is by definition equal to $\dE_{-z\mathbf{a}}\Gr_{\mathbf{-a}}(_{\mathbf{-a}}\cE^\dagger)$ (by \cite[11.2.1, 11.2.3]{Mo2}
and the same argument as above, i.e., the special behavior of the KMSS-spectrum for $\cE^\dagger$ in the current
situation), so that ${^\uul}\cG^\dagger_{(\mathbf{-a},\boldsymbol{0})}\cong j^*\Gr{_\mathbf{-a}}({_\mathbf{-a}}\cF(j^*\overline{\gamma^*\cH}))$ .
\end{enumerate}
\end{proof}
\begin{proof}[Proof of theorem \ref{theoLimitTwistor}]
Part 1. has already been shown in lemma \ref{lemKMSSforTERP}. In order to prove the remaining parts of the theorem,
we will show that for any $r\in D_\uul$, there is an isomorphism of twistors
\begin{equation}\label{eqIdentTwistors}
\widehat{\cH(\uul)}_{|\dP^1\times\{r\}}
\cong
\bigoplus_{\mathbf{a}\;\textup{mod}\;\dZ^l} S^{can}_{(\mathbf{a},\boldsymbol{0})}(E_{|\pi_Y^{-1}(r)}),
\end{equation}
where $S^{can}_{(\mathbf{a},\boldsymbol{0})}(E_{|\pi_Y^{-1}(r)})$ is defined in \cite[11.3.4]{Mo2}.
As in loc.cit., we will make the assumption that $D_\uul=\{0\}$, so that actually
we have to show that $\widehat{\cH(\uul)}\cong \oplus_{\mathbf{a}\;\textup{mod}\;\dZ^l} S^{can}_{(\mathbf{a},\boldsymbol{0})}(E)$.
Recall that $S^{can}_{(\mathbf{a},\boldsymbol{0})}(E)$ is obtained by
gluing ${^\uul}\cG_{(\mathbf{a},\boldsymbol{0})}$ with
${^\uul}\cG^\dagger_{(\mathbf{-a},\boldsymbol{0})}$ via the gluing map
\begin{equation}\label{eqGluing-Moch}
\xymatrix@!0
{
\left[{^\uul}\cG_{(\mathbf{a},\boldsymbol{0})}\right]_{|\dC^*}
\ar[rrrrrrr]^{\Phi^{\dagger,can}_{(\mathbf{-a},\boldsymbol{0})}\circ(\Phi^{can}_{(\mathbf{a},\boldsymbol{0})})^{-1}}
&&&&&&&
\left[{^\uul}\cG^\dagger_{(\mathbf{-a},\boldsymbol{0})}\right]_{|\dC^*}
}
\end{equation}
where the maps
$\Phi^{can}_{(\mathbf{a},\boldsymbol{0})}:{^\uul}\cG_{(\mathbf{a},\boldsymbol{0})}\cH\longrightarrow
\left[{^\uul}\cG_{(\mathbf{a},\boldsymbol{0})}\right]_{|\dC^*}$
resp. $\Phi^{\dagger,can}_{(\mathbf{-a},\boldsymbol{0})}:
{^\uul}\cG^\dagger_{(\mathbf{-a},\boldsymbol{0})}(\cH^\dagger)\longrightarrow
\left[{^\uul}\cG^\dagger_{(\mathbf{-a},\boldsymbol{0})}\right]_{|\dC^*}$
are defined in \cite[10.4]{Mo2}.
Moreover, there is an identification
\begin{equation}\label{eqIdentFlatSection-Moch}
{^\uul}\cG_{(\mathbf{a},\boldsymbol{0})}\cH\cong {^\uul}\cG^\dagger_{(\mathbf{-a},\boldsymbol{0})}(\cH^\dagger),
\end{equation}
so that the composition $\Phi^{can,\dagger}_{(\mathbf{-a},\boldsymbol{0})}\circ(\Phi^{can}_{(\mathbf{a},\boldsymbol{0})})^{-1}$
is well-defined.
It turns out that in the current situation we have
$$
\begin{array}{rcl}
{^\uul}\cG_{(\mathbf{a},\mathbf{0})}\cH & = & \cH'(\uul)_{e^{2\pi i \mathbf{a}}}, \\ \\
\Phi^{can}_{(\mathbf{a},\mathbf{0})} = \Phi': \cH'(\uul)_{e^{2\pi i \mathbf{a}}} & \longrightarrow & \Gr_\mathbf{a}({_\mathbf{a}}\cV)
= \Gr_\mathbf{a}({_\mathbf{a}}\cF)_{|\dC^*}, \\ \\
{^\uul}\cG^\dagger_{(\mathbf{-a},\mathbf{0})}(\cH^\dagger) & = & \overline{\gamma^*\cH'(\uul)_{e^{-2\pi i \mathbf{a}}}}, \\ \\
\Phi^{\dagger, can}_{(\mathbf{-a},\mathbf{0})} = \overline{\gamma^*\Phi'}  :
\overline{\gamma^*\cH'(\uul)_{e^{-2\pi i \mathbf{a}}}} & \longrightarrow & \overline{\gamma^*\Gr_\mathbf{-a}({_\mathbf{-a}}\cV)}
= j^*\Gr_{\mathbf{-a}}({_\mathbf{-a}}\cF(j^*\overline{\gamma^*\cH}))_{|\dC^*},
\end{array}
$$
and that the identification \eqref{eqIdentFlatSection-Moch} is just the map
$$
\tau:\cH'(\uul)_{e^{2\pi i \mathbf{a}}} \longrightarrow \overline{\gamma^*\cH'(\uul)_{e^{-2\pi i \mathbf{a}}}}
$$
which is induced by the original map $\tau:\cH'_{e^{2\pi i \mathbf{a}}} \rightarrow \overline{\gamma^*\cH'_{e^{-2\pi i \mathbf{a}}}}$.
The twistor $\widehat{\cH(\uul)}_{e^{2\pi i \mathbf{a}}}$ is obtained by gluing $\cH(\uul)_{e^{2\pi i \mathbf{a}}}$
and $\overline{\gamma^*\cH(\uul)_{e^{-2 \pi i \mathbf{a}}}}$ via this new morphism $\tau$. Thus we get an isomorphism
of twistors
$$
\widehat{\cH(\uul)}_{e^{2\pi i \mathbf{a}}} \cong S^{can}_{(\mathbf{a},\mathbf{0})}(E).
$$

Next we have to identify the pairings and the nilpotent maps.
Mochizuki establishes pairings on ${^\uul}\cG_{\mathbf{u}}$,
essentially via \eqref{eqElementary-r-4}. In our situation, this
boils down to the pairings $\whh S$ which are induced on
$\cH'(\uul)$ from those on $\cH'$. The pairing $S$ in
\cite[theorem 12.22]{Mo2} coincides with $\whh S$ on
$\whh\cH(\uul)$.

In \cite[11.3.6]{Mo2} the morphisms
$\cN_j^\Delta:S_{\mathbf{({\mathbf{a}},\boldsymbol{0})}}\to
S_{\mathbf{({\mathbf{a}},\boldsymbol{0})}}\otimes
\cO_{\dP^1}\cC^{an}_{D_\uul}(1,1)$
are defined via extension of the nilpotent parts of the
residue endomorphisms $(iz)^{-1}[zr_j\nabla_{r_j}]$.
The nilpotent parts of the residue endomorphisms $[zr_j\nabla_{r_j}]$
correspond by formula \eqref{eqPiecesElementary-3}
to $z\frac{-N_j}{2\pi i }$. Therefore the pull back of
$\cN_j^\Delta$ to $\cH'(\uul)$ is equal to $\frac{N_j}{2\pi }$.
The tuple $\mathbf{N}^\Delta=(\cN_1^\Delta,\ldots,\cN_l^\Delta)$
in \cite[theorem 12.22]{Mo2} thus corresponds to the tuple
$(\frac{N_1}{2\pi },\ldots,\frac{N_l}{2\pi }).$

Theorem 12.22 in \cite{Mo2} says that
$(S^{can}_{(\mathbf{a},\boldsymbol{0})}(E),W,\mathbf{N}^\Delta,S)$ is a polarized mixed
twistor of weight 0 in $l$ variables.
By definition \cite[definition 3.50]{Mo2} this means
that for any ${\mathbf{a}}\in (\dR^+)^l$ and
$N_{\mathbf{a}}:=\sum_{j\in\uul}a_jN_j$
the tuple $(S_{({\mathbf{a}},\boldsymbol{0})},N_{\mathbf{a}}^\Delta,S)$
is a polarized mixed
twistor of weight 0 and that the weight filtration $W$ is
independent of the choice of  ${\mathbf{a}}\in (\dR^+)^l$.
This includes that
the maps $zN_{\mathbf{a}}$ and $z^{-1}N_{\mathbf{a}}$ extend to
$\{0\}\times D_\uul$
respectively to $\{\infty\}\times D_\uul$,
and that they have everywhere the same Jordan normal form
so that together they give a global weight filtration.
As a conclusion, we obtain part 2. of the theorem, and
also part 3. Finally, part 4. is an easy consequence of
lemma \ref{lemTERPWeightFilt}.
\end{proof}

\section{Rigidity} \label{secRigidity}

As an application of the discussion on extensions of TERP-structures,
we prove here a generalized version of a conjecture of Sabbah concerning
the rigidity of integrable variations of twistor structures on quasi-projective varieties.
It was stated in \cite[conjecture 7.2.9]{Sa6}
for non-compact curves, but using the results of \cite{Mo3}, we
can actually prove it in this more general situation. We show the corresponding
statement for TERP-structures, the original formulation in \cite{Sa6}
can be easily obtained by a slight modification of our proof.

A quite simple but essential ingredient is the following lemma.
\begin{lemma}
\label{lemCompatibility}
Let, as in section \ref{secLimitTwistor}, $X=\Delta^n$, $Y=(\Delta^*)^l\times \Delta^{n-l}$ and $D=X\backslash Y$. Let
$\TERP$ be a variation of pure polarized TERP-structures on $Y$, tame along $D$.
Let $\cK:=\cH/z\cH$.
Consider the parabolic filtration on $j^{(0)}_*\cK$ (where $j^{(0)}:Y\hookrightarrow X$),
defined, analogously to formula \eqref{eqParabolicExt} by:
$$
_{\mathbf{a}}\cK:=\left\{s\in j^{(0)}_*\cK\,|\,|s|_h\in
O (\prod_{i\in I}|r|^{-\varepsilon-a_i})\;\;\forall \varepsilon>0\right\}.
$$
for any $\mathbf{a}\in\dR^\uul$, where $h$ is the hermitian metric induced
on $\cK$ by $z^{-w}P(-,\tau-)$ on $p_*\widehat{\cH}$.
Then the endomorphism $\cU:=[z^2\nabla_z]\in\End_{\cO_Y}(\cK)$
is compatible with this parabolic filtration,
i.e., for any $\mathbf{a}\in\dR^I$ it extends to an element in
$\End_{\cO_X}(_\mathbf{a}\cK)$.
\end{lemma}
\begin{proof}
This is a direct consequence of the results of the preceding sections:
Recall that we defined ${_\mathbf{a}}\cF:=j^{(1)}_*{_\mathbf{a}}\cV \cap j^{(2)}_* \cH\in\VB_\cX$.
By lemma \ref{lemExtensionSheaves_aF} we have that
$$
\nabla: {_\mathbf{a}}\cF \longrightarrow {_\mathbf{a}}\cF \otimes z^{-1}\Omega^1_\cX\left(\log (\cD\cup(\{0\}\times X))\right),
$$
in particular, $(z^2\nabla_z)({_\mathbf{a}}\cF) \subset {_\mathbf{a}}\cF$.
By lemma \ref{lemKMSSforTERP}, 3., we have the equality ${_\mathbf{a}}\cF\cong
{_\mathbf{a}}\cE$, from which it follows that
$_\mathbf{a}\cK={_\mathbf{a}}\cF/z\cdot{_\mathbf{a}}\cF$, so that we obtain
$\cU({_\mathbf{a}\cK})\subset {_\mathbf{a}}\cK$, as required.
\end{proof}
The following theorem is the generalization of \cite[corollary 7.2.8]{Sa6} to the
higher-dimensional quasi-projective case.
\begin{theorem}
\label{theoRigidity}
Let $X$ be a projective manifold and $Y:=X\backslash D$ where $D$ is a divisor with normal crossings.
Let $\TERP$ be a variation of pure polarized TERP-structures on $Y$, tame along $D$,
and denote by $(E,\overline{\partial},\theta,h)$ the corresponding harmonic bundle.
Then there is a decomposition of $\cAn_Y$-bundles
$E=E_w\oplus E_{w-1}$ where $E_w$ resp. $E_{w-1}$ underlies
a variation of pure polarized Hodge structures of weight $w$ resp. $w-1$.
\end{theorem}
\begin{proof}
The variation of TERP-structures $\TERP$ corresponds by \cite[theorem 2.19]{He4}
to a CV$\oplus$-structure. The latter consists of the harmonic bundle $E$ enriched by the endomorphisms
$\cU,\cQ\in\End_{\cAn_Y}(\cAn(E))$ and the real structure given by $\tau$ subject to
a couple of compatibility conditions, from which we quote only the following ones
(\eqref{eqInt1} and \eqref{eqInt2} are called integrability equations
in \cite{Sa6}).
\begin{eqnarray}
\label{eqInt1}[\theta,\cU]=0,&&D''(\cU)=0,\\
\label{eqInt2}D'(\cU)-[\theta,\cQ]+\theta=0,&&D'(\cQ)+[\theta,\tau\cU\tau]=0,\\
\label{eqIntTERP1}h(\cU-,-)=h(-,\tau\cU\tau-),&&h(\cQ-,-)=h(-,\cQ-),\\
\label{eqIntTERP2}\tau\theta\tau=\overline{\theta},&&\cQ=-\tau\cQ\tau.
\end{eqnarray}
Consider, as before, the bundle $\cK=\cH/z\cH\in\VB_Y$ and the extensions
$_\mathbf{a}\cK\in\VB_X$. By \cite[proposition 5.1]{Mo3}, any $_\mathbf{a}\cK$ is a $\mu_L$-polystable
Higgs bundle ($L$ being some fixed ample line bundle on $X$), and we have a canonical decomposition
\begin{equation}\label{eqDecompHiggs}
(_\mathbf{a}\cK,\theta) \cong \oplus_{i=1}^m (_\mathbf{a}\cK_i,\theta_i)\otimes \dC^{p_i}
\end{equation}
where each $_\mathbf{a}\cK_i$ is a $\mu_L$-stable Higgs bundle,
and any two $({_\mathbf{a}}\cK_i,\theta_i)$ and
$({_\mathbf{a}}\cK_j,\theta_j)$ are non-isomorphic for $i\neq j$.
In particular, $({_\mathbf{a}}\cK_i,\theta_i)$ is simple, i.e., any Higgs endomorphism respecting
the filtration $_\bullet\cK_i$ of $_\mathbf{a}\cK_i$ is of the form
$c_i\cdot\Id_{\cK_i}$ with $c_i\in\dC$.
(This follows from a standard argument: By loc.cit., lemma 3.10, any non-trivial endomorphism
is actually an isomorphism, which makes $\End_{\cO_X}(_\mathbf{a}\cK_i)$ into a skew field
which is a finite dimensional $\dC$-vector space as $X$ is compact. Any such isomorphism
is then necessarily a multiplication by a constant for otherwise
it would generate a commutative and finite dimensional subalgebra of $\End_{\cO_X}(_\mathbf{a}\cK_i)$, i.e.,
a proper algebraic field extension of $\dC$, a contradiction).
By restriction to $Y$, we obtain a decomposition $(\cK,\theta)\cong \oplus_{i=1}^m (\cK_i,\theta_i)\otimes
\dC^{p_i}$ and similarly a decomposition $\cAn(E)\cong\oplus_{i=1}^m \cAn(E_i)\otimes\dC^{p_i}$
of the corresponding sheaf of $\cAn_Y$-sections. Moreover,
proposition 5.1 of loc.cit. also gives that the hermitian metric $h$
decomposes as $h=\sum_{i=1}^m h_i\otimes g_i$, where $g_i$ is a constant
hermitian metric on $\dC^{p_i}$. This implies that the $(1,0)$-part
$D'$ of the Chern connection decomposes as $D'=\sum_{i=1}^m D'_i\otimes\partial$, where
$D'_i:\cAn(E_i)\rightarrow\cAn(E_i)\otimes\cA^{1,0}_Y$.
Notice that it follows from equation \eqref{eqDecompHiggs} that
$\theta=\sum_{i=1}^m\theta_i\otimes \Id_{\dC^{p_i}}$.

Consider $\cU$ and $\overline{\cU}:=\tau\cU\tau$ as elements
in $\End_{\cAn_Y}(\cAn(E))$. It follows from the previous
lemma that $\cU$ is an endomorphism of the $\mu_L$-polystable Higgs bundle
$\cK$, so that it decomposes as $\cU=\sum_{i=1}^m\Id_{\cAn(E_i)}\otimes U_i$,
where $U_i\in\End_\dC(\dC^{p_i})$. The first part of equation \eqref{eqIntTERP1} shows that
$\overline{\cU}=\sum_{i=1}^m\Id_{\cAn(E_i)}\otimes \overline{U}_i$,
$\overline{U}_i$ being the adjoint of $U_i$ with respect to $g_i$. It follows
that the commutators $[D',\cU]$ and $[\theta,\overline{\cU}]$ vanish
so that the integrability equations \eqref{eqInt2} reduce to
\begin{equation}\label{eqReducedIntEqu1}
D'(\cQ)=0\quad\textup{and}\quad\theta=[\theta,\cQ],
\end{equation}
and by adjunction with respect to $h$ we also obtain
\begin{equation}\label{eqReducedIntEqu2}
D''(\cQ)=0\quad\textup{and}\quad\overline{\theta}=-[\overline{\theta},\cQ].
\end{equation}
The remaining part of the proof is exactly the
same as the proof of \cite[lemma 3.4 and theorem 3.1]{He4}, where the stronger
assumption of $\cU=0$ was made. For the readers convenience, we briefly recall
how to obtain the desired conclusion. First define the following real analytic
bundles on $X$:
$$
\begin{array}{rcccccl}
\cAn(E_w^{p,w-p}) & := & \bigoplus\limits_{
\begin{array}{c}
\SSC \lfloor\alpha+\frac{w+1}{2}\rfloor=p\\
\SSC \alpha\notin\frac{w+1}{2}+\dZ
\end{array}} {\cK}\!er(\cQ-\alpha\Id) & , & E_w & := & \bigoplus_{p} E_w^{p,w-p}\\ \\
\cAn(E_{w-1}^{p,w-1-p}) & := &
{\cK}\!er(\cQ-(p-\frac{w-1}{2})\Id) & , & E_{w-1} & := & \bigoplus_{p} E_{w-1}^{p,w-1-p}.
\end{array}
$$
Moreover, we put
$$
\begin{array}{rccccccl}
F_w^p & := & \bigoplus_{q\geq p} E_w^{q,w-q} & ; & F_{w-1}^p & := & \bigoplus_{q\geq p} E_{w-1}^{q,w-1-q}.
\end{array}
$$
{}From the equations \eqref{eqReducedIntEqu1}
and \eqref{eqReducedIntEqu2} we deduce that the bundles $F^p_w$ resp. $F^p_{w-1}$
are holomorphic (for the operator $D''+\overline{\theta}$) and satisfy
Griffiths transversality, i.e.,
$\nabla \cF^p_w \subset \cF^{p-1}_w\otimes \Omega^1_Y$ resp. $\nabla \cF^p_{w-1} \subset \cF^{p-1}_{w-1}\otimes \Omega^1_Y$,
where $\nabla$ is the integrable operator $\nabla=D'+D''+\theta+\overline{\theta}$.
{}From the second part of equation \eqref{eqIntTERP2} we deduce that
$\tau{\cK}\!er(\cQ-\alpha\Id)={\cK}\!er(\cQ+\alpha\Id)$. Using that for any $\alpha\in
\dR\backslash\left(\frac{w+1}{2}+\dZ\right)$, the
equality
$$
-\left(	\left[\alpha+\frac{w+1}{2}\right]-p\right)=\left[-\alpha+\frac{w+1}{2}\right]-(w-p)
$$
holds, this implies
$E_w^{p,w-p}=\overline{E_w^{w-p,p}}$ resp. $E_{w-1}^{p,w-1-p}=\overline{E_{w-1}^{w-1-p,p}}$,
which yields the desired result.

%In particular, we have
%$\cAn(E_w) :=  \bigoplus_{p} \cAn(E_w^{p,w-p})$ resp.
%$\cAn(E_{w-1}) :=  \bigoplus_{p} \cAn(E_{w-1}^{p,w-1-p})$, where

\end{proof}
As an application, we obtain a generalization of \cite[corollary 4.5]{HS2}.
Notice that the reasoning is completely different, namely, we do not use the
curvature computation of loc.cit.
\begin{corollary}
Let $H$ be a variation of pure polarized TERP-structures $H$ on
$\dC^n$, tame along $\dP^n\backslash \dC^n$.
Then it is trivial, i.e., $\nabla_X(\cH)\subset \cH$ for all $X\in p^{-1}\cT_{\dC^n}$,
where $p:\dC\times\dC^n \rightarrow \dC^n$ is the projection.
If $H$ is regular singular (then by proposition \ref{propRegSingTame} tameness follows as soon as we suppose
that $H$ is pure polarized), then the period map $\dC^n \rightarrow \MBLp$, which will
be defined in section \ref{subsecApplications}, is constant.
\end{corollary}
\begin{proof}
Let $Y:=\dC^n\subset X=\dP^n$, then the assumptions ensure that $p_*\widehat{\cH}$ underlies
a tame harmonic bundle on $Y$. The last theorem gives that it also underlies a
sum of two variations of pure polarized Hodge structures.
A classical result (see, e.g., \cite[13.4.3]{CarlStachPeters}) shows that
the corresponding period maps to the classifying spaces of
pure polarized Hodge structures are constant, which implies that $\theta=0$,
so that the variation of TERP-structure itself is flat in parameter direction.
The last statement follows directly from the construction of the period map in
lemma \ref{lemPeriodMaps}.
\end{proof}

%\clearpage

\section{The compact classifying space}
\label{secCompactClassSpace}

The fact that special members of families of regular singular
TERP-structures can have different spectral numbers than
the general member of the family is reflected in the geometry
of classifying spaces of such TERP-structures. We construct and study
in this section several versions of these classifying spaces.
In all cases, the topological data of a family of regular singular
TERP-structures are fixed. If we fix moreover the spectral pairs
(as introduced in definition \ref{defHodgeSpectrum}) and a
corresponding reference regular singular mixed TERP-structure,
then all possible regular singular TERP-structures
with these data are classified by a complex manifold $\cbl$, which
was defined in \cite{He2} and further studied in \cite{HS2}. It has the
structure of an affine fibre bundle over a complex homogeneous manifold
$\DcPMHS$ parameterizing certain Hodge type filtrations. More precisely,
$\DcPMHS$ contains an open submanifold $\DPMHS$ which parameterizes
polarized mixed Hodge structures with the same fixed topological data from above,
in particular, with a fixed weight filtration. Passing to the induced
filtration on the graded parts of this weight filtration defines
a structure of an affine fibre bundle $\DcPMHS \rightarrow \DcPHS$, where the
latter is a projective manifold (the product of classifying
spaces of Hodge-like filtrations). This map restricts
to $\DPMHS \rightarrow \DPHS$, and $\DPHS$ is a product
of classifying spaces of pure polarized Hodge structures.
The following diagram shows how these manifolds are related.
\begin{equation}
\begin{array}{ccccc}
\cbl & \longrightarrow & \DcPMHS &
\longrightarrow & \DcPHS\\
\cup & & \cup & & \cup \\
D_\mathit{BL} & \longrightarrow & \DPMHS &
\longrightarrow &  \DPHS.
\end{array}
\end{equation}
We refer to \cite{He2} and \cite[chapter 2]{HS2} for more details about these classifying
spaces. $\cbl$ carries a tautological bundle $\cL\in\VB_{\dC\times\cbl}$
of regular singular TERP-structures, i.e., $\cL_{|\dC\times\{x\}}$ is the TERP-structure
which corresponds to the point $x\in\cbl$. $\cL$ underlies a family of TERP-structures
in the sense of definition \ref{defFamilyTERP}, but not a variation in general.

In order to capture the jumping phenomena of the spectrum as seen
in the examples in section \ref{secBasics}, we will construct a new classifying space which
parameterizes all regular singular TERP-structures where only the range
for the spectral numbers has been fixed. We will
see that this is a projective variety and that it contains the
classifying spaces $\cbl$ for fixed spectral numbers as
locally closed subvarieties. It also contains other strata, these
correspond to families of TERP-structures with fixed spectral pairs,
but where no element is mixed TERP.

All along this section, we fix the following topological data:
a real vector space $H^\infty_\dR$ of dimension $\mu$,
equipped with an automorphism $M\in\mathit{Aut}(H^\infty_\dR)$, an integer $w$ and
a non-degenerate bilinear pairing $S:H^\infty_\dR\times H^\infty_\dR\rightarrow \dR$.
$S$ is required to be invariant under $M$, and to have the following symmetry
property: Denote by $H^\infty$ the complexification of $H^\infty_\dR$,
by $H^\infty_\lambda$ the generalized eigenspaces of $M$, then $S$ is
$(-1)^{w-1}$-symmetric on $H^\infty_{\arg\neq 0}$ and $(-1)^w$-symmetric on
$H^\infty_{\arg=0}$, where, as before $H^\infty_{\arg=0}:=\oplus_{\arg\lambda=0} H^\infty_\lambda$
and $H^\infty_{\arg \neq 0}:=\oplus_{\arg\lambda \neq 0} H^\infty_\lambda$.
By \cite[lemma 5.1]{HS1}, these data correspond to a flat vector bundle
$H'\in\VB^\nabla_{\dC^*}$ with a flat real subbundle $H'_\dR$ of maximal rank, and a flat $(-1)^w$-symmetric
non-degenerate pairing $P:\cH'\otimes j^*\cH'\rightarrow \cO_{\dC^*}$ which takes values
in $i^w\dR$ on $H'_\dR$.

We denote, by abuse of notation, both of the two inclusions $\dC^*\hookrightarrow
\dC$ and $\dP^1\backslash\{0\}\hookrightarrow \dP^1$ by $i$,
and similarly by $\widetilde{i}$ either of the two inclusions $\dC^*\hookrightarrow \dP^1\backslash\{0\}$
or $\dC\hookrightarrow\dP^1$. We consider
the Deligne extensions $V^\alpha, V^{>\alpha}\subset i_*\cH'$ as defined in section \ref{secBasics}.
We also consider the corresponding Deligne extensions
$V_\alpha, V_{<\alpha}\subset \widetilde{i}_*\cH'$ at infinity, where the indices
are chosen so that they form an increasing filtration.
Finally, we work with the meromorphic bundles $V^{>-\infty}\subset i_*\cH'$,
$V_{<\infty}\subset \widetilde{i}_*\cH'$ and with
the sheaf $\widetilde{i}_* V^{>-\infty} \cap  i_* V_{<\infty}$
(here the intersection takes place in $\widehat{i}_*\cH'$,
where $\widehat{i}:\dC^*\hookrightarrow\dP^1$),
which is an algebraic vector bundle over $\dC^*$. We write
$W$ for its space of global sections, then
$W$ is a free $\dC[z,z^{-1}]$-module of rank $\mu$.
Denote for any $\alpha,\beta\in\dC$ the intersection $\widetilde{i}_* V^\alpha \cap i_* V_\beta\in\VB_{\dP^1}$
of subsheaves of $\widehat{i}_*\cH'$ by $V^\alpha_\beta$
(and similarly $V^{>\alpha}_\beta$, $V^\alpha_{<\beta}$ etc.).
For any $\alpha\geq \beta$, we have the following exact sequence
\begin{equation}
\label{eqDeligneLattices}
0\longrightarrow V^{>\alpha}_{\alpha}
\longrightarrow  V^{\beta}_{\alpha}
\longrightarrow  V^{\beta}_\alpha  / V^{>\alpha}_\alpha
\longrightarrow  0
\end{equation}
Obviously, $V^{>\alpha}_{\alpha}$ is semi-stable of weight $-1$,
so that $H^0(\dP^1,V^{>\alpha}_{\alpha})=
H^1(\dP^1,V^{>\alpha}_{\alpha})=0$. This implies that
we have a canonical isomorphism
from $W^\beta_\alpha:=H^0(\dP^1,V^{\beta}_{\alpha})$ to the skyscraper sheaf
$V^{\beta}_\alpha/ V^{>\alpha}_\alpha$ (which we identify with $V^\beta/V^{>\alpha}$).
This isomorphism will be used implicitly many times in the sequel.
The restriction of $P$ to $W$ can be written as
$$
P=\sum_{k\in\dZ} P^{(k)} z^k :W \otimes j^*W \rightarrow \dC[z,z^{-1}]
$$
In particular, $P^{(k)}$ is $(-1)^{w+k}$-symmetric and induces a pairing
\begin{equation}
\label{eqPairingPCoupling1}
P^{(k)}: W^{\alpha}_\alpha \otimes W^\beta_\beta \rightarrow \delta_{\beta,k-\alpha}\dC
\end{equation}
which is non-degenerate for $\alpha+\beta=k$. For any fixed $\alpha,\beta\in\dC$ with $\alpha\leq\beta$ and $\alpha+\beta=k$,
we obtain a non-degenerate $(-1)^{w+k}$-symmetric pairing
\begin{equation}
\label{eqPairingPCoupling2}
P^{(k)}:W^\alpha_\beta\otimes W^\alpha_\beta\longrightarrow \dC
\end{equation}
which is the sum
$$
\sum_{\gamma\in[\alpha,\beta]} \left(P^{(k)}:W^{\gamma}_\gamma \otimes W^{\alpha+\beta-\gamma}_{\alpha+\beta-\gamma} \longrightarrow \dC\right).
$$

We choose $\alpha_1\in \dC$ satisfying the following conditions:
$e^{-2\pi i \alpha_1}$ is required to be an eigenvalue
of $M$ and $\alpha_1\leq\frac{w}{2}$. Then
we put $\alpha_\mu:=w-\alpha_1\geq\frac{w}{2}\geq\alpha_1$
and $n:=\lfloor \alpha_\mu-\alpha_1\rfloor\geq 0$.

The classifying space we are going to consider in this section will
represent a certain functor of families of TERP-structures with
trivial monodromy in parameter direction. We first define this
functor.
\begin{definition}\label{defFuncMBL}
Fix $H^\infty_\dR, S, M, w, \alpha_1$ or the equivalent data
$H', H'_\dR, \nabla, P, w, \alpha_1$ from above, and consider the Deligne
extensions $V^\alpha$ and $V^{>\alpha}$ of $H'$.
Define the functor
$\MBLfun^{H^\infty_\dR,S,M,w,\alpha}$ (which we denote usually by $\MBLfun$
if no confusion can occur) from the category of complex
spaces to the category of sets by
$$
\begin{array}{rcl}
\MBLfun(X) & := &\Big\{(\cL,\varphi)\,|\,\cL\in\VB_{\dC^*\times X}, \varphi:i^*\cL \stackrel{\cong}{\longrightarrow}
  (p')^*H',(z^2\varphi^*\nabla_z)\cL\subset\cL,\\ \\
 & & \;\;\varphi^*P:\cL\otimes j^*\cL\rightarrow z^w\cO_{\dC\times
   X}\textup{ is non-degenerate}, \\ \\
 & & \;\;\cL\subset p^*V^{\alpha_1}\textup{ as subsheaves of } i_*i^*\cL =i_*\varphi^{-1}((p')^*H')\Big\}
\end{array}
$$
here $p:\dC\times X\rightarrow \dC$ resp. $p':\dC^*\times X\rightarrow
\dC^*$ are the projections and
$i:\dC^*\times X\hookrightarrow \dC\times X$ is the inclusion.
For any morphism $f:Y\rightarrow X$ of complex spaces and any
element $\cL\in\MBLfun(X)$ we set
$\MBLfun(f)(\cL):=f^*(\cL,\varphi)$.
\end{definition}

One easily checks that this is a functor which has the property
of being a sheaf for the classical topology. Notice also
that by definition, for any $X$, and any $(\cL,\varphi)\in\MBLfun(X)$,
the sheaf $\cL$ underlies a family of TERP-structures on $X$ in the sense of
definition \ref{defFamilyTERP}.
In order to study this functor, we will compare it to some other functor
which is easier to understand as it is simply a closed
subfunctor of some Grassmannian. We define it in several steps which corresponds to
the conditions imposed on the elements of $\MBLfun(X)$.

\begin{deflemma}
Consider the fixed data $H^\infty_\dR$, $M$, $S$, $w$, $\alpha_1$ from above as well
as the various Deligne extensions and the associated (finite dimensional) spaces $W^\alpha_\beta$.
Then $W^{\alpha_1}_{\alpha_\mu-1}$ is a symplectic vector space with respect to
the (class of the) anti-symmetric form $P^{(w-1)}$, we denote this
form by $\omega$ and write $W^\omega:=W^{\alpha_1}_{\alpha_\mu-1}$ for
short. In particular, the dimension $\dim_\dC(W^\omega)$
is even, and denoted by $2m$. Moreover, define nilpotent endomorphisms
$$
\left.
\begin{array}{rcl}
b& := & [z\cdot] \\
a& := & [z^2\nabla_z ] \\
\end{array}
\right\}
\in \End_\dC(W^\omega)=\End_\dC(V^{\alpha_1}/V^{>\alpha_\mu-1}).
$$
For any complex space $X$, consider the
$\cO_X$-locally free sheaf $\cW_X^\omega:=W^\omega\otimes_\dC \cO_X$ of rank $2m$.
$(\cW_X^\omega,\omega)$ is a symplectic bundle over $X$, with $\cO_X$-linear
operators $a$ and $b$.

Let $\cG(m,W^\omega)$ be the usual Grassmaniann functor, seen
as defined on the category of complex spaces. More precisely,
for any complex space $X$, let $\cG(m,W^\omega)(X)$ be the
set of rank $m$ locally free subsheaves
$\cG$ of $\cW^\omega_X$
that are locally direct summands. Clearly, if $f:Y\rightarrow X$ is a
morphism, then $\cG(m,W^\omega)(f)(\cG):=f^*\cG$.
Consider the following closed subfunctors:
$$
\begin{array}{rcl}
\cL\cG(W^\omega)(X)&:=&\left\{\cG\in\cG(m,W^\omega)(X)\,|\,\omega_{|\cG}=0
\right\}, \\ \\
\cL\cG_b(W^\omega)(X)&:=&\left\{\cG\in\cL\cG(W^\omega)(X)\,|\,b(\cG)\subset\cG
\right\},\\ \\
\cL\cG_{a,b}(W^\omega)(X)&:=&\left\{\cG\in \cL\cG_b(W^\omega)(X)\,|\,a(\cG)\subset\cG
\right\}
\end{array}
$$
These functors are represented by complex
spaces $G(m,W^\omega)$, $\Lambda(W^\omega)$, $\LL$
and $\Lambda_{a,b}(W^\omega)$, respectively. The first two of these
spaces are complex homogeneous, in particular, smooth, and all of them
have the structure of a projective variety.

% $$
% \begin{array}{rcl}
% \LGRabfun:(\mathit{ComplexSpaces}) & \longrightarrow & (\mathit{Sets}) \\ \\
% X & \longmapsto & \textup{rank $n$ free summands $ of } W^\omega\otimes_\dC \cO_X \\ \\
% f:Y\rightarrow X  & \longmapsto & f^(...)}
% \end{array}
% $$

% We define further
% $$
% \begin{array}{rcl}
% \Lambda(W^\omega) & := & \left\{L\subset W^\omega\,|\,\omega_{|L}=0\textup{ and }\dim_\dC(L)=\frac12\dim(W^\omega)\right\}, \\ \\
% \widetilde{\Lambda}_b(W^\omega) & := & \left\{L\in\Lambda(W^\omega)\,|\,b(L)\subset L\right\}, \\ \\
% \widetilde{\cM}_{BL} & := & \left\{L\in\Lambda(W^\omega)\,|\,b(L)\subset L\textup{ and } a(L)\subset L\right\}.
% \end{array}
% $$
% Then $\Lambda(W^\omega)$ is a complex homogeneous space and
% $\Lambda(W^\omega), \widetilde{\Lambda}_b(W^\omega)$ and $\widetilde{\cM}_{BL}$ are projective varieties.

\end{deflemma}
\begin{proof}
That $\omega\in(\bigwedge^2 W^\omega)^*$ defines a symplectic structure
is an immediate consequence of formula
\eqref{eqPairingPCoupling2} by putting $\alpha:=\alpha_1$,
$\beta:=\alpha_\mu-1$ and $k=w-1$. Moreover, both $z\cdot$ and $z^2\nabla_z$
map $V^\alpha$ to $V^{\alpha+1}$ and especially leave $V^{>\alpha_\mu-1}$
invariant, therefore $b$ and $a$ are well-defined and nilpotent.

Concerning the second part, first notice that again all of these functors
are sheaves of sets for the classical topology, hence, the universal
property needs only to be checked locally.
That the classical Grassmannian $G(m,W^\omega)$ represents the functor $\cG(m,W^\omega)$ is
well known. The subspace
$\Lambda(W^\omega):=\{L\subset W^\omega\,|\,\omega_{|L}=0\textup{ and }\dim_\dC(L)=m\}$
(sometimes called Lagrangian Grassmannian) represents $\cL\cG(W^\omega)$ and it is known that
it is complex homogeneous (in particular, smooth projective). Finally,
for any vector space $Y$ and any endomorphism $A\in\End_\dC(Y)$,
the subspace in $G(l,Y)$ of $A$-invariant $l$-dimensional subspaces
is easily seen to be closed, hence, the spaces $\LL
:=\{L\in\Lambda(W^\omega)\,|\,b(L)\subset L\}$ resp.
$\Lambda_{a,b}(W^\omega)
:=\{L\in\LL\,|\,a(L)\subset L\}$
are closed subvarieties of $\Lambda(W^\omega)$ and represent
$\cL\cG_b(W^\omega)$ resp. $\cL\cG_{a,b}(W^\omega)$.
\end{proof}

% Computation of examples shows that the spaces $\widetilde{\Lambda}_b(W^\omega)$ and $\widetilde{\cM}_{BL}$
% are non-reduced in general. As we are going to consider hermitian metrics
% on them, we are not interested in non-reduced parts of these spaces, therefore,
% define $\MBL$ resp. $\LL$ to be the complex analytic spaces obtained from
% $\widetilde{\cM}_{BL}$ resp. $\widetilde{\Lambda}_b(W^\omega)$ by taking the reduced structure sheaf.

In order to make use of these simplified functors, we have to compare them
to $\MBLfun$ which is our primary object of interest. This is done by the
following theorem.

\begin{theorem}\label{theoIsomorphismFunctors}
The natural transformation $\Phi:\MBLfun\longrightarrow \cL\cG_{a,b}(W^\omega)$ which
sends $(\cL,\varphi)\in\MBLfun(X)$ to $\cL/p^*V^{>\alpha_\mu-1}\in\cL\cG_{a,b}(W^\omega)(X)$ is
an isomorphism of functors. Hence $\MBLfun$ is represented by
$\Lambda_{a,b}(W^\omega)$, which we denote by $\MBL$.
\end{theorem}
\begin{proof}
First let us check that $\Phi$ is indeed well defined: We have that
$p^*V^{>\alpha_\mu-1}$ is a subsheaf of $\cL$ by lemma \ref{lemRegSingBounded}.
Moreover, $\cL/p^*V^{>\alpha_\mu-1}\subset p^*(V^{\alpha_1}/V^{>\alpha_\mu-1})$ by definition,
and the latter sheaf is isomorphic to
%$\cO_{\dC\times X}\otimes_{\cO_\dC}p^{-1}(V^{\alpha_1}/V^{>\alpha_\mu-1})
%\cong\cO_X\otimes_\dC W^\omega$ as $V^{\alpha_1}/V^{>\alpha_\mu-1}$
%is $z$-torsion.
$\cW^\omega_X$.
For any $x\in X$, the proof of lemma \ref{lemRegSingBounded} shows
also
$$
\cL_{|\dC\times \{x\}}=\{\sigma\in V^{>-\infty}\ |\
P^{(w-1)}(\cL_{|\dC\times \{x\}},\sigma)=0\}.
$$
Therefore $\cL_{|\dC\times \{x\}}/V^{>\alpha_\mu-1}
\subset V^{\alpha_1}/V^{>\alpha_\mu-1}$ is a lagrangian subspace.
With the lemma of Nakayama one obtains that
$\cL/p^*V^{>\alpha_\mu-1}$ is a rank $m$ locally free subsheaf of $\cW^\omega_X$
and is locally a direct summand. Because of $P^{(w-1)}(\cL,\cL)=0$,
it is an element of $\cL\cG(W^\omega)(X)$.
Finally, the $b$ resp. $a$-invariance follows directly
from the fact that $\cL$ is an $\cO_\dC$-module resp. from $(z^2\varphi^*\nabla_z)\cL\subset\cL$.

In order to show that $\Phi$ is an isomorphism, let us define an inverse.
For any complex space $X$, write $\pi:p^*V^{\alpha_1}\twoheadrightarrow
p^*(V^{\alpha_1}/V^{>\alpha_\mu-1})$ for the projection. Let
$\cG\in\cL\cG_{a,b}(W^\omega)(X)$ be given, write $k:X\hookrightarrow \dC\times X$ and consider
$\cL:=\pi^{-1}(k_*\cG)$. Then $i^*\cL\cong i^*p^*V^{\alpha_1} =(p')^*H'$, and this
defines the isomorphism $\varphi:i^*\cL\rightarrow (p')^* H'$. Put $\Psi(\cG):=(\cL,\varphi)$.
We have to show that this gives an element in $\MBLfun(X)$.

All the properties to be shown are local, hence we can restrict to the
case where $X=(X,x)$ is a germ of a complex space.
First, $\cL_{|(\dC,0)\times \{x\}}$ is $\dC\{z\}$-free of rank $\mu$,
and a basis $v_1^x,...,v_\mu^x$ of it is also a $\dC\{z\}[z^{-1}]$-basis
of $(V^{>-\infty})_0$. Furthermore, one can choose $m$ elements
%$\sigma_j^x=\sum_{i=1}^\mu \kappa_{ji}(z)\cdot v_i^x
%\in \bigoplus_{i=1}^\mu\dC\{z\}\cdot v_i^x$
$\sigma_j^x\in \cL_{|(\dC,0)\times \{x\}}$ ($j=1,...,m$)
such that they represent a basis of
$\cG_{|\{x\}}=\cL_{|\dC\times \{x\}}/V^{>\alpha_\mu-1}$.
Any sections $v_1,...,v_\mu\in \cL_0$ with ${v_i}_{|(\dC,0)\times \{x\}}=v_i^x$
are an $\cO_{\dC\times X,(0,x)}[z^{-1}]$-basis of $(p^*V^{>-\infty})_0$
by the lemma of Nakayama. Therefore they generate a free
$\cO_{\dC\times X,(0,x)}$-module of rank $\mu$ called $\cL'_0$.
Obviously $\cL'_0\subset \cL_0$.
In order to see the inverse inclusion $\cL_0\subset \cL'_0$, consider
$m$ sections $\sigma_j\in\cL'_0$ which extend the $\sigma_j^x$.
By the lemma of Nakayama they generate a rank $m$ free submodule of
$\cW^\omega_{X,x}$ which is a direct summand, and which is contained
in $\cG$. Therefore it coincides with $\cG$, and thus $\cL'_0\supset \cL_0$.
This shows that $\cL$ gives a vector bundle on $\dC\times (X,x)$.

The $a$-invariance of $\cG$ translates into the fact
that the connection $\varphi^*\nabla$ has a pole of order at most two along
$\{0\}\times X$ on $\cL$.
What remains to be shown
is that $\varphi^*P$ has the correct pole order properties on $\cL$.
This will complete the proof, as $\Psi$ is obviously an inverse for $\Phi$.
Consider $\varphi^*P$ as a pairing
$$
\varphi^*P=\sum_{k\in\dZ}P^{(k)}z^k:p^*V^{>-\infty}\otimes j^*p^*V^{>-\infty} \longrightarrow \cO_{\dC\times X}[z^{-1}].
$$
This induces a pairing $P:\cL\otimes j^*\cL\rightarrow \cO_{\dC\times X}[z^{-1}]$. We have to show that
$P^{(w-k)}(\cL,\cL)=0$ for all $k>0$ and that $P^{(w)}$ induces a non-degenerate
pairing $[P^{(w)}]:\cL/z\cL\otimes\cL/z\cL\rightarrow \cO_X$.
For the first point, notice that we have $P^{(w-1)}(\cL,\cL)=0$ by construction,
as $\omega_{|\cG}=0$. Moreover, the linearity of $P$ implies
$P^{(w-k)}(a,b)=P^{w-1}(z^{k-1}a,b)$ for any two sections $a,b\in\cL$ and
$k\in \dZ$. This gives the vanishing of $P^{(w-k)}$ on $\cL$ for $k>0$.

For the second point, consider the space $W^{\alpha_1-1}_{\alpha_\mu}$.
Again by formula \eqref{eqPairingPCoupling2}, we obtain a symplectic form $\omega'$ induced
by the class of $P^{(w-1)}$ on $W^{\alpha_1-1}_{\alpha_\mu}$. The subspace
$G':=\cG_{|\{x\}} + W^{\alpha_\mu-1}_{\alpha_\mu}=(\cL/\cV^{>\alpha_\mu})_{|\{x\}}
\subset W^{\alpha_1-1}_{\alpha_\mu}$
is again lagrangian with respect to $\omega'$. Now suppose that there is
$a\in \cL_{|\dC\times\{x\}}\backslash z\cL_{|\dC\times\{x\}}$
such that $P^{(w)}(a,\cL_{|\dC\times\{x\}})=0$.
Formula \eqref{eqPairingPCoupling2} shows that there exists some
$\widetilde{a}\in W^{\alpha_1}_{\alpha_\mu}$
such that $a-\widetilde{a}\in V^{>\alpha_\mu}\subset z\cL_{|\dC\times\{x\}}$
and $P^{(w)}(\widetilde{a},L)=0$.
This implies
$z^{-1}\www a\in z^{-1}\cL_{|\dC\times\{x\}}\backslash \cL_{|\dC\times\{x\}}$
and $\omega'(z^{-1}\widetilde{a}, G')=0$. The first property
gives $z^{-1}\www a\notin G'$, the second property and
the maximal isotropy of $G'$ imply $z^{-1}\www a\in G'$.
This is a contradiction. Therefore $P^{(w)}$ is nondegenerate on
$\cL_{|\dC\times\{x\}}/z\cL_{|\dC\times\{x\}}$ and thus also on
$\cL/z\cL$.

\end{proof}
As a piece of notation, for any complex space $X$, we write $|X|$ for the underlying topological space,
so that $X=(|X|,\cO_X)$ as ringed spaces. We will use this in particular for $X=\MBL$,
notice that it can happen that $\MBL$ has a non-reduced structure, as shown by the first example in subsection \ref{subsecWPS}.
We also write $\cL\in\VB_{\dC\times\MBL}$ for the universal sheaf of TERP-structures on $\MBL$,
i.e., $(\cL,\varphi)\in\MBLfun(\MBL)$ is the image of $\id_{\MBL}$ under
the isomorphism $\mathit{Hom}_{\mathit{CplxSp}}(-,\MBL) \rightarrow \MBLfun$.
By the above construction, this universal sheaf is explicitly given
as $\cL=\pi^{-1}(k_*\cG)$, where $\cG$ is the restriction of the tautological bundle
on the Grassmannian $G(m,W^\omega)$ to the closed subspace $\MBL$
and $\pi$ and $k$ are as above for the special case $X=\MBL$.
In section \ref{secHermMetricsGroupAction}, we will also need to consider
the space $\LL$, then the same construction yields
a universal sheaf $\cL'\in\VB_{\dC\times\LL}$ which has all properties
of a family of TERP-structures, except that the connection operator $\nabla_z$ may have
a pole order higher than two along $\{0\}\times \LL$.

The next lemma shows a case in which $\MBL$ has particularly simple structure.
\begin{lemma}
\label{lemMBLsmooth}
Suppose that $n=\lfloor\alpha_\mu-\alpha_1\rfloor=1$.
Then we have $\MBL\cong\LL\cong\Lambda(W^\omega)$, $\MBL$ is smooth in this case.
\end{lemma}
\begin{proof}
It follows from $zV^{>\alpha_\mu-2} = V^{>\alpha_\mu-1}$
and $(z^2\nabla_z)V^{>\alpha_\mu-2}\subset V^{>\alpha_\mu-1}$ that
both $b$ and $a$ vanish on $V^{>\alpha_\mu-2}/V^{>\alpha_\mu-1}$.
The condition of the lemma implies that $V^{\alpha_1}\subset V^{>\alpha_\mu-2}$
so that $b$ and $a$ vanish on $W^\omega$, which implies $\MBL=\Lambda(W^\omega)$.
\end{proof}

The next step is to understand the subspace of
$\MBL$ of regular singular TERP-structures with fixed spectral pairs
$\Spp=\sum_{\alpha,l}d(\alpha,l)\cdot (\alpha,l)\in\dZ[\dC\times \dZ]$.
We will define a subfunctor of $\MBLfun$ of such families, and
we will prove that it is represented by some complex subspace $U_{\Spp}$ of $\MBL$.
It turns out that $U_{\Spp}$ is locally closed in $\MBL$, and the spaces
$|U_{\Spp}|$ form
a stratification of $\MBLr$. Some of these strata, but usually
not all of them, are the classifying spaces $\cbl$ from \cite[section 2]{HS2}.

The definition of the spectral pairs $Spp(x)=\sum_{\alpha,l}d(\alpha,l)(x)$
for the regular singular TERP-structure $\cL_{|\dC\times\{x\}}$ for $x\in \MBLr$
in definition \ref{defHodgeSpectrum} can be rephrased as
$$
d(\alpha,l)(x)=\dim_\dC\frac{\Gr^\alpha_V\cL_{|x}\cap W_{l-(w-1)}W^\alpha_\alpha}
{\Gr^\alpha_V z\cL_{|x}\cap W_{l-(w-1)}W^\alpha_\alpha +
\Gr^\alpha_V\cL_{|x}\cap W_{l-1-(w-1)}W^\alpha_\alpha}.
$$
Here $W_\bullet $ is the weight filtration of the nilpotent endomorphism $N_z$ (the
logarithm of the unipotent part of the fixed automorphism $M\in\Aut(H^\infty_\dR)$) acting on $W^\alpha_\alpha$, centered
at $0$.

We first define the functor $U_{\Spp}$ alluded to above. For any complex space $X$, we will use, as
in definition \ref{defFuncMBL}, the pullbacks under
$p:\dC\times X\rightarrow \dC$ resp. $p':\dC^*\times X\rightarrow \dC$
of the flat bundle $H'$ and of the Deligne extensions $V^\alpha$ and $V^{>\alpha}$.
For simplicity of the notations, we write $\cV^\alpha:=p^* V^\alpha$ and $\cV^{>\alpha}:=p^* V^{>\alpha}$.
\begin{definition}\label{t7.7}
Fix the data $H^\infty_\dR,S,M,w$ as in definition \ref{defFuncMBL} and fix
moreover a tuple $\Spp$ of spectral pairs such that there is $x\in \MBLr$ with
$\Spp(x)=\Spp$.
%$\emptyset\neq |U_{\Spp}|\subset \MBLr$.
The functor $\cU^{H^\infty_\dR,S,M,w,\alpha_1}_{\Spp}$ ($\cU_{\Spp}$ for short)
from the category of complex spaces to the category of sets is defined by
\begin{eqnarray*}
\cU_{\Spp}(X)&:=& \Bigg\{(\cG,\phi)\in \MBLfun(X)\, |\,
\Spp(x)=\Spp\;\;\forall\ x\in X,\\ \\
&& \quad \forall\ (\alpha,l): \ \cV^\alpha\cap\cG, \cV^{>\alpha}\cap \cG\in \VB_{\dC\times X},
\ \Gr^\alpha_\cV \cG\in \VB_X,\
 \Gr^\alpha_\cV\cG\cap W_{l-(w-1)}W^\alpha_\alpha\in \VB_X,\\ \\
&&\quad \frac{\Gr^\alpha_\cV\cL \cap W_{l-(w-1)}W^\alpha_\alpha}
{\Gr^\alpha_\cV z\cL \cap W_{l-(w-1)}W^\alpha_\alpha +
\Gr^\alpha_\cV\cL \cap W_{l-1-(w-1)}W^\alpha_\alpha}\in \VB_X
\textup{ and its rank is }d(\alpha,l)\Bigg\}.
\end{eqnarray*}
\end{definition}

We have the following (not so surprising) statements about the functor $\cU_{\Spp}$.
\begin{theorem}\label{t7.6}
Each stratum $|U_{\Spp}|$ is the set of points of a complex space $U_{\Spp}$
with the following properties. It represents the functor $\cU_{\Spp}$ of families of
TERP-structures with constant spectral pairs in definition \ref{t7.7},
it carries a universal family, and the canonical map $U_{\Spp}\to \MBL$
is a locally closed embedding.
\end{theorem}
\begin{proof}[Proof of theorem \ref{t7.6}.]
As for Grassmannians, $\MBL$ can be covered by open affine subspaces,
each related to a TERP-structure with a fixed basis, and each consisting of
TERP-structures with bases obtained by deforming the given basis.
In fact, the choice of a basis could be reduced to the choice of an opposite filtration
in $H^\infty$, but here we stick to the more explicit bases.

We will first construct an affine chart of $\MBL$ together with its universal family.
This will be done without using the correspondence of lemma \ref{theoIsomorphismFunctors},
in other words, we will directly describe the universal family of TERP-structures.
Then an affine chart of $U_{\Spp}$ can easily be described as closed subspace
of this affine chart of $\MBL$, and hence we obtain a universal family on $U_{\Spp}$.

Write $\Spp=\sum_{j=1}^\mu (\beta_j,l_j)$ with $\alpha_1\leq \beta_1\leq ...\leq\beta_\mu\leq\alpha_\mu$.
Choose $x\in \MBL$ with $\Spp(x)=\Spp$. Choose an $M_s$-invariant common splitting of
$F^\bullet H^\infty$ and $W_\bullet H^\infty$. Choose a basis $b_1,...,b_\mu$
of $H^\infty$ which respects this splitting and such that $b_j$ corresponds to
the spectral pair $(\beta_j,l_j)$. This means that
$$
s_j := z^{\beta_j-\frac{N_z}{2\pi i}}b_j
\in \Gr^{\beta_j}_V\cL_{|x}\cap W_{l_j-(w-1)}W^{\beta_j}_{\beta_j}
$$
and that the classes of those $s_j$ with $(\beta_j,l_j)=(\alpha,l)$ form a basis of
$$
\frac{\Gr^\alpha_V\cL_{|x}\cap W_{l-(w-1)}W^\alpha_\alpha}
{\Gr^\alpha_V z\cL_{|x}\cap W_{l-(w-1)}W^\alpha_\alpha +
\Gr^\alpha_V\cL_{|x}\cap W_{l-1-(w-1)}W^\alpha_\alpha}.
$$
Now we define (finitely many) variables $c_{ij}^{(p)}$ for $i,j\in\{1,...,\mu\}$
and $p\in\dN$ with $\beta_j-p\geq \alpha_1$ and make the Ansatz
$$
v_i = s_i +\sum_{j,p}c_{ij}^{(p)}\cdot z^{-p}s_j.
$$
The requirement is that $(v_1,...,v_\mu)$ shall be a basis of the restriction of the universal family $\cL$
of $\MBL$ on a yet to be determined affine chart of $\MBL$. This chart will be defined
by the (analytic) spectrum of the quotient of $\dC[c_{ij}^{(p)}]$ by the ideal generated
by polynomial equations between the $c_{ij}^{(p)}$ which are defined by the
properties of the pairing $P$ and pole of order at most $2$.

The pairing $P$ gives the equations
$$
0=P^{(k)}(v_i,v_j)\qquad \textup{ for }
\left\lceil 2\alpha_1\right\rceil\leq k\leq w-1
$$
($P^{(k)}(v_i,v_j)=0$ anyway for $k<\left\lceil 2\alpha_1\right\rceil$
because of \eqref{eqPairingPCoupling1}).

Remember that $z^2\nnn_z s_j =\beta_j\cdot zs_j+\frac{-N_z}{2\pi i}zs_j$ and
$W^{\alpha_j+1}_{\alpha_j+1}=\bigoplus_{k,p:p\in \dZ,\alpha_k+p=\alpha_j+1}
\dC\cdot z^p s_k$.
Using this, one finds unique coefficients $\gamma_{ij}^{(q)}(q\geq 0),
\delta_{ij}^{(r)}(r\geq 1)\in \dC[c_{ij}^{(p)}]$ with
$$
z^2\nnn_z v_i -\sum_{j,q}\gamma_{ij}^{(q)}\cdot z^qv_j =
\sum_{j,r}\delta_{ij}^{(r)}\cdot z^{-r}s_j.
$$
This gives the equations
$$
\delta_{ij}^{(r)}=0.
$$
Now the affine chart of $\MBL$ we are looking for
is $\Specan(\dC[c_{ij}^{(p)}]/(P^{(k)}(v_i,v_j),\delta_{ij}^{(r)}))$
and $v_1,...,v_\mu$ form a basis of the universal family $\cL$ of
TERP-structures on this chart.

Of course, $x$ is in this chart, and the numbers $c_{ij}^{(p)}(x)$ satisfy
$c_{ij}^{(p)}(x)=0$ for $\alpha_j-p\leq \alpha_i$, i.e.
${v_i}_{|\dC\times \{x\}}-s_i\in V^{>\alpha_i}$.

The subfamily of TERP-structures with spectral pairs equal to $\Spp$
is simply obtained by the additional equations
\begin{eqnarray*}
&&c_{ij}^{(p)}=0\qquad \textup{ for }\alpha_j-p<\alpha_i,\\
&&c_{ij}^{(p)}=0\qquad \textup{ for }\alpha_j-p =\alpha_i\textup{ and }l_j>l_i.
\end{eqnarray*}
This defines an affine closed subspace of the affine chart of $\MBL$.

One obtains a locally closed
(but most often not closed) subspace $U_{\Spp}$ of $\MBL$.
The restriction of $\cL$ to $U_{\Spp}$  is the universal
family of TERP-structures with fixed spectral pairs, and hence
$U_{\Spp}$ represents the functor $\cU_{\Spp}$.
\end{proof}

The next result is more surprising than theorem \ref{t7.6}, and illustrates
that the spaces $\cbl$ behave better than an arbitrary stratum $U_{\Spp}$.
\begin{theorem}\label{t7.8}
Suppose that the spectral pairs of a stratum $U_{\Spp}$
coincide with the spectral pairs of a certain classifying space $\cbl$.
Then $U_{\Spp}=\cbl$ as complex spaces, in particular, $U_{\Spp}$ is reduced
and smooth in this case.
\end{theorem}

\begin{proof}
The first part of the proof shows $|U_{\Spp}|=|\cbl|$, the second part
discusses the complex structures.

\textbf{First part:} Suppose that we are given two regular singular TERP-structures with the same topological
data and the same spectral pairs $\Spp$, where one is a reference TERP structure
whose filtration $\widetilde{F}^\bullet_0$ (see definition \ref{defHodgeSpectrum}) is
part of a PMHS, whereas the other one is arbitrary and induces the filtration $\widetilde{F}^\bullet$.
We have to show that the second TERP-structure is an element of the classifying space
$\cbl$ defined by the first one.

{}From the construction of $\cbl$ as a bundle over $\DcPMHS$
\cite{He2} it follows that it is sufficient to show that $\widetilde{F}^\bullet$
lies in the classifying space $\DcPMHS$ which contains
$\widetilde{F}^\bullet_0$.
The definition of $\DcPMHS$ was rewritten in \cite{HS2} lemma 2.5 (i).
We refer to the notations used in loc.cit.
That the conditions $N(\widetilde{F}^p)\subset \widetilde{F}^{p-1}$ and those concerning $S$
(for $H^\infty_{\neq 1}$ with $w-1$ instead of $w$, for $H^\infty_1$ with $w$) hold is
clear from the construction of $\widetilde F^\bullet$, see
\cite[definition 2.3]{HS2}. It remains to show that for any eigenspace $H^\infty_\lambda$,
the conditions concerning $N,\ \widetilde{F}^\bullet$
and the primitive part $P_l\subset \Gr^W_l$ hold, namely,
\begin{eqnarray*}
\dim \widetilde{F}^pP_l=\dim \widetilde{F}^p_0P_l,\qquad \widetilde{F}^p N^j P_l=N^j \widetilde{F}^{p+j}P_l,\qquad
\widetilde{F}^p \Gr^W_l=\bigoplus_{j\geq 0} \widetilde{F}^pN^jP_{l+2j}.
\end{eqnarray*}
Notice that these conditions are closely related to the strictness of the powers of $N$.

Recall that the spectral pairs consist of finitely many sequences
of pairs $(\alpha,w-1+k),(\alpha-1,w-1+k-2),\ldots,(\alpha-k,w-1-k)$.
Each sequence corresponds to one Jordan block of $M$.
We can read off
the dimensions $\dim \Gr^W_l$ and $\dim P_l$
from the second entries of the spectral pairs only.

To show the above conditions, we will work inductively.
We fix $\lambda$, i.e. we consider only $\alpha$ with $e^{-2\pi i\alpha}=\lambda$.
Consider the set of sequences where $-2\cdot \alpha+(w-1+k)$ is minimal and choose
a sequence in this set where $k$ is maximal.
Then
$$
0 \neq \widetilde{F}^{\lfloor w-\alpha \rfloor}\Gr^W_k H^\infty_\lambda
\quad\quad\mbox{and}\quad\quad\quad
0 = \widetilde{F}^{\lfloor w-\alpha+1-j \rfloor}\Gr^W_{k-2j} \;\;\textup{ for all }j\in\dZ.
$$
Choose an element $v_1\in \widetilde{F}^{\lfloor w-\alpha \rfloor}\Gr^W_k H^\infty_\lambda \backslash \{0\}$.
Then also $N^j(v_1)\in \Gr^W_{k-2j}$ are non-vanishing for $j=0,1,\ldots,k$.
They must correspond to some spectral pairs where the second indices are
$w-1+k-2j$. But we know more, they actually correspond
to the spectral pairs
in the sequence which starts with $(\alpha,w-1+k)$, this follows from
$N(\widetilde{F}^p)\subset \widetilde{F}^{p-1}$ and from the vanishing
property from above. This gives strictness of the powers of $N$ with
respect to $v_1$.

Now by dividing out the subspaces $\dC\cdot N^j(v_1)$ from the
spaces $\Gr^W_{k-2j}$ we can erase
this sequence of spectral pairs and consider the next one.
We obtain in the quotients by the subspaces $\dC\cdot N^j(v_1)$ an element
$\widetilde{v}_2$ and images $N^j(\widetilde{v}_2)$ which correspond
to the next sequence of spectral pairs. It is possible to lift $\widetilde{v}_2$ and
all $N^j \widetilde{v}_2$ to elements $v_2$ and $N^j(v_2)$ such that they still
correspond to the next sequence of spectral pairs. Here one has possibly to adjust
a first lift of $\widetilde{v}_2$ by a multiple of a suitable image $N^j(v_1)$.
By induction, the strictness of all powers of $N$ and the
equalities above concerning $N,\ P_l$ and $\widetilde{F}^\bullet$ can be proved.

\medskip
\textbf{Second part:}
We will show that the main work has actually already been done in \cite[chapter 5]{He2}
and \cite[proof of theorem 12.8]{He3}.
In \cite[chapter 5]{He2}, the fibers of the bundle $\pi_{BL}:\cbl\to\DcPMHS$ are analyzed.
This is done for Brieskorn lattices, but these are Fourier-Laplace dual to
TERP-structures, and everything in loc.cit. can easily be rewritten with
TERP-structures. In this language \cite[chapter 5]{He2} does the following.

For a fixed filtration $\www F^\bullet_0\in\DcPMHS$ families of TERP-structures
inducing this filtration considered. An Ansatz is made as in the proof of theorem \ref{t7.6}
above to construct a universal family. As $\www F^\bullet_0\in \DcPMHS$,
the sections $s_j$ can be chosen to have very good properties with respect to
$P$, $\www F^\bullet_0$ and $N_z$. Then it is shown that the equations from $P$ and
$z^2\nnn_z$ yield a smooth complex space $\pi_{BL}^{-1}(\www F^\bullet_0)$ isomorphic
to $\dC^{N_1}$ for some $N_1\in\dN$. By construction it represents a functor
of families of TERP-structures with fixed filtration $\www F^\bullet_0$,
which is defined analogously to the functors $\cU_{\Spp}$.

In \cite[chapter 2]{He2} the spaces $\DcPMHS$ are constructed, but as homogeneous spaces,
not as spaces representing a functor of families of PMHS-like filtrations
with fixed spectral pairs. In particular, \cite[chapter 2]{He2} does not discuss
local coordinates for the smooth spaces $\DcPMHS$.
Though this is done in \cite[proof of theorem 12.8]{He3}.
It is easy to lift that discussion to a proof that $\DcPMHS$ with its smooth structure
represents a functor of families of PMHS-like filtrations with fixed spectral pairs.

One can combine the discussions in \cite[chapter 5]{He2}
and \cite[proof of theorem 12.8]{He3} in the following Ansatz similar to that in
the proof of theorem \ref{t7.6}, in order to show that $\cbl$ represents the
functor $\cU_{\Spp}$.

First choose $s_i\in W^{\beta_i}_{\beta_i}$ as in \cite[chapter 5]{He2},
fitting to $\www F^\bullet_0$. Make the Ansatz
$$
\www s_i = s_i +\sum_{j,p:p\geq 1, \beta_j-p=\beta_i} b_{ij}^{(p)}\cdot z^{-p}s_j
$$
for the family of filtrations and the Ansatz (generalizing that in \cite[chapter 5]{He2})
$$
v_i = \www s_i + \sum_{j,p:p\geq 1,\beta_j-p>\beta_i} c_{ij}^{(p)}\cdot z^{-p}\www s_j
$$
for a basis of sections of a family of TERP-structures.
The conditions for the $v_i$ from $P$ and $z^2\nnn_z$ contain the
analogous conditions for the $\www s_i$, and by \cite[lemma 5.7]{HS1} these
are equivalent to $N_z(\www F^p)\subset \www F^{p-1}$ and conditions
for $\www F^\bullet$ from the pairing $S$ from \cite[formula (5.1)]{HS1},
alluded to in definition \ref{defTopData}.
Therefore, combining the discussions in \cite[chapter 5]{He2}
and \cite[proof of theorem 12.8]{He3},
one obtains a smooth affine chart of $\cbl$
and a (restriction of a) universal family on it.
Consequently, $\cbl$ then represents the functor $\cU_{\Spp}$.
We leave the details to the reader.
\end{proof}

\begin{remark}\label{t7.9}
The closed embedding $\cbl=U_{\Spp}\hookrightarrow \MBL$ factors through
$\MBL^{red}$, as $\cbl$ is smooth.
For any component of $\MBL^{red}$, there exists a tuple $Spp$ of
spectral pairs such that $U_{\Spp}^{red}$ is the generic stratum.
Then $U_{\Spp}^{red}$ is open in (this component of) $\MBL^{red}$,
but $U_{\Spp}$ is not necessarily open in $\MBL$.
Subsection \ref{subsecWPS} gives an example where the generic stratum is a
space $\cbl$ and where $\MBL$ is not smooth on $|\cbl|\subset \MBLr$.
\end{remark}

%\clearpage

\section{Hermitian metrics and $G_\dZ$-action}
\label{secHermMetricsGroupAction}

In this section we define and study certain subspaces $\MBLp$ resp. $\LLp$ of $\MBL$ resp. $\LL$ which we
call pure polarized. This is reminiscent to the subspace $\cblp$ of $\cbl$
considered in \cite{HS2}. Using the twistor construction, we obtain positive definite
hermitian metrics on the tangent sheaves of these subspaces. The first main result
of this section is that the induced distances are complete. This is in sharp contrast to
the distance on the space $\cblp$ (see the example in the end of section 4 of \cite{HS2})
and motivates the construction of the compact classifying space
$\MBL$. In the second part, we study the action of a discrete group
on $\MBLp$ under the condition that there is an $M$-invariant lattice
in the vector space $H^\infty_\dR$ we started with. This will yield
quotients $\MBLp/G_\dZ$ and $(\overline{\check{D}}_{\mathit{BL}}\cap\MBLp)/G_\dZ$, which are
well suited as targets of period maps of variations of TERP-structures
over non-simply connected parameter spaces (see subsection \ref{subsecApplications})

Let us first give the definition of the pure polarized parts of $\MBL$ and $\LL$.
The universal locally free $\cO_{\dC\times\MBL}$-module $\cL$ underlies a family of TERP-structures
on $\MBL$. The construction in definition \ref{defExtInfinity} yields
an extension to a real-analytic family of holomorphic $\dP^1$-bundles, denoted by
$\widehat{\cL}$. Similarly, we obtain a locally free $\cO_{\dP^1}\cC^{an}_{\LL}$-module
$\widehat{\cL}'$ (remember that $\cL'$ was the universal sheaf on $\dC\times\LL$ constructed
in the same way as the sheaf $\cL$ on $\dC\times\MBL$).
For each $x\in\MBLr$ resp. $x\in\LLr$, the anti-linear involution
$\tau$ acts on $H^0(\dP^1,\widehat{\cL}_{|\dP^1\times\{x\}})$ resp.
$H^0(\dP^1,\widehat{\cL}'_{|\dP^1\times\{x\}})$. It induces a hermitian
form $h(-,-):=z^{-w}P(-,\tau-)$ on these spaces (which takes values in $\dC$).

\begin{definition}
\label{defMBLPurePolarized}
Define the following subspaces of $\MBLr$ resp. $\LLr$.
$$
\begin{array}{c}
|\MBL^{pure}| :=  \left\{x\in\MBLr\;\big|\;\widehat{\cL}_{|\dC\times\{x\}}\cong\cO_{\dP^1}^\mu\right\}
\quad,\quad
|\LL^{pure}| :=  \left\{x\in\LLr\;\big|\;\widehat{\cL}'_{|\dC\times\{x\}}\cong\cO_{\dP^1}^\mu\right\}, \\ \\
\MBLpr  :=  \left\{x\in|\MBL^{pure}|\;\big|\; h \textup{ is positive definite on }H^0(\dP^1,\widehat{\cL}_{|\dP^1\times\{x\}})\right\}, \\ \\
\LLpr  :=  \left\{x\in|\LL^{pure}|\;\big|\;h \textup{ is positive definite on }H^0(\dP^1,\widehat{\cL}'_{|\dP^1\times\{x\}})\right\},
\end{array}
$$
All of these subspaces are endowed with the canonical complex structures
defined by the restriction of $\cO_{\MBL}$ resp. $\cO_{\LL}$.
\end{definition}
\begin{lemma}\label{lemCharacterizationsPurePolarized}
$\MBL^{pure}$ and $\LL^{pure}$ are complements
of real-analytic subvarieties, the pure polarized parts
$\MBLp$ and $\LLp$ are unions of connected components of these complements,
and we have the following characterization of these subspaces.
\begin{eqnarray}
\label{eqHNonDegPure1}
|\MBL^{pure}| =  \left\{x\in\MBLr \;\big|\; h \textup{ is a non-degenerate on }H^0(\dP^1,\widehat{\cL}_{|\dP^1\times\{x\}})\right\}, \\
\nonumber \\
\label{eqHNonDegPure2}
|\LL^{pure}| =  \left\{x\in\LLr \;\big|\; h \textup{ is a non-degenerate on }H^0(\dP^1,\widehat{\cL}'_{|\dP^1\times\{x\}})\right\}, \\
\nonumber \\
\label{eqHNonDegPure3}
\MBLpr  =  \left\{x\in\MBLr \;\big|\; h \textup{ is positive definite on }H^0(\dP^1,\widehat{\cL}_{|\dP^1\times\{x\}})\right\}, \\
\nonumber \\
\label{eqHNonDegPure4}
\LLpr  =  \left\{x\in\LLr \;\big|\; h \textup{ is positive definite on }H^0(\dP^1,\widehat{\cL}'_{|\dP^1\times\{x\}})\right\}.
\end{eqnarray}
\end{lemma}
\begin{proof}
As $\widehat{\cL}$ resp $\widehat{\cL}'$ depend real-analytically on the
parameters, and as triviality of vector bundles is an open condition,
the first statement is clear. Equation \eqref{eqHNonDegPure1}
follows directly from lemma \ref{lemGeneralitiesTERP}, 5.,
and the same argument also gives equation \eqref{eqHNonDegPure2}. The remaining
equations \eqref{eqHNonDegPure3} and \eqref{eqHNonDegPure4} are then obvious.
\end{proof}
For fixed
spectral pairs $\Spp$, we recover the pure polarized part
of $\cbl$ as $\cblp =  \MBLp\cap \cbl$. We put
$\cL^{sp}:=p_*\widehat{\cL}_{|\dP^1\times\MBLp}\in\VBan_{\MBLp}$ resp.
$\cL'^{sp}:=p'_*\widehat{\cL'}_{|\dP^1\times\LLp}\in\VBan_{\LLp}$, by definition
these sheaves come equipped with positive definite
hermitian metrics defined by $h$.

Considering the tangent maps of the closed embeddings
$$
\MBL\hookrightarrow G(m,W^\omega)\quad
\textup{and}\quad
\LL\hookrightarrow G(m,W^\omega)
$$
gives inclusions
$$
%\begin{array}{c}
k_*\Theta_{\MBL}\subset k_*{\cH}\!om_{\cO_{\MBL}}(\cG,(\cO_{\MBL}\otimes W^\omega)/\cG) \cong {\cH}\!om_{\cO_{\dC\times\MBL}}(\cL,\cV^{\alpha_1}/\cL)
\subset {\cH}\!om_{\cO_{\dC\times\MBL}}(\cL,z^{-n}\cL/\cL)
$$
and similarly
$
%\quad
%\textup{resp.} \\ \\
(k')_*\Theta_{\LL}
%(k')_*{\cH}\!om_{\cO_{\LL}}(\cG',(\cO_{\LL}\otimes W^\omega)/\cG') \\ \\ \cong
%{\cH}\!om_{\cO_{\dC\times\LL}}(\cL',\cV^{\alpha_1}/\cL')
\subset {\cH}\!om_{\cO_{\dC\times\LL}}(\cL',z^{-n}\cL'/\cL'),
%\end{array}
$
where $k:\MBL\hookrightarrow \dC\times\MBL$, $x\mapsto(0,x)$ resp.
$k':\LL\hookrightarrow \dC\times\LL$, $x\mapsto(0,x)$.
Moreover, we have splittings
$k^{-1}(\cO_\dC\cAn_{\MBLp}\otimes \cL) =
\cL^{sp}\oplus k^{-1}(\cO_\dC\cAn_{\MBLp}(z\cL))$ resp.
$(k')^{-1}(\cO_\dC\cAn_{\LLp}\otimes \cL') =
\cL'^{sp}\oplus (k')^{-1}(\cO_\dC\cAn_{\LLp}(z\cL'))$
(these are equalities of $\cAn_{\MBLp}$ resp.
$\cAn_{\LLp}$-modules).
This yields
$$
\begin{array}{l}
\cAn_{\MBLp}\otimes \Theta_{\MBLp} \subset
{\cH}\!om_{\cAn_{\MBLp}}\left(\cL^{sp},\oplus_{i=1}^n z^{-i}\cL^{sp}\right)
\quad\quad\textup{resp.}\\ \\
\cAn_{\LLp}\otimes \Theta_{\LLp} \subset
{\cH}\!om_{\cAn_{\LLp}}\left(\cL'^{sp},\oplus_{i=1}^n z^{-i}\cL'^{sp}\right)
\end{array}
$$
which defines positive definite hermitian metrics (both denoted by $h$) on $\Theta_{\MBLp}$ resp. $\Theta_{\LLp}$.
Here positive definite means that for any point $x\in\MBLpr$ resp. $x\in\LLpr$, the induced metrics on the fibres
$$
\Theta_{\MBLp}/\mathbf{m}_x\Theta_{\MBLp}
\quad\textup{resp.}\quad
\Theta_{\LLp}/\mathbf{m}_x\Theta_{\LLp},
$$
i.e., on the Zariski tangent spaces of $\MBLp$ resp. $\LLp$ at $x$, are
positive definite. Consider the linear spaces $T_{\MBLp}$ resp.
$T_{\LLp}$ associated to $\Theta_{\MBLp}$ resp. $\Theta_{\LLp}$, which are the unions of
the Zariski tangent spaces at all points of $\MBLpr$ resp. $\LLpr$.
The hermitian metric $h$ defines length functions $l_h:T_{\MBLp}\rightarrow \dR_{\geq 0}$ resp.
$l_h:T_{\LLp}\rightarrow \dR_{\geq 0}$.
Following \cite[section 2.3]{Ko3}, these length functions
define by integration of piecewise $C^1$-curves
distance functions, both denoted by $d_h$ on $\MBLp$ resp. $\LLp$.
Notice that the fact that $\MBLp$ resp $\LLp$ may be non-reduced
does not affect this construction, as the linear spaces
$T_{(\MBLp)^{red}}$ resp. $T_{(\LLp)^{red}}$ are contained in
$T_{\MBLp}$ resp. $T_{\LLp}$.
By definition, the distances $d_h$ are inner distances (see cit.loc., proposition 1.1.8 and chapter
2.3), so that they induce the standard topology of $\MBLpr$ resp. $\LLpr$.
In particular, they are weakly complete, that is, for any point $x\in\MBLpr$
resp. $x\in \LLpr$, there is an $\varepsilon\in\dR_{>0}$ such that
the closed ball $B^{d_h}_\varepsilon(x)$ is compact in $\MBLpr$ resp. in $\LLpr$.

We are going to
show that the distance $d_h$ on $\LLpr$ is in fact strongly complete,
that is, that there is a uniform $\varepsilon$ with this property.
For this purpose, we will construct for any $L\in \LLpr$ a metric embedding of $\LLpr$ into a larger
space $|\Lambda^{pp}_b(W^\omega_L)|$ and show that all $|\Lambda^{pp}_b(W^\omega_L)|$ are isometric.
This will prove that $\LLpr$ is strongly complete. It follows that $\MBLpr$ is strongly complete,
as it is a closed subspace of $\LLpr$. By standard arguments (see \cite[proposition 1.1.9]{Ko3}),
a strongly complete space is Cauchy complete, i.e., any Cauchy sequence has a limit.

Consider, as before, the free $\dC[z,z^{-1}]$-module $W:=H^0(\dP^1,\widetilde{i}_* V^{>-\infty}
\cap i_*V_{<\infty})$. The anti-linear involution $\tau:\cH'\rightarrow \gamma^*\overline{\cH'}$ extends to
$\widetilde{i}_* V^{>-\infty} \cap  i_* V_{<\infty}$
and therefore defines an anti-linear automorphism of $W$. As we already remarked in the last section,
the pairing $P$ is defined on $W$. We obtain a pairing
$\widehat{S}_W:W\times W\rightarrow \dC[z,z^{-1}]$ by putting $\widehat{S}_W(-,-):=z^{-w}P(-,\tau-)$.
It satisfies $\widehat{S}_W(a,b)=-\widehat{S}_W(za,zb)=-\widehat{S}_W(z^{-1}a,z^{-1}b)$.
We write $h_W$ for the hermitian pairing $W\times W\rightarrow \dC$ defined by composing
$\widehat{S}_W$ with the natural projection $\dC[z,z^{-1}]\rightarrow \dC$ onto the $z^0$-component.
For any $L\in \LLpr$, we write $L^{sp}$
for the fibre of $(\cL')^{sp}$ at the point $L$. The purity of $\widehat{\cL}'$
on $\LLp$ gives that we have a decomposition $W = \oplus_{i\in \dZ} z^i L^{sp}$,
which is $h_W$-orthogonal (notice that $h$ equals $\widehat{S}_W$ and $h_W$
on $L^{sp}$).
\begin{deflemma}
For $L\in \LLpr$ put
$$
W^\omega_L:=
z^{-n}\dC[z]L^{sp} \cap z^{n-1} \dC[z^{-1}]L^{sp} = \bigoplus\limits_{k=-n}^{n-1} z^k L^{sp}
=z^{-n}\left(L\oplus W^{>\alpha_\mu-1}_{<\infty}\right)\cap z^{n-1}\tau\left(L\oplus W^{>\alpha_\mu-1}_{<\infty}\right)
\subset W.
$$
Then $W^\omega\subset W^\omega_L$.
$P^{(w-1)}$ is non-degenerate on $W^\omega_L$, and we write
as before $\Lambda(W^\omega_L)$ for the Lagrangian Grassmannian
of half-dimensional subspaces of $W^\omega_L$ on which
$\omega=[P^{(w-1)}]$ vanishes. Similarly to the situation considered before, we define
$$
\begin{array}{rcl}
|\Lambda_b(W^\omega_L)| & := & \left\{G\in|\Lambda(W^\omega_L)|\;\big|\;b(G)\subset G\right\},\\ \\
|\Lambda^{pp}_b(W^\omega_L)| & := & \left\{G\in|\Lambda_b(W^\omega_L)|\;\big|\;h_W\textup{ is positive definite on }
G^{sp} \right\}.
\end{array}
$$
Here
$
G^{sp}:=(G\oplus z^n\dC[z]L^{sp})\cap\tau(G\oplus z^n\dC[z]L^{sp})\subset W
$
for any $G\in|\Lambda_b(W^\omega_L)|$. These spaces are equipped with canonical complex structures,
as they are defined as subspaces of a Grassmannian. Again we have
the universal sheaves $\cG\in\VB_{\Lambda^{pp}_b(W^\omega_L)}$ (where
$\cG_{|G}=G$), $\cK:=(\pi^L)^{-1}(k^L_*\cG)\in\VB_{\dC\times\Lambda^{pp}_b(W^\omega_L)}$,
where $k^L:\Lambda^{pp}_b(W^\omega_L)\hookrightarrow \dC\times\Lambda^{pp}_b(W^\omega_L), x\mapsto(0,x)$,
$$
\pi^L:
z^{-n}\cO_{\dC\times\LLp}\otimes L^{sp} \twoheadrightarrow
\frac{z^{-n}\cO_{\dC\times\LLp}\otimes L^{sp}}{z^n\cO_{\dC\times\LLp}\otimes L^{sp}}
\cong k'_*\left(\bigoplus_{k=-n}^{n-1} \cO_{\LLp} \otimes z^k L^{sp}\right),
$$
and $\cK^{sp}\in\VBan_{\Lambda^{pp}_b(W^\omega_L)}$ (with $\cK^{sp}_{|G}=G^{sp}$).
The latter sheaf comes equipped with a positive definite
hermitian metric, which induces
a hermitian metric $h_L$ on the tangent sheaf $\Theta_{\Lambda^{pp}_b(W^\omega_L)}$.
We write $d_{h_L}$ for the induced distance function on $|\Lambda^{pp}_b(W^\omega_L)|$.
\end{deflemma}
\begin{proof}
The first statement simply follows from the fact that by construction,
we have $W^{\alpha_1}_{<\infty} \subset z^{-n}\dC[z]L^{sp}$
and, consequently,
$W^{>-\infty}_{\alpha_\mu-1}\subset \tau(z^{-n+1}\dC[z]L^{sp}) =z^{n-1} \dC[z^{-1}]L^{sp}$ for any $L\in\LLp$. Moreover,
$P^{(w)}$ is non-degenerate on $L^{sp}$, so that
$P^{(w-1)}:z^{-i}L^{sp}\times z^{i-1}L^{sp}\rightarrow \dC$ is non-degenerate and thus it
induces a symplectic form on $W^\omega_L$.

Notice that by the same argument as in the proof of
\ref{lemCharacterizationsPurePolarized} (that is, essentially by lemma \ref{lemGeneralitiesTERP}, 5.),
the points of $|\Lambda^{pp}_b(W^\omega_L)|$ parameterizes those $b$-invariant subspaces $G\in|\Lambda(W_L^\omega)|$
such that $G^{sp}$ defines an extension of $G\oplus z^n\dC[z]L^{sp}$ to a trivial (algebraic) bundle
over $\dP^1$, on which the pairing $h_W$ is positive definite.

The hermitian metric on the tangent sheaf of $\Lambda^{pp}_b(W^\omega_L)$ is defined as before:
{}From the definition of $\cK$ we know that $z^n\cO_{\dC\times\Lambda^{pp}_b(W^\omega)}\otimes L^{sp}\subset \cK$,
hence $z^{-n} \cO_{\dC\times\Lambda^{pp}_b(W^\omega)}\otimes L^{sp}\subset z^{-2n}\cK$
(as subsheaves of the $\cO_{\dC\times\Lambda^{pp}_b(W^\omega_L)}[z^{-1}]$-module $\cV^{>-\infty}$).
This yields an inclusion
$$
\cAn_{\Lambda^{pp}_b(W^\omega_L)} \otimes
\Theta_{\Lambda^{pp}_b(W^\omega_L)} \subset {\cH}\!om_{\cAn_{\Lambda^{pp}_b(W^\omega_L)}}
\left(\cK^{sp},\oplus_{i=1}^{2n} z^{-i}\cK^{sp}\right)
$$
which defines the hermitian metric on $\Theta_{\Lambda^{pp}_b(W^\omega_L)}$ (denoted by $h_L$)
by restriction.
\end{proof}
\begin{lemma}
\label{lemIsometry}
Let $L_1,L_2\in|\LLp|$. Then there is an isomorphism $A:L^{sp}_1\rightarrow L^{sp}_2$
which induces an isometry
$$
\Phi_A:(|\Lambda^{pp}_b(W^\omega_{L_1})|,d_{h_{L_1}}) \longrightarrow (|\Lambda^{pp}_b(W^\omega_{L_2})|,d_{h_{L_2}}).
$$
\end{lemma}
\begin{proof}
Choose bases $\underline{w}_1$ of $L^{sp}_1$ and $\underline{w}_2$ of $L^{sp}_2$ such that
$\tau(\underline{w}_i)=\underline{w}_i$ and $z^{-w}P(\underline{w}^{tr}_i,\underline{w}_i)=\de_\mu$ ($i=1,2$)
and define $A$ by putting $A(\underline{w}_1):=\underline{w}_2$. From
$W=\oplus_{i\in\dZ}z^iL^{sp}$ we see that $A$ can be extended
$z$-linearly to an automorphism of $W$ which respects both $P$ and
$\tau$ and thus also $\widehat{S}_W$ and $h_W$.
In particular, $A(W^\omega_{L_1}) \subset W^\omega_{L_2}$, and as $A^*\omega=\omega$ we obtain an
induced mapping $\Phi_A:\Lambda_b(W^\omega_{L_1})\rightarrow \Lambda_b(W^\omega_{L_2})$.
{}From $A^*h_W=h_W$ we conclude that $\Phi_A:\Lambda^{pp}_b(W^\omega_{L_1})\stackrel{\cong}{\rightarrow}
\Lambda^{pp}_b(W^\omega_{L_2})$, and the definition of the hermitian metrics of
these spaces gives that $\Phi_A$ is an isometry.
\end{proof}
\begin{lemma}
\label{lemEmbedMetricSpace}
For any $L\in\LLp$, there is a canonical closed embedding
$$
i_L:\LL \hookrightarrow \Lambda_b(W^\omega_L)
$$
which sends $\LLp$ to $\Lambda^{pp}_b(W^\omega_L)$.
Moreover $i_L^*h_L = h$ and consequently $i_L^*d_{h_L}\leq d_h$.
\end{lemma}
\begin{proof}
We define $i_L(G):=(G \oplus W_{<\infty}^{>\alpha_\mu-1})\cap W^\omega_L$. Then the image of $i_L$ is given as
$$
\textup{Im}(i_L)=
\left\{
\widetilde{G}\subset\Lambda_b(W^\omega_L)\,|\,
W^\omega_L\cap W_{<\infty}^{>\alpha_\mu-1} \subset \widetilde{G}
\right\}
$$
which shows that it is closed, an inverse map is given by $\widetilde{G}\mapsto G:=\widetilde{G}\cap W^\omega$.
That $i_L(\LLp)\subset\Lambda^{pp}_b(W^\omega_L)$ follows from the fact that
$(i_L(G))^{sp}=G^{sp}$.
Finally, for any $\xi\in\Theta_{\LLp}$, we have
$$
|\xi|_h = |(i_L)_*(\xi)|_{h_L}
$$
where $(i_L)_*:\Theta_{\LLp} \hookrightarrow i_L^{-1}\Theta_{\Lambda^{pp}_b(W^\omega_L)} \otimes \cO_{\LLp}$ is the tangent
map. This follows directly from the definition of the hermitian metrics $h$ and $h_L$: Consider the
diagram
$$
\xymatrix{
\cAn_{\LLp}\otimes\Theta_{\LLp} \ar@{^{(}->}[d] \ar@{^{(}->}[rrr]^{(i_L)_*}
                &  & & \cAn_{\LLp} \otimes i_L^{-1}\Theta_{\Lambda^{pp}_b(W^\omega_L)} \ar@{^{(}->}[d]  \\
\bigoplus\limits_{i=1}^n{\cH}\!om_{\cAn_{\LLp}}(\cL'^{sp}, z^{-i}\cL'^{sp})
 \ar@{^{(}->}[rrr]^g
                & & & \bigoplus\limits_{i=1}^{2n}\left(\cAn_{\LLp}\otimes i_L^{-1}{\cH}\!om_{\cAn_{\Lambda^{pp}_b(W^\omega_L)}}
                (\cK^{sp}, z^{-i}\cK^{sp})\right)
}
$$
then it follows from $i_L^{-1}\cK^{sp}=\cL'^{sp}$ that $g$ is simply defined by
$g(\phi_1,\ldots,\phi_n):=(\phi_1,\ldots,\phi_n,0,\ldots,0)$. In particular, the two
metrics induced from $\cK^{sp}$ and $\cL'^{sp}$ are compatible and therefore
$i_L^*h_L=h$.
\end{proof}
With all these preparations, we can state and prove the following theorem, which is the
first main result of this section.
\begin{theorem}
\label{theoCompleteness}
The distances $d_h$ on $|\LLp|$ and $|\MBLp|$ are strongly complete,
and so are the induced distances on the closed subspaces
$|\overline{U}_{\Spp} \cap \MBLp|$ for any fixed spectral pairs $\Spp$.
\end{theorem}
\begin{proof}
As $\MBLp$ and $\overline{U}_{\Spp} \cap \MBLp$ are closed analytic
subspaces of $\LLp$, it is sufficient to prove the completeness of the latter. We have to show that
there is an $\varepsilon>0$ such that for any $L\in\LLpr$, the
closed ball $B^{d_h}_\varepsilon(L)\subset \LLpr$ is compact. Using the
closed embedding $i_L:\LLp\hookrightarrow \Lambda^{pp}_b(W^\omega_L)$ and the estimate
$i_L^*d_{h_L}\leq d_h$ of lemma
\ref{lemEmbedMetricSpace}, it will be sufficient to show
that there is a uniform $\varepsilon$, such that for each $L\in\LLpr$,
$B^{d_{h_L}}_\varepsilon(L)\subset |\Lambda^{pp}_b(W^\omega_L)|$ is compact.
For any $L\in |\LLp|$ there is an $\varepsilon_L\in\dR_{>0}$ with this property
(this is exactly the property of being weakly complete,
always satisfied for locally compact spaces), but by lemma \ref{lemIsometry}
all spaces $|\Lambda_b^{pp}(W^\omega_L)|$ are isometric, so that $\varepsilon_L$ does not depend on $L$.
\end{proof}
Notice that the result applies in particular to the partial compactifications
$\overline{\check{D}}_{\mathit{BL}}\cap\MBLp$ of the
pure polarized classifying spaces $\cblp$ from \cite{HS2}.

In order to obtain a suitable target for period maps of variations of TERP-structures
over non simply-connected parameter spaces, we have to study quotients of
the classifying space $\MBLp$ by certain discrete groups. First we consider the real Lie group
$G_\dR:=\Aut(H^\infty_\dR,S,M)$.
\begin{lemma}
$G_\dR$ acts on $\MBL$, this action respects the strata $U_{\Spp}$ and their
closures $\overline{U}_{\Spp}$.
It acts by isometries on $\MBLp$ and on the intersection
$\overline{U}_{\Spp}\cap\MBLp$, in particular, on
$\overline{\check{D}}_{\mathit{BL}}\cap \MBL^{pp}$.
\end{lemma}
\begin{proof}
We first describe how to define the action of $G_\dR$ on $\MBL$.
Consider for any $\beta\in\dC$ the isomorphism $H^\infty\stackrel{es}{\rightarrow} W^{\beta}_{\beta+1}$
given by $es=\sum_{\alpha\in[\beta,\beta+1)}es_\alpha$,
where $es_\alpha:H^\infty_\lambda \rightarrow  W^\alpha_\alpha$ is defined by
$A \mapsto z^{\alpha \mathit{Id}-\frac{N}{2\pi i}} A$, here $e^{-2\pi i \alpha}=\lambda$.
Then $G_\dR\subset\Aut(H^\infty)$ acts on any $W^\beta_{\beta+1}$, and thus on
$W=\oplus_{k\in\dZ}z^k W^\beta_{\beta+1}$ and on
$W^\omega=W^{\alpha_1}_{\alpha_\mu-1}$.
This action commutes with the endomorphisms $b$ and $a$.
As $G_\dR$ respects the bilinear form $S$,
the action on $W^\omega$ respects the pairing $P$ and the symplectic form $\omega$, which
yields an action on $\Lambda(W^\omega)$. It induces
an action on $\MBL$ and an equivariant action on the
universal sheaf $\cL$. We see that the $V$-filtration
is stable under $G_\dR$, so that the spectral numbers do not change under this action
(neither do the spectral pairs, as $G_\dR$ respects $M$ and thus $N$),
which gives that $G_\dR(U_{\Spp})\subset U_{\Spp}$. Moreover, both
the involution $\tau$ and the pairing $P$ are respected by $G_\dR$,
so that we obtain finally an action on $\MBLp$ and a compatible
action on $\cL^{sp}$. The hermitian form $h$ is $G_\dR$-equivariant,
and so is the induced form on $\cAn_{\MBLp}\otimes \Theta_{\MBLp}$.
{}From this we conclude that $G_\dR\subset \textup{Isom}(\MBLpr,d_h)$
and $G_\dR\subset \textup{Isom}(|\overline{U}_{\Spp}\cap \MBL|,d_h)$.
\end{proof}

{}From now on and until the end of this section, we make the following additional
assumption which is virtually always satisfied for variations of TERP-structures
defined by families of geometric objects: There is a lattice $H^\infty_\dZ\subset H^\infty_\dR$ such that
$M\in\Aut(H^\infty_\dZ)$. Then we put $G_\dZ:=\Aut(H^\infty_\dZ,S,M)$.
In this situation, we have the following result.
\begin{theorem}
$G_\dZ$ acts properly discontinuously on $\MBLp$ and on $\overline{U}_{\Spp} \cap \MBLp$
so that the quotients $\MBLp/G_\dZ$ and $(\overline{U}_{\Spp} \cap \MBLp)/G_\dZ$
have the structure of complex spaces (this holds in particular
for the spaces $\overline{\check{D}}_{\mathit{BL}}\cap\MBLp$).
$\MBLp/G_\dZ$ resp. $(\overline{U}_{\Spp} \cap \MBLp)/G_\dZ$ are normal
if $\MBLp$ resp. $\overline{U}_{\Spp}\cap \MBLp$ are smooth.
\end{theorem}
Before entering into the proof of this theorem, we state and show the
following simple fact.
\begin{lemma}
\label{lemBaseChange}
Consider the free $\dC[z,z^{-1}]$-module $W$,
and let $\underline{v}^{(1)}$ and $\underline{v}^{(2)}$ be two bases of $W$ such that
$v^{(i)}_j\in W^{\alpha_1}_{\alpha_\mu}$ for $i\in\{1,2\}$ and all $j\in\{1,\ldots,\mu\}$.
Then the base change matrix between $\underline{v}^{(1)}$ and $\underline{v}^{(2)}$,
i.e. the matrix $M\in\Gl(\mu,\dC[z,z^{-1}])$ satisfying $\underline{v}^{(2)}=\underline{v}^{(1)} M$
can be written as $M=\sum\limits_{i=-n\mu}^{n\mu}M^{(k)} z^k$, where $M^{(k)}\in M(\mu\times\mu,\dC)$.
\end{lemma}
\begin{proof}
Choose a $\dC[z,z^{-1}]$-basis $\underline{v}^{(0)}\in(W^{\alpha_1}_{\alpha_1+1})^\mu$ of $W$, then we have matrices
$M_i \in\Gl(\mu,\dC[z,z^{-1}])$, ($i=1,2$) with coefficients in $\dC[z]_{\leq n}$ such that
$\underline{v}^{(i)}=\underline{v}^{(0)}\cdot M_i$.
This implies that $\det(M_i)=cz^k$ for some $c\in\dC^*$ and $k\in\{0,\ldots,\mu\cdot n\}$.
It follows that the coefficients of $M_1^{-1}$ are in
$\oplus_{k\in \dZ\cap[-\mu\cdot n,\mu\cdot (n-1)]}\dC z^k$ and the assertion
of the lemma is a consequence of $\underline{v}^{(2)}=\underline{v}^{(1)}\cdot M_1^{-1}\cdot M_2$.
\end{proof}
\begin{proof}[Proof of the theorem]
We fix once and for all a basis $\underline{A}$ of $H^\infty_\dZ$
which realizes $G_\dZ$ resp. $G_\dR$ as subgroups of $\Gl(\mu,\dZ)$
resp. $\Gl(\mu,\dR)$. We will show the following fact which implies that $G_\dZ$ acts
properly discontinuously: For any compact set $K\subset\MBLpr$ the
set
$$
\left\{a\in G_\dZ\,|\,a(K)\cap K\neq\emptyset\right\}
$$
is finite. As $G_\dZ$ is a discrete subgroup of $G_\dR$,
it is equivalent to show that the set
$$
\left\{a\in G_\dR\,|\,a(K)\cap K\neq\emptyset\right\}
$$
is compact. This can be reformulated by saying that there exists $R>0$ such
that for all $a\in G_\dR$ with $a(K)\cap K \neq \emptyset$ and all $i,j\in\{1,\ldots,\mu\}$
we have $|a^{mat}_{ij}|\leq R$, where $a^{mat}\in \Gl(\mu,\dR)$ is the
matrix of the automorphism $a$ with respect to the fixed basis $\underline{A}$.
We denote by $\underline{s}=(s_1,\ldots,s_\mu)\in W^{\alpha_1}_{\alpha_1+1}$ the basis which
corresponds to $\underline{A}$ under the isomorphism $es$. Then for any $a\in G_\dR$, the induced
action on $W^{\alpha_1}_{\alpha_1+1}$ is simply given by $\underline{s}\mapsto \underline{s}\cdot a^{mat}$.

Consider the hermitian bundle $\cL^{sp}\in\VBan_{\MBLp}$ and denote by
$U\!\cL^{sp}$ the total space of its associated bundle of $h$-orthonormal frames.
The projection $\pi:U\!\cL^{sp}\rightarrow \MBLp$ is proper (with fibres isomorphic to
$U(\mu)$). It follows that $\pi^{-1}(K)$ is compact.
Now we fix any element $\underline{v}^{(0)}\in\pi^{-1}(K)$, and write
$\underline{v}^{(0)}=\underline{s} \cdot \Gamma$ for some $\Gamma\in\Gl(\mu,\dC[z,z^{-1}])$.
Let $\underline{v}^{(1)}, \underline{v}^{(2)}\in\pi^{-1}(K)$ such that there is $a\in G_\dR$
with $a(\underline{v}^{(1)})=\underline{v}^{(2)}$. Then by lemma \ref{lemBaseChange}
there exist matrices $M^{(1)}=\sum_{k=-\mu n}^{\mu n} M^{(1)}_k z^k$ and
$M^{(2)}=\sum_{k=-\mu n}^{\mu n} M^{(2)}_k z^k$
such that $\underline{v}^{(0)}=\underline{v}^{(1)} M^{(1)}$ and
$\underline{v}^{(2)}=\underline{v}^{(0)} M^{(2)}$. This yields
$$
a^{mat}=\Gamma \cdot M^{(2)} \cdot M^{(1)} \cdot \Gamma^{-1} \in \Gl(\mu,\dR)
$$
The coefficients of the matrices $M^{(i)}_k$ are bounded by the compactness of $\pi^{-1}(K)$,
this implies that the coefficients of $a^{mat}$ are bounded by some positive real number $R$, as required.

It is obvious from the proof of the last lemma that the action of $G_\dZ$ on $\MBL$
respects $\overline{U}_{\Spp}$, so that it acts properly discontinuously
on $\overline{U}_{\Spp}\cap\MBLp$.
\end{proof}

\section{Examples and Applications}
\label{secExApp}

In this final section we first discuss in some detail
the geometry of several examples of the
classifying space $\MBL$ or of $\MBLred$. This illustrates the behavior
of the families of TERP-structures at boundary points,
in particular the jumping of the spectral pairs.
At this point it seems rather unclear which kind of varieties can appear as
these classifying spaces. We calculate in particular the limit TERP-structures on the boundary strata,
and discuss the $tt^*$-geometry as well as the relation to
the classifying spaces $\cbl$.
In most of the examples we care only about $\MBLred$.
Only in the first example in subsection \ref{subsecWPS} we care about $\MBL$ and indeed find
$\MBL\neq \MBLred$.
Finally, we use all the results proved so far to give some applications
for the study of period maps associated to variations of regular singular
TERP-structures.

\subsection{Smooth compactifications}
\label{subsecSmoothCompact}

In this first example the compact classifying space $\MBL$ (with its canonical complex structure)
is smooth, namely, it is the whole Lagrangian Grassmannian
(see lemma \ref{lemMBLsmooth}).
Consider the following initial data:
$H_\dR^\infty:=\oplus_{i=1}^4 \dR B_i$, $A_1:=B_1+iB_4$,
$A_2:=B_2+iB_3$, $\overline{A}_1=A_4$, $\overline{A}_2=A_3$,
$M(A_i):=e^{-2\pi i \alpha_i}A_i$, where we choose $\alpha_1$ and $\alpha_2$
with $-1<\alpha_1<\alpha_2<-\frac12$ and $\alpha_3=-\alpha_2$, $\alpha_4=-\alpha_1$. Moreover, we put
$w=0$ and
$$
S(\underline{B}^{tr},\underline{B}):=
\begin{pmatrix}
0 & 0 & 0 & -\gamma_1 \\
0 & 0 & -\gamma_2 & 0\\
0 & \gamma_2 & 0 & 0 \\
\gamma_1 & 0 & 0 & 0
\end{pmatrix}
$$
where $\gamma_1:=\frac{-1}{4\pi}\Gamma(\alpha_1+1)\Gamma(\alpha_4)$ and
$\gamma_2:=\frac{-1}{4\pi}\Gamma(\alpha_2+1)\Gamma(\alpha_3)$
(so that $S(A_1,A_4)=2i\gamma_1$ and $S(A_2,A_3)=2i\gamma_2$).
Let $s_i:=z^{\alpha_i}A_i$ then the relation between the pairing
$S$ and the pairing $P$ as expressed by \cite[formulas (5.4) and (5.5)]{HS1}
yields that
$$
P(\underline{s}^{tr},\underline{s}):=
\begin{pmatrix}
0 & 0 & 0 & 1 \\
0 & 0 & 1 & 0\\
0 & 1 & 0 & 0 \\
1 & 0 & 0 & 0
\end{pmatrix}.
$$
We define $\cH:=\oplus_{i=1}^4\cO_{\dC^4}v_i$, where
$$
\begin{array}{rclcrcl}
v_1 & := & s_1 + r z^{-1} s_3 + q z^{-1} s_4 & ; &
v_2 & := & s_2 + p z^{-1} s_3 + r z^{-1} s_4\\
v_3 & := & s_3 & ; &
v_4 & := & s_4
\end{array}
$$
Then $\TERP$ is a variation of TERP-structures of weight zero on $\dC^3$ with constant spectrum
$\Sp=(\alpha_1,\ldots,\alpha_4)$.
It is in fact the universal family of the classifying space $\cbl$ associated
to the initial data $(H^\infty,H^\infty_\dR, S, M, \Sp)$. The induced (constant) filtration
(recall definition \ref{defHodgeSpectrum}) is
$$
\{0\}=F^1\subsetneq F^0:=\dC A_1 \oplus \dC A_2 \subsetneq F^{-1} = H^\infty
$$
from which one checks that $(H^\infty, H_\dR^\infty, S, F^\bullet)$ is a pure polarized Hodge structure
of weight $-1$. Thus $\DcPMHS=\DPMHS=\{pt\}$ and $\cbl=D_{\mathit{BL}}$ in this case.
The situation is visualized in the following diagram, where each column represents
a space generated by elementary sections (that is, a space isomorphic to $W^\alpha_\alpha\cong V^\alpha/V^{>\alpha}$).
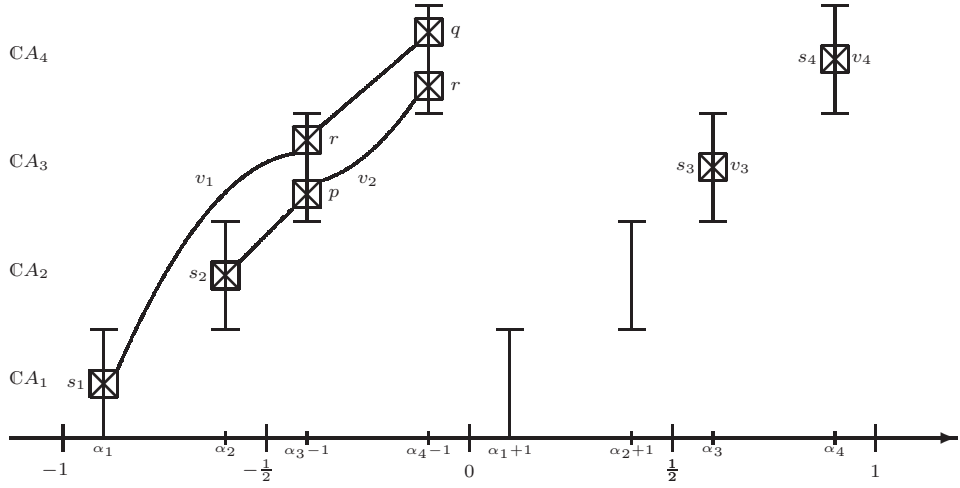
\begin{center}
\setlength{\unitlength}{0.36cm}
\begin{picture}(31,20)
\thicklines

\put(-2,2){\vector(1,0){35}}

\put(-2,4){\makebox{$\scriptstyle \dC A_1$}}
\put(-2,8){\makebox{$\scriptstyle \dC A_2$}}
\put(-2,12){\makebox{$\scriptstyle \dC A_3$}}
\put(-2,16){\makebox{$\scriptstyle \dC A_4$}}

\put(0,1.5){\line(0,1){1}}
\put(-0.3,0.8){\makebox(0,0){$\scriptstyle -1$}}
\put(7.5,1.5){\line(0,1){1}}
\put(7.2,0.8){\makebox(0,0){$\scriptstyle -\frac12$}}
\put(15,1.5){\line(0,1){1}}
\put(15,0.8){\makebox(0,0){$\scriptstyle 0$}}
\put(22.5,1.5){\line(0,1){1}}
\put(22.5,0.8){\makebox(0,0){$\scriptstyle \frac12$}}
\put(22.5,1.5){\line(0,1){1}}
\put(22.5,0.8){\makebox(0,0){$\scriptstyle \frac12$}}
\put(30,1.5){\line(0,1){1}}
\put(30,0.8){\makebox(0,0){$\scriptstyle 1$}}

\multiput(1.5,2)(15,0){2}{\line(0,1){4}}
\put(1.5,1.5){\makebox(0,0){$\scriptscriptstyle \alpha_1$}}
\put(16.5,1.5){\makebox(0,0){$\scriptscriptstyle \alpha_1+1$}}
\multiput(1,6)(15,0){2}{\line(1,0){1}}
\multiput(6,6)(15,0){2}{\line(0,1){4}}\multiput(6,1.75)(15,0){2}{\line(0,1){0.5}}
\put(6,1.5){\makebox(0,0){$\scriptscriptstyle \alpha_2$}}
\put(21,1.5){\makebox(0,0){$\scriptscriptstyle \alpha_2+1$}}
\multiput(5.5,6)(15,0){2}{\line(1,0){1}}
\multiput(5.5,10)(15,0){2}{\line(1,0){1}}

\multiput(9,10)(15,0){2}{\line(0,1){4}}
\multiput(9,1.75)(15,0){2}{\line(0,1){0.5}}
\put(9,1.5){\makebox(0,0){$\scriptscriptstyle \alpha_3-1$}}
\put(24,1.5){\makebox(0,0){$\scriptscriptstyle \alpha_3$}}
\multiput(8.5,10)(15,0){2}{\line(1,0){1}}
\multiput(8.5,14)(15,0){2}{\line(1,0){1}}

\multiput(13.5,14)(15,0){2}{\line(0,1){4}}
\multiput(13.5,1.75)(15,0){2}{\line(0,1){0.5}}
\put(13.5,1.5){\makebox(0,0){$\scriptscriptstyle \alpha_4-1$}}
\put(28.5,1.5){\makebox(0,0){$\scriptscriptstyle \alpha_4$}}
\multiput(13,14)(15,0){2}{\line(1,0){1}}
\multiput(13,18)(15,0){2}{\line(1,0){1}}

\multiput(1,3.5)(4.5,4){2}{\line(1,0){1}}\multiput(1,3.5)(4.5,4){2}{\line(0,1){1}}
\multiput(2,3.5)(4.5,4){2}{\line(0,1){1}}\multiput(1,4.5)(4.5,4){2}{\line(1,0){1}}
\multiput(1,3.5)(4.5,4){2}{\line(1,1){1}}\multiput(1,4.5)(4.5,4){2}{\line(1,-1){1}}
\qbezier(2,4.5)(5.25,12)(8.5,12.5)
\qbezier(9.5,13.5)(9.5,13.5)(13,16.5)
\qbezier(6.5,8.5)(6.5,8.5)(8.5,10.5)
\qbezier(9.5,11.5)(11.25,12)(13,14.5)
\put(0.5,4){\makebox(0,0){$\scriptstyle s_1$}}
\put(5,8){\makebox(0,0){$\scriptstyle s_2$}}
\put(5.25,11.5){\makebox(0,0){$\scriptstyle v_1$}}
\put(11.25,11.5){\makebox(0,0){$\scriptstyle v_2$}}

\multiput(8.5,10.5)(4.5,4){2}{\line(1,0){1}}\multiput(8.5,10.5)(4.5,4){2}{\line(0,1){1}}
\multiput(9.5,10.5)(4.5,4){2}{\line(0,1){1}}\multiput(8.5,11.5)(4.5,4){2}{\line(1,0){1}}
\multiput(8.5,10.5)(4.5,4){2}{\line(1,1){1}}\multiput(8.5,11.5)(4.5,4){2}{\line(1,-1){1}}
\put(10,11){\makebox(0,0){$\scriptstyle p$}}
\put(10,13){\makebox(0,0){$\scriptstyle r$}}
\put(14.5,15){\makebox(0,0){$\scriptstyle r$}}
\put(14.5,17){\makebox(0,0){$\scriptstyle q$}}

\multiput(8.5,12.5)(4.5,4){2}{\line(1,0){1}}\multiput(8.5,12.5)(4.5,4){2}{\line(0,1){1}}
\multiput(9.5,12.5)(4.5,4){2}{\line(0,1){1}}\multiput(8.5,13.5)(4.5,4){2}{\line(1,0){1}}
\multiput(8.5,12.5)(4.5,4){2}{\line(1,1){1}}\multiput(8.5,13.5)(4.5,4){2}{\line(1,-1){1}}

\multiput(23.5,11.5)(4.5,4){2}{\line(1,0){1}}\multiput(23.5,11.5)(4.5,4){2}{\line(0,1){1}}
\multiput(24.5,11.5)(4.5,4){2}{\line(0,1){1}}\multiput(23.5,12.5)(4.5,4){2}{\line(1,0){1}}
\multiput(23.5,11.5)(4.5,4){2}{\line(1,1){1}}\multiput(23.5,12.5)(4.5,4){2}{\line(1,-1){1}}
\put(23,12){\makebox(0,0){$\scriptstyle s_3$}}
\put(27.5,16){\makebox(0,0){$\scriptstyle s_4$}}
\put(25,12){\makebox(0,0){$\scriptstyle v_3$}}
\put(29.5,16){\makebox(0,0){$\scriptstyle v_4$}}
\end{picture}
\end{center}
In this example it is quite easy to describe the space $\MBL$:
As $\lfloor\alpha_\mu-\alpha_1\rfloor=1$ we have by lemma \ref{lemMBLsmooth} that
$\MBL=\Lambda(W^{\alpha_1}_{\alpha_4-1}, [P^{(-1)}])$,
where $[P^{(-1)}]=[z^{-1}s_4]^*\wedge [s_1]^* + [z^{-1}s_3]^*\wedge [s_2]^* \in \bigwedge^2 (W^{\alpha_1}_{\alpha_4-1})^*$.
Using
the Pl{\"u}cker embedding, one checks that this Lagrangian Grassmannian
is a hyperplane section of the Pl\"ucker quadric in $\dP^5$, i.e., a smooth quadric in $\dP^4$.
Such a smooth quadric is isomorphic neither to $\dP^3$ nor to $\dP^1\times\dP^2$.
It is also clear that $\MBL$ is indeed the closure of the three-dimensional
affine classifying space $\cbl$ considered above.
The stratification of the boundary of $\MBL$ is as follows:
$\MBL \backslash \cbl=\overline{U}_{\Sp_2}$ where
$\Sp_2=(\alpha_1,\alpha_3-1, \alpha_2+1, \alpha_4)$,
$\overline{U}_{\Sp_2}\backslash U_{Sp_2}=
\overline{U}_{\Sp_1}$, where
$\Sp_1=(\alpha_2,\alpha_4-1, \alpha_1+1, \alpha_3)$ and
$\overline{U}_{\Sp_1}\backslash U_{Sp_1}=\overline{U}_{Sp_0}=U_{Sp_0}=\{pt\}$,
where $\Sp_0=(\alpha_3-1,\alpha_4-1,\alpha_1+1,\alpha_2+1)$.
In all cases we have $\dim(U_{\Sp_i})=i$.

Let us describe these strata with some more details:
We have $U_{\Sp_2} \cong \dC^2=\Spec\dC[x,y]$, namely, the universal family is
$\cH^{(2)}=\oplus_{i=1}^4\cO_{\dC^3}v^{(2)}_i$ where
$$
\begin{array}{rclcrcl}
v^{(2)}_1 & := & s_1-xs_2+yz^{-1}s_4 & ; &
v^{(2)}_2 & := & z^{-1}s_3+xz^{-1}s_4\\
v^{(2)}_3 & := & z s_2 &; &
v^{(2)}_4 & := & s_4.
\end{array}
$$
This is shown in the following diagram.
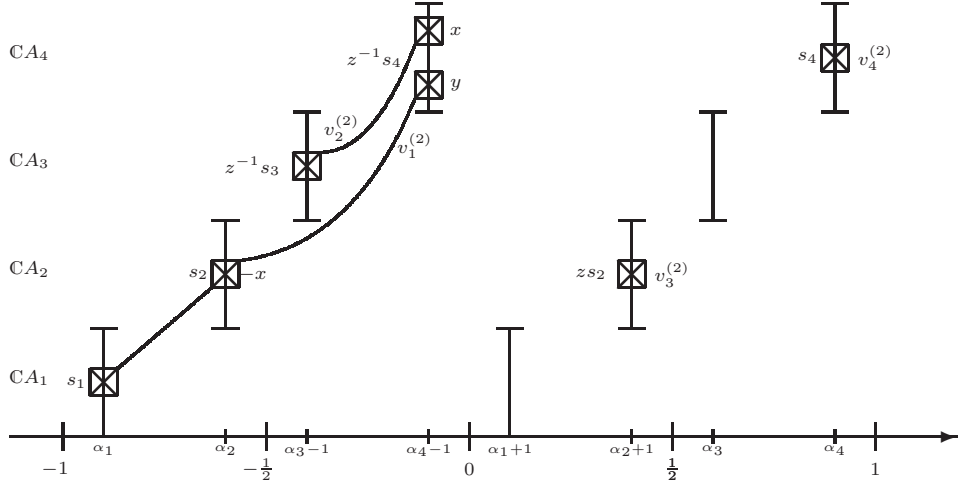
\begin{center}
\setlength{\unitlength}{0.36cm}
\begin{picture}(31,20)
\thicklines

\put(-2,2){\vector(1,0){35}}

\put(-2,4){\makebox{$\scriptstyle \dC A_1$}}
\put(-2,8){\makebox{$\scriptstyle \dC A_2$}}
\put(-2,12){\makebox{$\scriptstyle \dC A_3$}}
\put(-2,16){\makebox{$\scriptstyle \dC A_4$}}

\put(0,1.5){\line(0,1){1}}
\put(-0.3,0.8){\makebox(0,0){$\scriptstyle -1$}}
\put(7.5,1.5){\line(0,1){1}}
\put(7.2,0.8){\makebox(0,0){$\scriptstyle -\frac12$}}
\put(15,1.5){\line(0,1){1}}
\put(15,0.8){\makebox(0,0){$\scriptstyle 0$}}
\put(22.5,1.5){\line(0,1){1}}
\put(22.5,0.8){\makebox(0,0){$\scriptstyle \frac12$}}
\put(22.5,1.5){\line(0,1){1}}
\put(22.5,0.8){\makebox(0,0){$\scriptstyle \frac12$}}
\put(30,1.5){\line(0,1){1}}
\put(30,0.8){\makebox(0,0){$\scriptstyle 1$}}

\multiput(1.5,2)(15,0){2}{\line(0,1){4}}
\put(1.5,1.5){\makebox(0,0){$\scriptscriptstyle \alpha_1$}}
\put(16.5,1.5){\makebox(0,0){$\scriptscriptstyle \alpha_1+1$}}
\multiput(1,6)(15,0){2}{\line(1,0){1}}
\multiput(6,6)(15,0){2}{\line(0,1){4}}\multiput(6,1.75)(15,0){2}{\line(0,1){0.5}}
\put(6,1.5){\makebox(0,0){$\scriptscriptstyle \alpha_2$}}
\put(21,1.5){\makebox(0,0){$\scriptscriptstyle \alpha_2+1$}}
\multiput(5.5,6)(15,0){2}{\line(1,0){1}}
\multiput(5.5,10)(15,0){2}{\line(1,0){1}}

\multiput(9,10)(15,0){2}{\line(0,1){4}}
\multiput(9,1.75)(15,0){2}{\line(0,1){0.5}}
\put(9,1.5){\makebox(0,0){$\scriptscriptstyle \alpha_3-1$}}
\put(24,1.5){\makebox(0,0){$\scriptscriptstyle \alpha_3$}}
\multiput(8.5,10)(15,0){2}{\line(1,0){1}}
\multiput(8.5,14)(15,0){2}{\line(1,0){1}}

\multiput(13.5,14)(15,0){2}{\line(0,1){4}}
\multiput(13.5,1.75)(15,0){2}{\line(0,1){0.5}}
\put(13.5,1.5){\makebox(0,0){$\scriptscriptstyle \alpha_4-1$}}
\put(28.5,1.5){\makebox(0,0){$\scriptscriptstyle \alpha_4$}}
\multiput(13,14)(15,0){2}{\line(1,0){1}}
\multiput(13,18)(15,0){2}{\line(1,0){1}}

\multiput(1,3.5)(4.5,4){2}{\line(1,0){1}}\multiput(1,3.5)(4.5,4){2}{\line(0,1){1}}
\multiput(2,3.5)(4.5,4){2}{\line(0,1){1}}\multiput(1,4.5)(4.5,4){2}{\line(1,0){1}}
\multiput(1,3.5)(4.5,4){2}{\line(1,1){1}}\multiput(1,4.5)(4.5,4){2}{\line(1,-1){1}}
\qbezier(2,4.5)(2,4.5)(5.5,7.5)
\qbezier(6.5,8.5)(10.75,9)(13,14.5)
\qbezier(9.5,12.5)(11.5,12.5)(13,16.5)
\put(0.5,4){\makebox(0,0){$\scriptstyle s_1$}}
\put(5,8){\makebox(0,0){$\scriptstyle s_2$}}
\put(7,8){\makebox(0,0){$\scriptstyle -x$}}
\put(7,12){\makebox(0,0){$\scriptstyle z^{-1}s_3$}}
\put(10.3,13.4){\makebox(0,0){$\scriptstyle v^{(2)}_2$}}
\put(13,12.75){\makebox(0,0){$\scriptstyle v^{(2)}_1$}}
\put(11.5,16){\makebox(0,0){$\scriptstyle z^{-1}s_4$}}

\put(13,14.5){\line(1,0){1}}\put(13,14.5){\line(0,1){1}}
\put(14,14.5){\line(0,1){1}}\put(13,15.5){\line(1,0){1}}
\put(13,14.5){\line(1,1){1}}\put(13,15.5){\line(1,-1){1}}
\put(14.5,15){\makebox(0,0){$\scriptstyle y$}}
\put(14.5,17){\makebox(0,0){$\scriptstyle x$}}

\put(8.5,11.5){\line(1,0){1}}\put(8.5,11.5){\line(0,1){1}}
\put(9.5,11.5){\line(0,1){1}}\put(8.5,12.5){\line(1,0){1}}
\put(8.5,11.5){\line(1,1){1}}\put(8.5,12.5){\line(1,-1){1}}

\put(13,16.5){\line(1,0){1}}\put(13,16.5){\line(0,1){1}}
\put(14,16.5){\line(0,1){1}}\put(13,17.5){\line(1,0){1}}
\put(13,16.5){\line(1,1){1}}\put(13,17.5){\line(1,-1){1}}

\multiput(20.5,7.5)(7.5,8){2}{\line(1,0){1}}\multiput(20.5,7.5)(7.5,8){2}{\line(0,1){1}}
\multiput(21.5,7.5)(7.5,8){2}{\line(0,1){1}}\multiput(20.5,8.5)(7.5,8){2}{\line(1,0){1}}
\multiput(20.5,7.5)(7.5,8){2}{\line(1,1){1}}\multiput(20.5,8.5)(7.5,8){2}{\line(1,-1){1}}
\put(19.5,8){\makebox(0,0){$\scriptstyle zs_2$}}
\put(27.5,16){\makebox(0,0){$\scriptstyle s_4$}}
\put(22.5,8){\makebox(0,0){$\scriptstyle v^{(2)}_3$}}
\put(30,16){\makebox(0,0){$\scriptstyle v^{(2)}_4$}}
\end{picture}
\end{center}
The induced filtration $F^\bullet$ has changed, it is now given by
$$
\{0\}=F^1\subsetneq F^0:=\dC A_1 \oplus \dC A_3 \subsetneq F^{-1} = H^\infty.
$$
One checks that this is still a pure Hodge structure of weight $-1$ on $H^\infty$,
but the Hodge metric has signature $(+,-,-,+)$.

The compactification $\overline{U}_{\Sp_2}$ can be calculated in a rather direct way.
Namely, we take the basis $\underline{\widetilde{v}}^{(2)}$ of the restriction $\cH^{(2)}_{|x\neq 0}$,
given by
$$
\begin{array}{rclcl}
\widetilde{v}^{(2)}_1 & := &x^{-1}v^{(2)}_1-y\cdot x^{-2}\cdot v^{(2)}_2 & =  &
x^{-1}s_1 - s_2 -yx^{-2}z^{-1}s_3 \\
\widetilde{v}^{(2)}_2 & := & x^{-1}v^{(2)}_2 & = &
x^{-1}z^{-1}s_3+z^{-1}s_4\\
\widetilde{v}^{(2)}_3 & := & zv^{(2)}_1+xv^{(2)}_3-yv^{(2)}_4 & = &zs_1\\
\widetilde{v}^{(2)}_4 & := & zv^{(2)}_2-xv^{(2)}_4  & = & s_3
\end{array}
$$
For a fixed parameter $w\in\dC$, consider the restriction
$\cH^{(2)}_{|y-w x^2=0,x\neq 0}$ which has a basis
$$
\left(\widetilde{v}^{(2)}_1,\widetilde{v}^{(2)}_2,\widetilde{v}^{(2)}_3,\widetilde{v}^{(2)}_4\right)_{|y-w x^2=0,x\neq 0}
=\left(x^{-1} s_1-s_2-wz^{-1}s_3, x^{-1} z^{-1} s_3+z^{-1} s_4,
zs_1,s_3\right).
$$
This family extends to $x=\infty$, namely
$$
\lim_{x\rightarrow \infty}\cH^{(2)}_{|y-w x^2=0,x\neq 0}
=\cO_\dC(s_2+w z^{-1} s_3) \oplus
\cO_\dC z^{-1} s_4 \oplus
\cO_\dC z s_1 \oplus
\cO_\dC s_3.
$$
These extensions are pairwise non-isomorphic for different parameters $w$, so
that the closure $U_{\Sp_2}$ is isomorphic to
$\dP(1,1,2)=\Proj\dC[X,Y,Z]$, where $\deg(X)=1, \deg(Y)=1$ and $\deg(Z)=2$
and where the embedding into the compactification is given by
$$
\begin{array}{rcl}
U_{\Sp_2} & \hookrightarrow &
\overline{U}_{\Sp_2} \subset \MBL\\ \\
(x,y)& \longmapsto & (1,x,y)=(X,Y,Z).
\end{array}
$$
An explicit calculation using the Pl\"ucker embedding
$$
\textup{Gr}(2,W^\omega)\cong V(AF-BE+CD) \hookrightarrow \Proj\dC[A,B,C,D,E,F]=\Proj\textup{Sym}^\bullet\,(\Lambda^2(W^\omega))
$$
shows that $\overline{U}_{\Sp_2}$ is the variety
$V(AF-BE+CD, D+C, A)\subset\dP^5$, which gives
again that $\overline{U}_{\Sp_2} \cong \dP(1,1,2)$.

It follows that the
boundary $\overline{U}_{\Sp_2}\backslash U_{\Sp_2}$ is isomorphic
to $\dP(1,2)\cong\dP^1$. In particular, the interior $U_{\Sp_1}$
of this boundary is isomorphic to $\dC$, with the universal family given by
$\cH^{(1)}
=\oplus_{i=1}^4 \cO_{\dC^2}v^{(1)}_i$
$$
\begin{array}{c}
v^{(1)}_1 := s_2+wz^{-1}s_3\quad ; \quad v^{(1)}_2 := z^{-1}s_4 \quad ; \quad
v^{(1)}_3 := z s_1 \quad ; \quad v^{(1)}_4 := s_3
\end{array}
$$
This situation looks as follows.
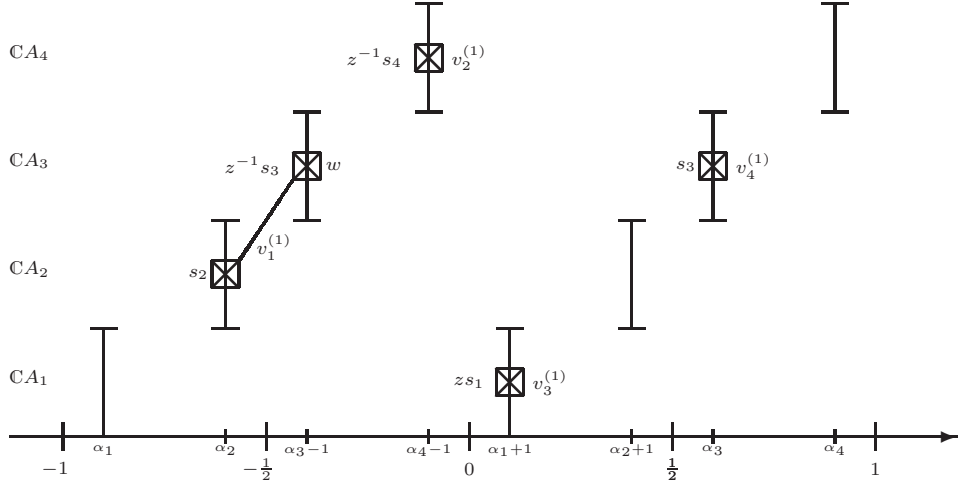
\begin{center}
\setlength{\unitlength}{0.36cm}
\begin{picture}(31,20)
\thicklines

\put(-2,2){\vector(1,0){35}}

\put(-2,4){\makebox{$\scriptstyle \dC A_1$}}
\put(-2,8){\makebox{$\scriptstyle \dC A_2$}}
\put(-2,12){\makebox{$\scriptstyle \dC A_3$}}
\put(-2,16){\makebox{$\scriptstyle \dC A_4$}}

\put(0,1.5){\line(0,1){1}}
\put(-0.3,0.8){\makebox(0,0){$\scriptstyle -1$}}
\put(7.5,1.5){\line(0,1){1}}
\put(7.2,0.8){\makebox(0,0){$\scriptstyle -\frac12$}}
\put(15,1.5){\line(0,1){1}}
\put(15,0.8){\makebox(0,0){$\scriptstyle 0$}}
\put(22.5,1.5){\line(0,1){1}}
\put(22.5,0.8){\makebox(0,0){$\scriptstyle \frac12$}}
\put(22.5,1.5){\line(0,1){1}}
\put(22.5,0.8){\makebox(0,0){$\scriptstyle \frac12$}}
\put(30,1.5){\line(0,1){1}}
\put(30,0.8){\makebox(0,0){$\scriptstyle 1$}}

\multiput(1.5,2)(15,0){2}{\line(0,1){4}}
\put(1.5,1.5){\makebox(0,0){$\scriptscriptstyle \alpha_1$}}
\put(16.5,1.5){\makebox(0,0){$\scriptscriptstyle \alpha_1+1$}}
\multiput(1,6)(15,0){2}{\line(1,0){1}}
\multiput(6,6)(15,0){2}{\line(0,1){4}}\multiput(6,1.75)(15,0){2}{\line(0,1){0.5}}
\put(6,1.5){\makebox(0,0){$\scriptscriptstyle \alpha_2$}}
\put(21,1.5){\makebox(0,0){$\scriptscriptstyle \alpha_2+1$}}
\multiput(5.5,6)(15,0){2}{\line(1,0){1}}
\multiput(5.5,10)(15,0){2}{\line(1,0){1}}

\multiput(9,10)(15,0){2}{\line(0,1){4}}
\multiput(9,1.75)(15,0){2}{\line(0,1){0.5}}
\put(9,1.5){\makebox(0,0){$\scriptscriptstyle \alpha_3-1$}}
\put(24,1.5){\makebox(0,0){$\scriptscriptstyle \alpha_3$}}
\multiput(8.5,10)(15,0){2}{\line(1,0){1}}
\multiput(8.5,14)(15,0){2}{\line(1,0){1}}

\multiput(13.5,14)(15,0){2}{\line(0,1){4}}
\multiput(13.5,1.75)(15,0){2}{\line(0,1){0.5}}
\put(13.5,1.5){\makebox(0,0){$\scriptscriptstyle \alpha_4-1$}}
\put(28.5,1.5){\makebox(0,0){$\scriptscriptstyle \alpha_4$}}
\multiput(13,14)(15,0){2}{\line(1,0){1}}
\multiput(13,18)(15,0){2}{\line(1,0){1}}

\put(5.5,7.5){\line(1,0){1}}\put(5.5,7.5){\line(0,1){1}}
\put(6.5,7.5){\line(0,1){1}}\put(5.5,8.5){\line(1,0){1}}
\put(5.5,7.5){\line(1,1){1}}\put(5.5,8.5){\line(1,-1){1}}
\qbezier(6.5,8.5)(6.5,8.5)(8.5,11.5)
\put(5,8){\makebox(0,0){$\scriptstyle s_2$}}
\put(7.8,9){\makebox(0,0){$\scriptstyle v^{(1)}_1$}}
\put(7,12){\makebox(0,0){$\scriptstyle z^{-1}s_3$}}
\put(10,12){\makebox(0,0){$\scriptstyle w$}}
\put(15,16){\makebox(0,0){$\scriptstyle v^{(1)}_2$}}
\put(11.5,16){\makebox(0,0){$\scriptstyle z^{-1}s_4$}}
\put(15,4){\makebox(0,0){$\scriptstyle zs_1$}}
\put(18,4){\makebox(0,0){$\scriptstyle v^{(1)}_3$}}

\multiput(8.5,11.5)(15,0){2}{\line(1,0){1}}\multiput(8.5,11.5)(15,0){2}{\line(0,1){1}}
\multiput(9.5,11.5)(15,0){2}{\line(0,1){1}}\multiput(8.5,12.5)(15,0){2}{\line(1,0){1}}
\multiput(8.5,11.5)(15,0){2}{\line(1,1){1}}\multiput(8.5,12.5)(15,0){2}{\line(1,-1){1}}

\put(13,15.5){\line(1,0){1}}\put(13,15.5){\line(0,1){1}}
\put(14,15.5){\line(0,1){1}}\put(13,16.5){\line(1,0){1}}
\put(13,15.5){\line(1,1){1}}\put(13,16.5){\line(1,-1){1}}

\put(16,3.5){\line(1,0){1}}\put(16,3.5){\line(0,1){1}}
\put(17,3.5){\line(0,1){1}}\put(16,4.5){\line(1,0){1}}
\put(16,3.5){\line(1,1){1}}\put(16,4.5){\line(1,-1){1}}
\put(23,12){\makebox(0,0){$\scriptstyle s_3$}}
\put(25.5,12){\makebox(0,0){$\scriptstyle v^{(1)}_4$}}
\end{picture}
\end{center}
Now the filtration is
$$
\{0\}=F^1\subsetneq F^0:=\dC A_2 \oplus \dC A_4 \subsetneq F^{-1} = H^\infty
$$
which is again pure of weight $-1$ but not polarized, the hermitian form has
signature $(-,+,+,-)$.

Finally, the stratum $U_{\Sp_0}=\{\lim_{w\rightarrow \infty}\cH^{(1)}\}$ is the single TERP-structure $\cH^{(0)}
=\oplus_{i=1}^4 \cO_{\dC}v^{(0)}_i$ where
$$
\begin{array}{c}
v^{(0)}_1 := z^{-1}s_3\quad ; \quad v^{(0)}_2 := z^{-1}s_4 \quad ; \quad
v^{(0)}_3 := z s_1 \quad ; \quad v^{(0)}_4 := z s_2
\end{array}
$$
and
the associated filtration is given by
$$
\{0\}=F^1\subsetneq F^0:=\dC A_3 \oplus \dC A_4 \subsetneq F^{-1} = H^\infty
$$
which is pure of weight $-1$, and the hermitian form is negative definite.

The following is a brief description of the $tt^*$-geometry on the different strata
$U_{\Sp_i}$. The simplest one is the zero-dimensional stratum $U_{\Sp_0}$: Its
single TERP-structure is generated by elementary sections, and we see that it is
pure but not polarized, the hermitian form $h$ on $\widehat{\cH}^{(0)}$ is equal
to the Hodge metric on $H^\infty$, which is negative definite.

The universal family $\cH^{(1)}$ on $U_{\Sp_1}$ is actually a sum of two TERP-structures,
namely the one generated by $v^{(1)}_1$ and $v^{(1)}_4$ and the one generated by
$v^{(1)}_2$ and $v^{(1)}_3$. The latter is generated by elementary sections and is
pure but with negative definite hermitian metric on the space of global sections,
the former is pure outside the hypersurface $|w|=1$ and polarized for $|w|<1$.
The two-dimensional family $\cH^{(2)}$ is pure outside the real-analytic hypersurface
$D=\{|y|=|x|^2+1\}$, we have
$$
(p_*\widehat{\cH}^{(2)})_{|U_{\Sp_2}\backslash D}=
\cO_{\dP^1}\cAn_{U_{\Sp_2}\backslash D}\, v^{(2)}_1\oplus\cO_{\dP^1}\cAn_{U_{\Sp_2}\backslash D}\, v^{(2)}_2\oplus
\cO_{\dP^1}\cAn_{U_{\Sp_2}\backslash D}\, \tau v^{(2)}_1\oplus\cO_{\dP^1}\cAn_{U_{\Sp_2}\backslash D}\, \tau v^{(2)}_2,
$$
and the hermitian form $h$ has signatures $(+-+-)$ resp. $(----)$ on
$\{|y|<|x|^2+1\}$ resp $\{|y|>|x|^2+1\}$. The following picture visualizes
the situation.

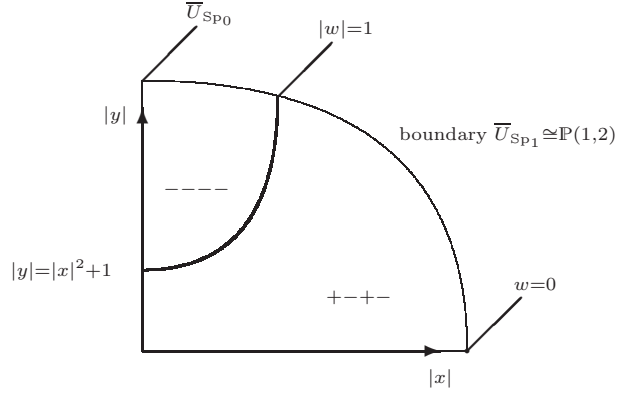
\begin{center}
\setlength{\unitlength}{0.36cm}
\begin{picture}(36,16)
\thicklines

\put(12,2){\vector(1,0){11}}
\put(12,2){\vector(0,1){9}}

\qbezier(12,5)(17,5)(17,11.4)

\put(23,1){\makebox(0,0){$\scriptstyle |x|$}}
\put(11,10.7){\makebox(0,0){$\scriptstyle |y|$}}
\put(9,5){\makebox(0,0){$\scriptstyle |y|=|x|^2+1$}}

\put(25.5,10){\makebox(0,0){$\scriptstyle \textup{boundary }\overline{U}_{\Sp_1}\cong\dP(1,2)$}}

\put(24,2){\circle{0}}

\linethickness{0.1mm}
\qbezier(12,12)(24,12)(24,2)

\put(12,12){\line(1,1){2}}
\put(14.5,14.5){\makebox(0,0){$\scriptstyle \overline{U}_{\Sp_0}$}}

\put(17,11.4){\line(1,1){2}}
\put(19.5,13.9){\makebox(0,0){$\scriptstyle |w|=1$}}

\put(24,2){\line(1,1){2}}
\put(26.5,4.5){\makebox(0,0){$\scriptstyle w=0$}}

\put(12,2){\line(1,0){12}}
\put(12,2){\line(0,1){10}}

\put(20,4){\makebox(0,0){$\scriptstyle +-+-$}}
\put(14,8){\makebox(0,0){$\scriptstyle ----$}}

\end{picture}
\end{center}

{}From these observations we see that the pure polarized part $\MBLp$ is contained
inside the generic stratum $\cbl$ of $\MBL$ and then it follows that it
must be relatively compact in $\cbl$, for otherwise theorem \ref{theoLimitTwistor} would give that there is
a pure polarized point in the boundary $\MBL\backslash \cbl =\overline{U}_{\Sp_2}$ (note that
for any curve approaching a boundary point through $\cblp$, the limit is
necessarily a pure polarized TERP-structure, as the corresponding monodromy is
reduced to the identity so that the weight filtration $W_\bullet$ appearing
in theorem \ref{theoLimitTwistor} is trivial).
Moreover, the fact that the filtration $F^\bullet=\left(\{0\}=F^1\subsetneq
F^0=\dC A_1\oplus \dC A_2 \subsetneq F^{-1}=H^\infty\right)$ gives rise
to a PHS shows that $\MBLp$ is non-empty, as it contains the
origin $r=p=q=0$ in $\cbl\cong\dC^3 \subset \MBL$.

\vspace*{1cm}

In the second example of this subsection, the compactification $\MBLred$ is also smooth
although $n=\lfloor \alpha_\mu-\alpha_1 \rfloor >1$
(here we consider only $\MBLred$, not $\MBL$). Let $H^\infty_\dR:=\dR B_1 \oplus
\dR B_2$, $S(\underline{B}^{tr},\underline{B})=\begin{pmatrix}0&-\pi^2\\\pi^2& 0\end{pmatrix}$.
A reference Hodge structure of weight one is given by
$$
\{0\}=F^3 \subsetneq F^2_0=\dC A = F^1_0 = F^0_0 \subsetneq F^{-1}_0:=H^\infty
$$
where $A:=B_1+iB_2$. Indeed, we have $F^2_0 \oplus \overline{F}_0^0
= F^1_0 \oplus \overline{F}^1_0 = F^0_0 \oplus \overline{F}^2_0 = H^\infty$ and
$i^{2-(1-2)}S(A,\overline{A})=-iS(A,\overline{A})=2S(B_2,B_1)>0$ (Note that
the isotropy condition $S(F^p, F^{2-p})=0$ is automatically satisfied as $S$ is symplectic).
The classifying spaces $\DPMHS\subsetneq \DcPMHS$ are well-known
(see, e.g., \cite[\S 5]{Sch}): $\DcPMHS\cong \dP^1$ and $\DPMHS=\dH$.
A point $(x:y)\in\dP^1$ corresponds to the filtration given by $F^2:=\dC(xB_1+yB_2)$,
so that $F^\bullet_0 = (1:i)\in \dH\subset \dP^1$. The complement
$\DcPMHS\backslash \DPMHS$ is the union $\overline{\dH}\cup \dP^1_\dR\subset\dP^1$,
where the points of the real projective line are non-pure Hodge filtrations, whereas
the points in $\overline{\dH}$ are pure but the Hodge metric
$h:=-iS(\cdot,\overline{\cdot})_{|H^{2,-1}\times H^{2,-1}}\oplus\,iS(\cdot,\overline{\cdot})_{|H^{-1,2}\times H^{-1,2}}$
is negative definite. The pairing $P$ is given by $P(s_i,s_j)=(-1)^{j+1} z \delta_{i+j,3}$,
where $s_1:=z^{1/2}A$ and $s_2:=z^{1/2}\overline{A}$.

We put $w=2$ and $\alpha=-\frac12$, then the classifying space $\cbl$ associated to the spectrum
$(\alpha_1,\alpha_2)=(-\frac12, \frac52)$ is the total space $\mathbf{V}(\cE)$ of a line bundle $\cE$
over $\DcPMHS$, the universal family over a fibre $\check{\pi}_{BL}^{-1}(F^\bullet)\cong\Spec\dC[r]$
of the projection $\check{\pi}_{BL}:\cbl\rightarrow \DcPMHS$ is given by
$$
\cH:=\cO_{\dC^2}\left[\underbrace{z^{-1/2}A_1+rz^{3/2}A_2}_{v_1}\right]\oplus\cO_{\dC^2}\underbrace{z^{5/2}A_2}_{v_2}.
$$
where $F^\bullet=\left(\{0\}\subsetneq F^2=\dC A_1=F^0\subsetneq F^{-1}=\dC A_1\oplus\dC A_2=H^\infty\right)$.
For any such family, letting $r$ tend
to infinity yields the limit structure $\cG:=\cH_{r\rightarrow \infty}=\cO_\dC g_1 \oplus \cO_\dC g_2$, where
$g_1=z^{1/2}A_1$ and $g_2=z^{3/2}A_2$. The stratum at infinity is $U_{(\frac12,\frac32)} \cong
\dP^1$, which shows that $\cbl$ is compactified to $\MBLred$ along the fibres
of $\check{\pi}_{BL}$, so that it must be a Hirzebruch surface $\Sigma_k:=
\mathbf{Proj}\,(\cO_{\dP^1}(k)\oplus\cO_{\dP^1})$, where $k=\deg(\cE)$. The following picture
shows the situation.
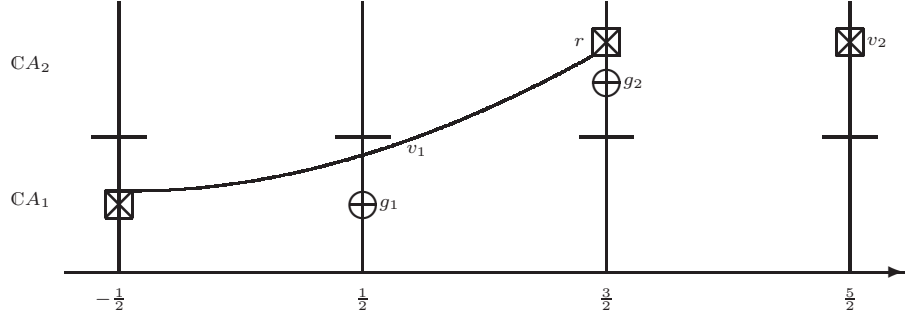
\begin{center}
\setlength{\unitlength}{0.36cm}
\begin{picture}(26,14)
\thicklines

\put(0,2){\vector(1,0){31}}

\put(-2,4.5){\makebox{$\scriptstyle \dC A_1$}}
\put(-2,9.5){\makebox{$\scriptstyle \dC A_2$}}

\multiput(2,2)(9,0){4}{\line(0,1){10}}
\multiput(1,7)(9,0){4}{\line(1,0){2}}

\multiput(1.5,4)(18,6){2}{\line(1,0){1}}\multiput(1.5,4)(18,6){2}{\line(0,1){1}}
\multiput(2.5,4)(18,6){2}{\line(0,1){1}}\multiput(1.5,5)(18,6){2}{\line(1,0){1}}
\multiput(1.5,4)(18,6){2}{\line(1,1){1}}\multiput(1.5,5)(18,6){2}{\line(1,-1){1}}

\put(28.5,10){\line(1,0){1}}\put(28.5,10){\line(0,1){1}}
\put(29.5,10){\line(0,1){1}}\put(28.5,11){\line(1,0){1}}
\put(28.5,10){\line(1,1){1}}\put(28.5,11){\line(1,-1){1}}

\qbezier(2.5,5)(10.5,5)(19.5,10)

\multiput(11,4.5)(9,4.5){2}{\circle{1}}
\multiput(10.5,4.5)(9,4.5){2}{\line(1,0){1}}

\put(13,6.5){\makebox(0,0){$\scriptstyle v_1$}}
\put(30,10.5){\makebox(0,0){$\scriptstyle v_2$}}
\put(19,10.5){\makebox(0,0){$\scriptstyle r$}}
\put(12,4.5){\makebox(0,0){$\scriptstyle g_1$}}
\put(21,9){\makebox(0,0){$\scriptstyle g_2$}}

\put(1.7,1){\makebox(0,0){$\scriptstyle -\frac12$}}
\put(11,1){\makebox(0,0){$\scriptstyle \frac12$}}
\put(20,1){\makebox(0,0){$\scriptstyle \frac32$}}
\put(29,1){\makebox(0,0){$\scriptstyle \frac52$}}
\end{picture}
\end{center}

The degree of $\cE$ is calculated as follows:
Let $(x:y)$ be homogeneous coordinates on $\cbl\cong\dP^1$,
denote by $\dC_0=\Spec\dC[y]$ resp. $\dC_\infty=\Spec\dC[x]$ the standard charts of $\dP^1$ at zero
and infinity, and write $\widetilde{\dC}_0:=\dC\times\dC_0$ resp. $\widetilde{\dC}_\infty:=\dC\times\dC_\infty$.
Then $\cbl=\widetilde{\dC}_0\cup_{\dC\times\dC^*}\widetilde{\dC}_\infty$, and the restrictions
of the universal family to the charts are
$$
\begin{array}{rcl}
\cH_0:=\cH_{|\dC\times\widetilde{\dC}_0}:=\cO_{\dC\times\widetilde{\dC}_0}\left[z^{-1/2}(A+y\overline{A})+r_0z^{3/2}\overline{A}\right]\oplus\cO_{\dC\times\widetilde{\dC}_0}z^{5/2}\overline{A} \\ \\
\cH_\infty:=\cH_{|\dC\times\widetilde{\dC}_\infty}:=\cO_{\dC\times\widetilde{\dC}_\infty}\left[z^{-1/2}(xA+\overline{A})+r_\infty z^{3/2}A\right]\oplus\cO_{\dC\times\widetilde{\dC}_\infty}z^{5/2}A
\end{array}
$$
On the intersection $\widetilde{\dC}^*:=\widetilde{\dC}_0\cap\widetilde{\dC}_\infty$, we have
$(\cH_0)_{|\dC\times\widetilde{\dC}^*}=
\cO_{\dC\times\widetilde{\dC}^*}\left[z^{-1/2}(y^{-1}A+\overline{A})+y^{-1}r_0 z^{3/2}\overline{A}\right]\oplus\cO_{\dC\times\widetilde{\dC}^*}z^{5/2}\overline{A}$,
which is equal to $(\cH_\infty)_{|\dC\times\widetilde{\dC}^*}$ iff $r_\infty=-r_0 x^2$:
write $v_1:=z^{-1/2}(y^{-1}A+\overline{A})+y^{-1}r_0 z^{3/2}\overline{A}$ and $v_2:=z^{5/2}\overline{A}$, then
\begin{equation}
\label{eqDegreeCalc}
\begin{array}{rcl}
(\cH_0)_{|\dC\times\widetilde{\dC}^*}& = &
\cO_{\dC\times\widetilde{\dC}^*}\left[v_1(1-z^2r_0y^{-1})+r^2_0y^{-2}z v_2\right]\oplus\cO_{\dC\times\widetilde{\dC}^*}z^{5/2}\overline{A}\\ \\
& = & \cO_{\dC\times\widetilde{\dC}^*}\left[z^{-1/2}(xA+\overline{A}) - (r_0 x^2) z^{3/2}A\right]\oplus\cO_{\dC\times\widetilde{\dC}^*}z^{5/2}A.
\end{array}
\end{equation}
We obtain the following result.
\begin{proposition}
The classifying space $\MBLred$ associated to the above topological data and
the spectral range $\alpha_1=-\frac12$, $w=2$ is the Hirzebruch surface $\Sigma_2$.
\end{proposition}
We also describe the $tt^*$-geometry on $\MBLred$: On the chart $\Spec\dC[r_0,y]$, the
locus where $\widehat{\cH}$ is non-pure is the real-analytic hypersurface given by
$(1-|y|^2)(|r_0|^2-(1-|y|^2)^2)=0$, the complement has four connected components,
on two of them $\widehat{\cH}$ is polarized, on the other two components
the metric is negative definite.
This is visualized in the following picture, where the thickened lines represent the
pure polarized limit TERP-structures $\cG$.
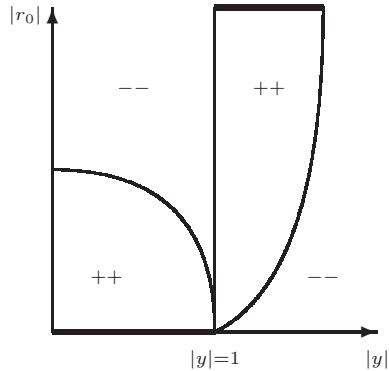
\begin{center}
\setlength{\unitlength}{0.36cm}
\begin{picture}(36,16)
\thicklines

\put(12,2){\vector(1,0){12}}
\put(12,2){\vector(0,1){12}}

\qbezier(12,8)(18,8)(18,2)
\put(18,2){\line(0,1){12}}
\qbezier(18,2)(22,4)(22,14)

\put(24,1){\makebox(0,0){$\scriptstyle |y|$}}
\put(11,13.7){\makebox(0,0){$\scriptstyle |r_0|$}}

\put(18,1){\makebox(0,0){$\scriptstyle |y|=1$}}

\put(14,4){\makebox(0,0){$\scriptstyle ++$}}
\put(15,11){\makebox(0,0){$\scriptstyle --$}}
\put(20,11){\makebox(0,0){$\scriptstyle ++$}}
\put(22,4){\makebox(0,0){$\scriptstyle --$}}

\linethickness{0.7mm}
\put(12,2){\line(1,0){6}}
\put(18,14){\line(1,0){4}}

\end{picture}
\end{center}

\subsection{Weighted projective spaces}
\label{subsecWPS}
We already encountered a weighted projective space
as the compactification of a (non-maximal) stratum
in a space $\MBLred$ which was itself smooth.
In the following two examples the whole compactification $\MBLred$
will be isomorphic to some weighted projective spaces.
After lemma \ref{t9.2}, which will conclude the discussion of $\MBLred$
and the variation of twistor structures of the first example,
we will come to $\MBL$ and find $\MBL\neq \MBLred$.

Consider a three-dimensional real vector space $H^\infty_\dR$,
its complexification $H^\infty:=H^\infty_\dR\otimes \dC$, choose a basis
$H^\infty=\oplus_{i=1}^3\dC A_i$ such that $\overline{A}_1=A_3$, $A_2\in H^\infty_\dR$.
Choose a real number $\alpha_1\in(-3/2,-1)$, put $\alpha_2:=0$, $\alpha_3:=-\alpha_1$
and let $M\in\mathit{Aut}(H^\infty_\dC)$ be given
by $M(\underline{A})=\underline{A}\cdot\diag(\lambda_1,\lambda_2,\lambda_3)$
where $\underline{A}:=(A_1,A_2,A_3)$ and $\lambda_i:=e^{-2\pi i \alpha_i}$ (then $M$ is actually an element in $\mathit{Aut}(H^\infty_\dR)$).
Let $(H',H_\dR',\nabla)$ be the flat holomorphic bundle on $\dC^*\times\dC^2$ with real flat subbundle corresponding
to $(H^\infty,H^\infty_\dR,M)$, and put $s_i:=z^{\alpha_i}A_i\in\cH'$.
Moreover, define the pairing $P:\cH'\otimes j^*\cH'\rightarrow \cO_{\dC^*\times \dC^2}$
by $P(\underline{s}^{tr}, \underline{s}):=(\delta_{i+j,4})_{i,j\in\{1,\ldots,3\}}$

Denote by $r,t$ coordinates on $\dC^2$,
and define $\cH:=\oplus_{i=1}^3 \cO_{\dC^3} v_i$, where
$$
\begin{array}{rcl}
v_1 & := & s_1 + r z^{-1} s_2 + \frac{r^2}{2} z^{-2} s_3 + t z^{-1} s_3\\
v_2 & := & s_2 + r z^{-1} s_3\\
v_3 & := & s_3
\end{array}
$$

Let $w:=0$, then it can be checked by direct calculations that
$\TERP$ is a variation of regular singular TERP-structures on $\dC^2$.
Moreover, the Hodge filtration induced on $H^\infty$ is constant in $r$ and $t$ and
gives a sum of pure polarized Hodge structures of weights $0$ and $-1$ on $H^\infty_1$ and $H^\infty_{\neq 1}$, namely, we
have that
$$
\{0\}=F^2 \subsetneq
F^1:=\dC A_1 \subsetneq
F^0 := \dC A_1 \oplus \dC A_2 =
F^{-1} \subsetneq
F^{-2}:= H^\infty
$$
The polarizing form $S$ is given by $P$ via \cite[formulas (5.4) and (5.5)]{HS1}:
$$
S(\underline{A}^{tr},\underline{A}):=
\begin{pmatrix}
0 & 0 & \gamma \\
0 & 1 & 0 \\
-\gamma & 0 &0
\end{pmatrix},
$$
where $\gamma:=\frac{-1}{2\pi i}\Gamma(\alpha_1+2)\Gamma(\alpha_3-1)$. In particular,
we have for $p=1$
$$
i^{p-(-1-p)}S(A_1,A_3)= (-1)iS(A_1,A_3)= \frac{\Gamma(\alpha_1+2)\Gamma(\alpha_3-1)}{2\pi} >0
$$
and for $p=0$
$$
i^{p-(-p)}S(A_2,A_2)=S(A_2,A_2)>0
$$
so that $F^\bullet$ indeed induces a pure polarized Hodge structure of weight $-1$ on $H^\infty_{\neq 1}=\dC A_1\oplus\dC A_2$
and a pure polarized Hodge structure of weight $0$ on $H^\infty_1=\dC A_2$.
This situation is shown in the following diagram.

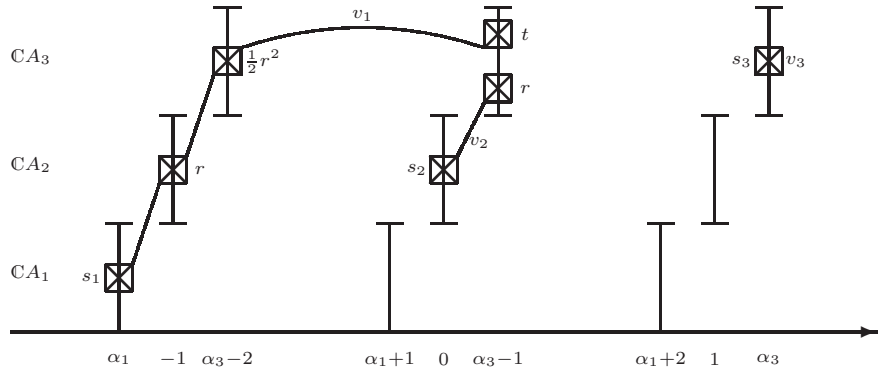
\begin{center}
\setlength{\unitlength}{0.36cm}
\begin{picture}(30,16)
\thicklines

\put(0,2){\vector(1,0){32}}

\put(0,4){\makebox{$\scriptstyle \dC A_1$}}
\put(0,8){\makebox{$\scriptstyle \dC A_2$}}
\put(0,12){\makebox{$\scriptstyle \dC A_3$}}

\multiput(4,2)(10,0){3}{\line(0,1){4}}

\multiput(6,6)(10,0){3}{\line(0,1){4}}

\multiput(8,10)(10,0){3}{\line(0,1){4}}

\multiput(3.5,6)(2,0){2}{\line(1,0){1}}
\multiput(5.5,10)(2,0){2}{\line(1,0){1}}

\multiput(13.5,6)(2,0){2}{\line(1,0){1}}
\multiput(15.5,10)(2,0){2}{\line(1,0){1}}

\multiput(23.5,6)(2,0){2}{\line(1,0){1}}
\multiput(25.5,10)(2,0){2}{\line(1,0){1}}

\multiput(7.5,14)(10,0){3}{\line(1,0){1}}

\put(3.5,3.5){\line(1,0){1}}\put(3.5,3.5){\line(0,1){1}}
\put(4.5,3.5){\line(0,1){1}}\put(3.5,4.5){\line(1,0){1}}
\put(3.5,3.5){\line(1,1){1}}\put(3.5,4.5){\line(1,-1){1}}
\put(3,4){\makebox(0,0){$\scriptstyle s_1$}}
\qbezier(4.5,4.5)(4.5,4.5)(5.5,7.5)

\multiput(5.5,7.5)(10,0){2}{\line(1,0){1}}\multiput(5.5,7.5)(10,0){2}{\line(0,1){1}}
\multiput(6.5,7.5)(10,0){2}{\line(0,1){1}}\multiput(5.5,8.5)(10,0){2}{\line(1,0){1}}
\multiput(5.5,7.5)(10,0){2}{\line(1,1){1}}\multiput(5.5,8.5)(10,0){2}{\line(1,-1){1}}
\put(15,8){\makebox(0,0){$\scriptstyle s_2$}}
\put(7,8){\makebox(0,0){$\scriptstyle r$}}
\qbezier(6.5,8.5)(6.5,8.5)(7.5,11.5)
\put(13,13.7){\makebox(0,0){$\scriptstyle v_1$}}
\put(17.3,9){\makebox(0,0){$\scriptstyle v_2$}}

\multiput(7.5,11.5)(20,0){2}{\line(1,0){1}}\multiput(7.5,11.5)(20,0){2}{\line(0,1){1}}
\multiput(8.5,11.5)(20,0){2}{\line(0,1){1}}\multiput(7.5,12.5)(20,0){2}{\line(1,0){1}}
\multiput(7.5,11.5)(20,0){2}{\line(1,1){1}}\multiput(7.5,12.5)(20,0){2}{\line(1,-1){1}}
\put(19,13){\makebox(0,0){$\scriptstyle t$}}
\put(9.3,12){\makebox(0,0){$\scriptstyle \frac{1}{2}r^2$}}
\put(27,12){\makebox(0,0){$\scriptstyle s_3$}}
\put(29,12){\makebox(0,0){$\scriptstyle v_3$}}

\qbezier(8.5,12.5)(13,14)(17.5,12.5)

\put(17.5,10.5){\line(1,0){1}}\put(17.5,10.5){\line(0,1){1}}
\put(18.5,10.5){\line(0,1){1}}\put(17.5,11.5){\line(1,0){1}}
\put(17.5,10.5){\line(1,1){1}}\put(17.5,11.5){\line(1,-1){1}}

\put(17.5,12.5){\line(1,0){1}}\put(17.5,12.5){\line(0,1){1}}
\put(18.5,12.5){\line(0,1){1}}\put(17.5,13.5){\line(1,0){1}}
\put(17.5,12.5){\line(1,1){1}}\put(17.5,13.5){\line(1,-1){1}}

\qbezier(16.5,8.5)(16.5,8.5)(17.5,10.5)
\put(19,11){\makebox(0,0){$\scriptstyle r$}}

\put(4,1){\makebox(0,0){$\scriptstyle \alpha_1$}}
\put(6,1){\makebox(0,0){$\scriptstyle -1$}}
\put(8,1){\makebox(0,0){$\scriptstyle \alpha_3-2$}}

\put(14,1){\makebox(0,0){$\scriptstyle \alpha_1+1$}}
\put(16,1){\makebox(0,0){$\scriptstyle 0$}}
\put(18,1){\makebox(0,0){$\scriptstyle \alpha_3-1$}}

\put(24,1){\makebox(0,0){$\scriptstyle \alpha_1+2$}}
\put(26,1){\makebox(0,0){$\scriptstyle 1$}}
\put(28,1){\makebox(0,0){$\scriptstyle \alpha_3$}}
\end{picture}
\end{center}

It is clear that $\DPHS=\DcPHS=\DPMHS=\DcPMHS=\{*\}$
for the given topological data and that the above family is indeed the
universal family, in particular, $D_{BL}=\cbl\cong\dC^2$.

Let us describe the variation of twistors associated to this example.
We have
$$
\tau v_1 :=  s_3+\overline{r}zs_2+\frac{\overline{r}^2}{2}z^2s_1 + \overline{t} zs_1
\quad;\quad
\tau v_2 :=  s_2+\overline{r}zs_1
\quad;\quad
\tau v_3 :=  s_1.
$$
The family of twistors $\widehat{\cH}$ is pure outside of the real-analytic hypersurface
$D:=\{(1-\rho)^4=\theta\}$, where $\rho=\frac12r\overline{r}$ and $\theta=t\overline{t}$.
Namely, write $U_1:=\cbl\backslash (D\cup\{\rho= 1\})$ and $U_2:=\cbl\backslash (D\cup\{\theta =0\})$,
then $\cbl\backslash D = U_1\cup U_2$. We have
$\widehat{\cH}_{|\dP^1\times U_1}=\oplus_{i=1}^3 \cO_{\dP^1}\cAn_{U_1} \, w_i$, where
$$
\begin{array}{rcl}
w_1 & := & s_1 + r z^{-1} s_2 + \frac{r^2}{2} z^{-2} s_3 -\frac{\frac12\overline{r}^2t}{1-\frac{(r\overline{r})^2}{4}}zs_1
+ \frac{t}{1-\frac{(r\overline{r})^2}{4}}z^{-1}s_3 \\ \\
w_2 & := & \overline{r} z s_1 +\left(1+\frac12r\overline{r}\right)s_2 +r z^{-1} s_3\\ \\
w_3 & := & \frac{\overline{t}}{1-\frac{(r\overline{r})^2}{4}}z s_1 -\frac{\frac12 r^2\overline{t}}{1-\frac{(r\overline{r})^2}{4}}z^{-1}s_3
+ \frac{\overline{r}^2}{2}z^2s_1 + \overline{r}zs_2+ s_3 = \tau(w_1)\\
\end{array}
$$
On the other hand, $\widehat{\cH}_{|\dP^1\times U_2}=\oplus_{i=1}^3 \cO_{\dP^1}\cAn_{U_2} \, \widetilde{w}_i$, where
$$
\begin{array}{rcl}
\widetilde{w}_1 & := & s_1 + r z^{-1} s_2 + \frac{r^2}{2} z^{-2} s_3 +tz^{-1}s_3
+\frac{t\overline{r}^2}{2\overline{t}} \left( \frac{\overline{r}^2}{2} z^{2} s_1
+ \overline{r}z s_2 + s_3\right) \\ \\
\widetilde{w}_2 & := & \overline{r} z s_1 +\left(1+\frac12r\overline{r}\right)s_2 +r z^{-1} s_3\\ \\
\widetilde{w}_3 & := & \frac{\overline{t}r^2}{2t}\left( s_1 +  r z^{-1}s_2 + \frac{r^2}{2}z^{-2}s_3\right) +
\overline{t}zs_1
+ \frac{\overline{r}^2}{2}z^2s_1 + \overline{r}zs_2+ s_3 = \tau(\widetilde{w}_1)\\
\end{array}
$$
This shows that $\widehat{\cH}$ is pure precisely outside $D$.
The complement of $D$ has three components. $\widehat{\cH}$ is polarized on two of them,
those which contain $\{(r,0)\ |\ |r|<\sqrt{2}\}$ and $\{(r,0)\ |\ |r|>\sqrt{2}\}$, respectively.
On the third component the metric on $p_*\widehat \cH$ has signature $(+,-,-)$.
This is visualized in the following picture.
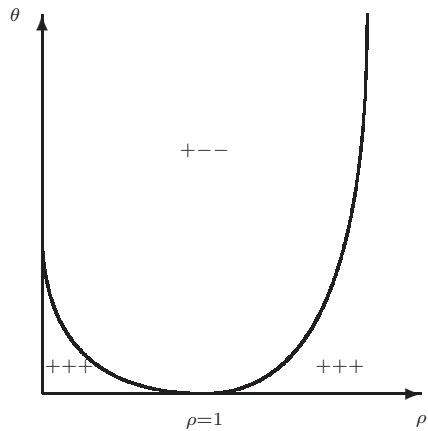
\begin{center}
\setlength{\unitlength}{0.36cm}
\begin{picture}(38,18)
\thicklines

\put(12,2){\vector(1,0){14}}
\put(12,2){\vector(0,1){14}}

\qbezier(12,8)(12,2)(18,2)
\qbezier(18,2)(24,2)(24,16)

\put(26,1){\makebox(0,0){$\scriptstyle \rho$}}
\put(11,16){\makebox(0,0){$\scriptstyle \theta$}}

\put(18,1){\makebox(0,0){$\scriptstyle \rho=1$}}

\put(13,3){\makebox(0,0){$\scriptstyle +++$}}
\put(18,11){\makebox(0,0){$\scriptstyle +--$}}
\put(23,3){\makebox(0,0){$\scriptstyle +++$}}

\end{picture}
\end{center}

We are going to compute the compact space $\MBLred$ as in subsection \ref{subsecSmoothCompact}.
First note that we have $\cH_{|r\neq 0}=\oplus_{i=1}^3\cO_{\dC\times\dC^*\times\dC}\widetilde{v}_i$, where
$$
\begin{array}{rllll}
\widetilde{v}_1 & := & r^{-2}(v_1-\frac{t}{r}v_2+2\frac{t}{r} \widetilde{v}_2)
&=& r^{-2}(s_1 + r z^{-1} s_2) + \frac{1}{2} z^{-2} s_3 + 2\frac{t}{r^4} zs_1\\
\widetilde{v}_2 & := & r^{-1}(zv_1-\frac{r}{2} v_2 - t v_3)  & = & r^{-1} z s_1+\frac12s_2\\
\widetilde{v}_3 & := & z^2v_1-rzv_2+\frac{r^2}{2}v_3-tzv_3 & = & z^2s_1
\end{array}
$$
Consider for any fixed parameter $u\in\dC$ the restriction
$\cL^{(u)}:=\cL_{|\{t-ur^4=0,r\neq 0\}}$. Then we have $\cL^{(u)}=\oplus_{i=1}^3\cO_{\dC\times\dC^*}\widetilde{w}_i$,
where
$$
\begin{array}{rcl}
\widetilde{w}_1 & := & r^{-2}(s_1 + r z^{-1} s_2) + \frac{1}{2} z^{-2} s_3 + 2 u zs_1\\
\widetilde{w}_2 & := & r^{-1} z s_1+\frac12s_2\\
\widetilde{w}_3 & := & z^2s_1
\end{array}
$$
It is clear that this basis defines an extension of $\cL^{(u)}$ to a locally free $\cO_{\dC\times(\dP^1\backslash\{0\})}$-module,
its fibre at $r=\infty$ is given by $\cO_{\dC}(z^{-2}s_3+4uzs_1)\oplus\cO_{\dC} s_2 \oplus\cO_{\dC} z^2 s_1$. These extensions are
non-isomorphic for any two $u_1\neq u_2$, which shows the following result.
\begin{lemma}\label{t9.2}
The space $\MBLred$ for the initial data from above is the
weighted projective space $\dP(1,1,4)$, where the embedding
$\cbl \hookrightarrow \dP(1,1,4)=\Proj\dC[X,Y,Z]$
is given by $X:=r, Y:=1,Z:=t$ (here $\deg(X)=1, \deg(Y)=1, \deg(Z)=4$).
The only singular point of $\dP(1,1,4)$ is $(0:0:1)$.
In particular, $(\MBLp)^{red}$ is smooth in this case.
\end{lemma}
\begin{proof}
The only thing to show is that in this case $\MBLred$ is indeed the closure of
the above classifying space $D_{BL} = \dC^2$. This follows from the fact that the only
possible spectral numbers for the range $[\alpha_1,\alpha_\mu]$ are
$(\alpha_1,\alpha_2,\alpha_3)$, $(\alpha_3-2,\alpha_2,\alpha_1+2)$ and
$(\alpha_1+1,\alpha_2,\alpha_3-1)$. The respective strata are
$U_{(\alpha_1,\alpha_2,\alpha_3)}=D_{BL}=\Spec\dC[r,t]$,
$U_{(\alpha_3-2,\alpha_2,\alpha_1+2)}=\Spec\dC[u]$ and
$U_{(\alpha_1+1,\alpha_2,\alpha_3-1)}=(0:0:1)\in\dP(1,1,4)$
which are all the possible strata in $\dP(1,1,4)$.
\end{proof}

We will determine precisely the part of $\MBL$ underlying the affine
chart of $\MBLred\cong \dP(1,1,4)$ with coordinates $(r,t)$
and make remarks about the other two standard charts.
We follow the Ansatz in the proof of theorem \ref{t7.6}.
A priori here we need nine coordinates in the Ansatz,
\begin{eqnarray*}
v_1&=& s_1 + r\cdot z^{-1}s_2 + r_2\cdot z^{-2}s_3 + t\cdot z^{-1}s_3,\\
v_2&=& s_2 + \varepsilon\cdot z^{-1}s_2 + r_3\cdot z^{-2}s_3 + r_4\cdot z^{-1}s_3,\\
v_3&=& s_3 + r_5\cdot z^{-1}s_2 + r_6\cdot z^{-2}s_3 + r_7\cdot z^{-1}s_3.
\end{eqnarray*}
The pairing $P$ gives the following seven equations,
\begin{eqnarray*}
0&=& P^{(w-2)}(v_1,v_3)=r_6,\\
0&=& P^{(w-1)}(v_2,v_3)=-r_5,\\
0&=& P^{(w-1)}(v_1,v_3)=-r_7,\\
0&=& P^{(w-2)}(v_1,v_2)=r_3,\\
0&=& P^{(w-1)}(v_1,v_2)=r-r_4,\\
0&=& P^{(w-2)}(v_2,v_2)=-\varepsilon^2,\\
0&=& P^{(w-2)}(v_1,v_1)=2r_2-r^2.
\end{eqnarray*}
These show
\begin{eqnarray*}
&&r_6=r_5=r_7=r_3=0,\ r_4=r,\ r_2=\frac{r^2}{2},\ \varepsilon^2=0,\\
v_1&=& s_1 + r\cdot z^{-1}s_2 + \frac{r^2}{2}\cdot z^{-2}s_3 + t\cdot z^{-1}s_3,\\
v_2&=& s_2 + \varepsilon\cdot z^{-1}s_2 + r\cdot z^{-1}s_3,\\
v_3&=& s_3.
\end{eqnarray*}
The pole of order 2 gives nothing from $v_2$ and $v_3$, but one equation from $v_1$:
\begin{eqnarray*}
z^2\nnn_z v_3 &=& \alpha_3 \cdot zv_3,\\
z^2\nnn_z v_2 &=& -\varepsilon \cdot zs_2 = -\varepsilon\cdot zv_2,\\
z^2\nnn_z v_1 &=& \alpha_1\cdot zs_1 + (-1)r\cdot s_2 + (\alpha_3-2)\frac{r^2}{2}\cdot z^{-1}s_3
+ (\alpha_3-1)t\cdot s_3\\
&=& \alpha_1\cdot zv_1 + (-1-\alpha_1)r\cdot v_2 + (\alpha_3-1)t\cdot v_3
+ (1+\alpha_1)r\cdot\varepsilon\cdot( z^{-1}s_2 + r\cdot z^{-2}s_3),\\
\textup{thus} & & r\cdot \varepsilon=0.
\end{eqnarray*}
We recover the variation of TERP-structures in $(r,t)$, but additionally there is
an obstructed deformation on the line $\{r=0\}$ with the parameter $\varepsilon$ with $\varepsilon^2=0$
and $r\cdot\varepsilon=0$. Obviously it does not preserve the spectrum.

On the affine chart with coordinates $(\www r,t)=(\frac{1}{r},t)$ one finds exactly the
same behavior, on the line $\{\www r=0\}$ there is an obstructed deformation with
a parameter $\www\varepsilon$ with $\www\varepsilon^2=0$ and $\www r\cdot\www\varepsilon=0$.
On the affine chart around $(0:0:1)\in \dP(1,1,4)\cong \MBLred$ one obtains with some more work
six coordinates and nine quadratic equations. The Zariski tangent spaces at $(0:0:1)$
satisfy
$$
\dim T_{(0:0:1)}\MBL = 6>5=\dim T_{(0:0:1)}\MBLred,
$$
so also at $(0:0:1)$ the canonical and the reduced complex structure differ.

\vspace*{1cm}

The next example also gives the weighted projective space $\dP(1,1,4)$ as
the final result, but in a completely different way, namely, the classifying space
$\cbl$ of the generic spectrum is a line bundle over $\dP^1$ of weight $4$,
and the compactification $\MBLred$ is obtained by adding a single point.

Let $H^\infty_\dR:=\oplus_{i=1}^3\dR B_i$,
$S(\underline{B}^{tr},\underline{B})=\diag(\frac12,1,\frac12),
$
$A_1:=(B_1+iB_3),A_2:=B_2,A_3:=(B_1-iB_3)$, $\alpha_1:=-1$ and $w=0$. Define $F^\bullet_0$ by
$$
\{0\} = F^2_0 \subsetneq F_0^1:=\dC A_1\subsetneq F_0^0:=\dC A_1\oplus \dC A_2\subsetneq F_0^{-1}=H^\infty.
$$
Then $(H^\infty, H_\dR^\infty, S, F^\bullet)$ is a pure Hodge structure of weight zero,
however, the Hodge metric
$$
h:=-S(\cdot,\overline{\cdot})_{|H_0^{1,-1}\times H_0^{1,-1}}\oplus\, S(\cdot,\overline{\cdot})_{|H_0^{0,0}\times H_0^{0,0}}
\oplus\,-S(\cdot,\overline{\cdot})_{|H_0^{-1,1}\times H_0^{-1,1}}
$$
(where $H_0^{p,-p}=F_0^p \cap \overline{F}_0^{-p}$) has signature $(1,2)$.
Consider the classifying space $\DcPMHS$ of all filtrations
$F^\bullet$ on $H^\infty$ satisfying $S(F^p,F^{1-p})=0$ and having the same Hodge numbers
as $F^\bullet_0$. Such a filtration is uniquely determined by $F^1$, for
$F^0$ is necessarily the $S$-orthogonal complement of $F^1$ in $H^\infty$. It
follows that $F^1$ must satisfy the isotropy condition $S(F^1,F^1)=0$. This is
the defining equation for a plane quadric $Q\subset V(a^2+c^2+2b^2)\subset\Proj\dC[a,b,c]=\Gr(1,3)$, where each point defines
$F^1:=\dC \widetilde{A}_1:=\dC(a B_1 + b B_2 + c B_3)$, $F^0:=(F^1)^{\bot,S}$.
We conclude that $\DcPMHS\cong Q$. The equation $a^2+c^2+2b^2=0$ has no real
solutions other than $(0,0,0)$, so that for any $F^1\in Q$, $(F^1)^{\bot,S}\cap \overline{F}^1=\{0\}$. This means that
$(\{0\}\subsetneq F^1 \subsetneq F^0:=(F^1)^{\bot,S}\subsetneq F^{-1} = H^\infty)$ is pure
with signature $(1,2)$. Thus in this case $\emptyset = \DPMHS\subsetneq \DcPMHS = Q \cong \dP^1$.

Consider the classifying space $\cbl$ associated to these given initial data.
The fibration $\cbl\rightarrow \DcPMHS$ has one-dimensional
affine fibres with a $\dC^*$-action, hence it is again the total space $\mathbf{V}(\cE)$ of
a line bundle $\cE$ on $\dP^1$: For fixed $F^\bullet=
(\{0\}\subsetneq F^1:=\dC \widetilde{A}_1\subsetneq F^0=\dC \widetilde{A}_1\oplus \dC \widetilde{A}_2\subsetneq H^\infty)$
in $\DcPMHS$, the fibre $\check{\pi}_{BL}^{-1}(F^\bullet)$
is given by $\cH:=\oplus_{i=1}^3 \cO_{\dC^2}v_i$, where
$$
v_1 = s_1 + r z ^{-1} s_3 \quad ; \quad
v_2 = s_2 \quad ; \quad
v_3 = s_3
$$
and $s_1:=z^{-1}\widetilde{A}_1$, $s_1:=\widetilde{A}_2$ and $s_3:=z \overline{\widetilde{A}}_1$. If $r$ tends to infinity,
we have $\lim_{r\rightarrow\infty}\cH=\oplus_{i=1}^3 \cO_\dC A_i$. The diagram
for this situation is as follows:
\begin{center}
\setlength{\unitlength}{0.36cm}
\begin{picture}(26,14)
\thicklines

\put(0,2){\vector(1,0){26}}

\put(-2,4){\makebox{$\scriptstyle \dC A_1$}}
\put(-2,8){\makebox{$\scriptstyle \dC A_2$}}
\put(-2,12){\makebox{$\scriptstyle \dC A_3$}}

\multiput(2,2)(11,0){3}{\line(0,1){12}}
\multiput(1,6)(11,0){3}{\line(1,0){2}}
\multiput(1,10)(11,0){3}{\line(1,0){2}}

\multiput(1.5,3.5)(11,4){3}{\line(1,0){1}}\multiput(1.5,3.5)(11,4){3}{\line(0,1){1}}
\multiput(2.5,3.5)(11,4){3}{\line(0,1){1}}\multiput(1.5,4.5)(11,4){3}{\line(1,0){1}}
\multiput(1.5,3.5)(11,4){3}{\line(1,1){1}}\multiput(1.5,4.5)(11,4){3}{\line(1,-1){1}}

\put(12.5,11.5){\line(1,0){1}}\put(12.5,11.5){\line(0,1){1}}
\put(13.5,11.5){\line(0,1){1}}\put(12.5,12.5){\line(1,0){1}}
\put(12.5,11.5){\line(1,1){1}}\put(12.5,12.5){\line(1,-1){1}}

\qbezier(2.5,4.5)(2.5,4.5)(12.5,11.5)

\put(1,4){\makebox(0,0){$\scriptstyle s_1$}}
\put(7,7){\makebox(0,0){$\scriptstyle v_1$}}
\put(12,8){\makebox(0,0){$\scriptstyle s_2$}}
\put(14,8){\makebox(0,0){$\scriptstyle v_2$}}
\put(23,12){\makebox(0,0){$\scriptstyle s_3$}}
\put(25,12){\makebox(0,0){$\scriptstyle v_3$}}
\put(14,12){\makebox(0,0){$\scriptstyle r$}}

\put(1.7,1){\makebox(0,0){$\scriptstyle -1$}}
\put(13,1){\makebox(0,0){$\scriptstyle 0$}}
\put(24,1){\makebox(0,0){$\scriptstyle 1$}}
\end{picture}
\end{center}
\begin{lemma}
The space $\MBLred$ for the topological data $(H^\infty,H_\dR^\infty,S,M=\mathit{id},
\alpha_1=-1,w=0)$ from above is the blow-down of the $\infty$-section of the Hirzebruch surface
$\Sigma_4=\mathbf{Proj}\,(\cO_{\dP^1}(4)\oplus\cO_{\dP^1})$, i.e., it is isomorphic to the weighted projective
space $\dP(1,1,4)$.
\end{lemma}
\begin{proof}
The degree of $\cE$ is seen to be four by a calculation similar (although more complicated)
to equation \eqref{eqDegreeCalc}, where the biregular parametrization
$$
\begin{array}{rcl}
\dP^1 & \longrightarrow & Q \\
(x:y) & \longmapsto & (x^2+y^2:i\sqrt{2}xy:i(y^2-x^2))
\end{array}
$$
is used. On the other hand, it is directly evident that
$\lim_{r\mapsto \infty} \cH = \cO_\dC z s_1 \oplus \cO_\dC s_2 \oplus \cO_\dC z^{-1} s_3
=V^0=\cO_\dC H^\infty$. In particular, for any two $F^\bullet_1,F^\bullet_2,\in \DcPMHS$,
and $\cH_i(r)\in\check{\pi}_{BL}^{-1}(F^\bullet_i)$, $\lim_{r\rightarrow \infty}\cH_1(r)=\lim_{r\rightarrow \infty}\cH_2(r)$.
Geometrically, this means that the fibration $\check{\pi}_{BL}:\mathbf{V}(\cE)\rightarrow \dP^1$
is compactified to $\Sigma_4\rightarrow \dP^1$, and $\MBLred$ is obtained by blowing down the $\infty$-section
of $\Sigma_4\rightarrow \dP^1$. Moreover, it is known (see, e.g., \cite{Dolg}) that the blow-up of
the singular points of $\dP(1,1,n)$ is exactly $\Sigma_n$.

We remark that it is also possible to calculate directly the local structure of $\MBLred$ in a neighborhood
of $\cG:=\cO_\dC H^\infty$, which yields $\cO_{(\MBLred,\cG)}\cong \dC\{a^4,a^3b,a^2b^2,ab^3,b^4\}$.
\end{proof}

\subsection{Reducible spaces}\label{subsecRedSpaces}
The following example shows that $\MBL$ and $\MBLred$
might have several components. Fix any number $n\in\dN_{>0}$ and
consider the topological data: $H^\infty_\dR:=\dR B_1 \oplus \dR B_2$, $M:=\Id\in\Aut(H^\infty_\dR)$,
and $S(\underline{B}^{tr},\underline{B})=(-1)^n\frac12\diag(1,1)$. Let $A_1:=B_1+iB_2$, $A_2:=\overline{A}_1=B_1-iB_2$
so that $S(A_i,A_j)=(-1)^n\delta_{i+j,3}$ and consider the reference filtration
$$
\{0\}= F_0^{n+1} \subsetneq F_0^n:=\dC A_1 = F_0^{n-1} = \ldots = F_0^{-n+1}\subsetneq F_0^{-n} = H^\infty.
$$
Let $w=0$, $\alpha_1=-n$, then
the classifying space is $\DcPMHS=\DPMHS=\DcPHS=\DPHS$,
it consists of two points, namely $F^\bullet_0$ and
$\overline{F}^\bullet_0$, which are both pure polarized Hodge structures of
weight zero with Hodge decomposition $H^{n,-n}\oplus H^{-n,n}$.
Put $s_1:=z^{-n}A_1$, $s_2:=z^n A_2$ then
\cite[formulas (5.4), (5.5)]{HS1} yields $P(s_i,s_j)=\delta_{i+j,3}$.
The classifying space $\cbl$ is a disjoint union of two affine lines, namely, the universal family
over the component above $F_0^\bullet$ is given by
$\cH_{-n}(r):=\cO_{\dC^2} v_1\oplus \cO_{\dC^2} v_2$, where
$v_1:=s_1+rz^{-1}s_2, v_2:=s_2$, and
the universal family over the other component is
$\cH_n(r):=\cO_{\dC^2} (z^{-n}A_2+rz^{n-1}A_1)\oplus \cO_{\dC^2} z^n A_1$.
The following diagram visualizes this situation.
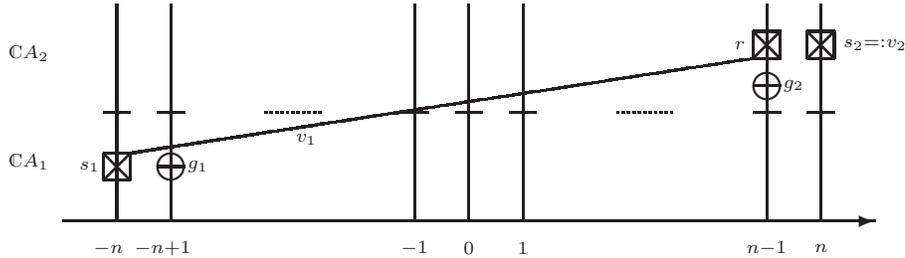
\begin{center}
\setlength{\unitlength}{0.36cm}
\begin{picture}(30,12)
\thicklines

\put(0,2){\vector(1,0){30}}

\put(-2,4){\makebox{$\scriptstyle \dC A_1$}}
\put(-2,8){\makebox{$\scriptstyle \dC A_2$}}

\multiput(2,2)(2,0){2}{\line(0,1){8}}
\multiput(13,2)(2,0){3}{\line(0,1){8}}
\multiput(13,2)(2,0){3}{\line(0,1){8}}
\multiput(26,2)(2,0){2}{\line(0,1){8}}

\multiput(1.5,6)(2,0){2}{\line(1,0){1}}
\multiput(12.5,6)(2,0){3}{\line(1,0){1}}
\multiput(25.5,6)(2,0){2}{\line(1,0){1}}
\multiput(7.5,6)(13,0){2}{\dashbox{0.1}(2,0.02){}}

\put(1.5,3.5){\line(1,0){1}}\put(1.5,3.5){\line(0,1){1}}
\put(2.5,3.5){\line(0,1){1}}\put(1.5,4.5){\line(1,0){1}}
\put(1.5,3.5){\line(1,1){1}}\put(1.5,4.5){\line(1,-1){1}}

\put(25.5,8){\line(1,0){1}}\put(25.5,8){\line(0,1){1}}
\put(26.5,8){\line(0,1){1}}\put(25.5,9){\line(1,0){1}}
\put(25.5,8){\line(1,1){1}}\put(25.5,9){\line(1,-1){1}}

\put(4,4){\circle{1}}\put(3.5,4){\line(1,0){1}}
\put(26,7){\circle{1}}\put(25.5,7){\line(1,0){1}}

\put(27.5,8){\line(1,0){1}}\put(27.5,8){\line(0,1){1}}
\put(28.5,8){\line(0,1){1}}\put(27.5,9){\line(1,0){1}}
\put(27.5,8){\line(1,1){1}}\put(27.5,9){\line(1,-1){1}}

\qbezier(2.5,4.5)(2.5,4.5)(25.5,8)

\put(1,4){\makebox(0,0){$\scriptstyle s_1$}}
\put(9,5){\makebox(0,0){$\scriptstyle v_1$}}
\put(5,4){\makebox(0,0){$\scriptstyle g_1$}}

\put(25,8.5){\makebox(0,0){$\scriptstyle r$}}
\put(30,8.5){\makebox(0,0){$\scriptstyle s_2=:v_2$}}
\put(27,7){\makebox(0,0){$\scriptstyle g_2$}}

\put(1.7,1){\makebox(0,0){$\scriptstyle -n$}}
\put(3.7,1){\makebox(0,0){$\scriptstyle -n+1$}}
\put(13,1){\makebox(0,0){$\scriptstyle -1$}}
\put(15,1){\makebox(0,0){$\scriptstyle 0$}}
\put(17,1){\makebox(0,0){$\scriptstyle 1$}}
\put(26,1){\makebox(0,0){$\scriptstyle n-1$}}
\put(28,1){\makebox(0,0){$\scriptstyle n$}}
\end{picture}
\end{center}
It is directly evident, that the ``limit TERP''-structure
(when $r$ approaches infinity), is given as $\cG_{-n+1}:=\cO_\dC g_1\oplus
\cO_\dC g_2$, where $g_1:=z s_1$ and $g_2:=z^{-1}s_2$. We see that $\cG_{-n+1}$
is the origin in one of the two components of the stratum $U_{(-n+1,n-1)}$.
The closure of this stratum is the classifying space associated to the
same topological data and to $\alpha_1=-n+1$, so that we get
$U_{(-n+1,n-1)}\cong\dC\coprod\dC$. Note however that
the two (conjugate) filtrations induced by $\cG_{-n+1}$ and $\cG_{n-1}$, respectively,
are not pure polarized: the Hodge metric is negative definite.
Taking the limit of the universal family for this classifying spaces yields
TERP-structures $\cG_{-n+2}$ and $\cG_{n-2}$, respectively, and we can continue
this procedure until we arrive at $\cG_{-1}$ and $\cG_1$.
The limits $\lim_{r\rightarrow \infty}\cH_{-1}(r)=\lim_{r\rightarrow \infty}\cH_1(r)$ are
both equal to the lattice $\cG_0=V^0$.
This shows that the space $\MBLred$ is a chain of $2n$ copies of $\dP^1$, where the Hodge filtration
gives pure polarized resp. negative definite pure Hodge structures on every other component of this chain.
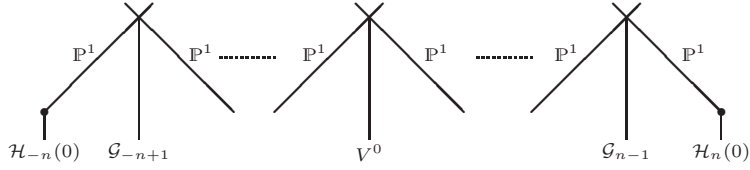
\begin{center}
\setlength{\unitlength}{0.36cm}
\begin{picture}(29,8)
\thicklines
\multiput(2,2)(18,0){2}{\line(1,1){4}}
\multiput(3.5,4.2)(22.5,0){2}{\makebox(0,0){$\scriptstyle \dP^1$}}
\multiput(5,6)(18,0){2}{\line(1,-1){4}}
\multiput(7.8,4.2)(13.5,0){2}{\makebox(0,0){$\scriptstyle \dP^1$}}
\put(10.5,2){\line(1,1){4}}
\put(12,4.2){\makebox(0,0){$\scriptstyle \dP^1$}}
\put(13.5,6){\line(1,-1){4}}
\put(16.5,4.2){\makebox(0,0){$\scriptstyle \dP^1$}}

\multiput(8.5,4)(9.5,0){2}{\dashbox{0.1}(2,0.02){}}
\thinlines
\qbezier(14,5.5)(14,5.5)(14,1)
\put(14,0.5){\makebox(0,0){$\scriptstyle V^{\scriptscriptstyle 0}$}}

\qbezier(5.5,5.5)(5.5,5.5)(5.5,1)
\put(5.5,0.5){\makebox(0,0){$\scriptstyle \cG_{\scriptscriptstyle -n+1}$}}

\qbezier(23.5,5.5)(23.5,5.5)(23.5,1)
\put(23.5,0.5){\makebox(0,0){$\scriptstyle \cG_{\scriptscriptstyle n-1}$}}

\put(2,2){\circle*{0.3}}
\qbezier(2,2)(2,2)(2,1)
\put(2,0.5){\makebox(0,0){$\scriptstyle \cH_{\scriptscriptstyle -n}(0)$}}
\put(27,2){\circle*{0.3}}
\qbezier(27,2)(27,2)(27,1)
\put(27,0.5){\makebox(0,0){$\scriptstyle \cH_{\scriptscriptstyle n}(0)$}}
\end{picture}
\end{center}

It is easy to calculate the associated twistors: For the original family $\cH_{-n}(r)$,
we have
$$
(\widehat{\cH}_{-n})_{||r|\neq 1}:=\cO_{\dP^1}\cAn_{\dC\backslash\{|r|=1\}}(\underbrace{s_1+r z^{-1} s_2}_{w_1})\oplus
\cO_{\dP^1}\cAn_{\dC\backslash\{|r|=1\}}(\underbrace{s_2+\overline{r} z s_1}_{w_2})
$$
and the metric is $h(w_1,w_2)=\diag(1-|r|^2)$, so that $\cH_{-n}(r)$ is a variation
of pure polarized TERP-structures on $\Delta^*$. If $r$ tends to zero,
it degenerates to a twistor generated by elementary sections which corresponds to the
pure polarized Hodge structure $(H^\infty,H^\infty_\dR, S, F^\bullet_0)$.
A similar statement holds for the family $\cH_n(r)$ which degenerates
to a twistor corresponding to $(H^\infty,H^\infty_\dR, S, \overline{F}^\bullet_0)$.

Note however that due to $P(g_1,g_2)=-1$, the variation of twistors on the second left- or rightmost $\dP^1$
is pure polarized on $\dP^1\backslash\overline{\Delta}$, where the origin is
the TERP-structure $\cG_{-n+1}$ (resp. $\cG_{n-1}$).
This means that in the above picture, the points of intersection on the lower level are pure polarized,
but not those on the upper level. In particular, $V^0$ is pure polarized precisely if $n$ is even (the above picture
already supposes that $n$ is odd), which can also be seen directly from $S(A_1,A_2)=(-1)^n$.

\textbf{Remark:} There is a common generalization of this example and the second one from
subsection \ref{subsecSmoothCompact} (where $\MBLred=\Sigma_2$), namely, if we consider the same topological data as in
\ref{subsecSmoothCompact}, but allow a larger spectral range: we put $\alpha_1:=-k-\frac12$
for some $k>0$ and, as before, $w=2$. Then $\cbl$ is still $\mathbf{V}(\cE)$ with $\cE\in\textup{Pic}(\dP^1)$, but
$\MBLred=\coprod_{k} (\Sigma_2)^{(k)}$ is a union of copies of $\Sigma_2$, which are glued along the
zero resp. infinity section of two successive such copies.

\subsection{Monodromy with Jordan block}
\label{subsecJordanBlock}

The following example has a geometric realization within the
1-parameter $\mu$-constant families of hyperbolic singularities
$T_{pqr}$. It is of rank two, and contrary to all the previous
examples, the monodromy $M$ is not semi-simple, but has a $2\times 2$-Jordan block.
We have $\check D_{BL}=\DcPMHS\cong \dC \neq \DcPHS=\{pt\}$.

Set $w=0$, $\alpha_1=-\frac{1}{2}$, $\alpha_2=\frac{1}{2}$, $H^\infty_\dR = \dR A_1\oplus \dR A_2$
and define $M\in\Aut(H^\infty_\dR)$ by
$M(A_1)=-A_1-A_2,\ M(A_2)=-A_2$, so that $N(A_1)=A_2,\  N(A_2)=0$. Moreover,
define the anti-symmetric form $S$ by
$$
S(A_i,A_i)=0, \ S(A_1,A_2)=-1
$$
and let $s_1=i\cdot z^{-\frac{1}{2}-\frac{N}{2\pi i}}A_1$, $s_2=z^{\frac{1}{2}}A_2$ which implies that
$\tau(s_1)=-z\cdot s_1$ and $\tau(s_2)=z^{-1}\cdot s_2$.

Using the relation between $S$ and $P$ from \cite[formulas 5.4, 5.5]{HS1} one calculates
$P(s_i,s_i)=0$ and $P(s_1,s_2)=(-2)\cdot S(A_1,A_2)=2$. The universal family on $\cbl\cong\dC$ is
given as $\cH:=\cO_{\dC^2} v_1\oplus\cO_{\dC^2} v_2$, where
$$
v_1 = s_1+r\cdot z^{-1}s_2\quad\textup{and}\quad v_2 = s_2.
$$
Note that both $v_1$ and $v_2$ are elementary sections, and
$r$ is a parameter on $\DcPMHS$, not on the fibres of
$\cbl\rightarrow \DcPMHS$ (which are single points in this case).
$\cH$ extends to a variation over $\dP^1$ where the fibre over $r=\infty$ is given by
$\cH(\infty)= \cO_\dC\cdot z^{-1}\cdot s_2\oplus \cO_\dC\cdot zs_1$.
It has constant spectrum $\Sp=(-\frac{1}{2},\frac{1}{2})$,
but the spectral pairs jump at $r=\infty$.

$\cH$ is pure outside $\{\Re(r)=0\}\cup\{\infty\}$.
For $\Re(r)\neq 0$ the space $H^0(\dP^1,\widehat{\cH}(r))$ is generated by
$v_1$ and $\tau(v_1)$.
For $\Re(r)>0$ the TERP-structure $\cH(r)$ is pure and polarized,
for $\Re(r)<0$ it is pure with negative definite metric $h$.

For $r\in\dC$ the data $(H^\infty,H^\infty_\dR,F^\bullet,S,-N)$ form
a PMHS of weight $-1$. Here
\begin{eqnarray*}
H^\infty &=& W_0 \supsetneq W_{-1}=W_{-2}=\dC\cdot A_2\supsetneq W_{-3}=\{0\},\\
H^\infty &=& F^{-1} \supsetneq F^0=\dC\cdot
(iA_1+r\cdot A_2) \supsetneq F^1=\{0\},\\ \\
H^{0,0}&=&\dC\cdot[A_1] = W_0/W_{-1},
\quad H^{-1,-1}=\dC\cdot A_2 = W_{-2}=W_{-2}/W_{-3},\\
&&i^{0-0}S([A_1],-N(\overline{[A_1]}))=S(A_1,-A_2)=1>0.
\end{eqnarray*}
For $\Re(r)>0$ it is simultaneously also a PHS of weight $-1$,
$$
i^{0-(-1)}S(iA_1+rA_2,-iA_1+\overline{r}A_2)=2\Re(r).
$$
But for $r=\infty$ we have
\begin{eqnarray*}
H^\infty &=& F^{-1}_\infty
\supset F^0_\infty=\dC\cdot A_2\supset F^1_\infty=\{0\},\\
W_0/W_{-1}&=&F^{-1}_\infty\Gr^W_0 \supset F^0_\infty \Gr^W_0 =\{0\},\\
W_{-2}/W_{-3}&=& W_{-2}=\dC\cdot A_2 =F^0_\infty
\supset F^1_\infty=\{0\},
\end{eqnarray*}
so here $W_0/W_{-1}$ carries a Hodge structure of weight $-2$,
and $W_{-2}/W_{-3}$ carries a Hodge structure of weight $0$.
Here $N$ is not strict,
$N(F^0_\infty)=\{0\}\neq N(H^\infty)\cap F^{-1}_\infty=\dC\cdot A_2=F^0_\infty.$
 So the filtration for $r=\infty$ is not at all part of a PMHS.

The classifying space
\begin{eqnarray*}
\DcPMHS = \{F^\bullet\subset H^\infty\ |\
F^1=0\subset F^0=\dC\cdot(iA_1+rA_2)\subset F^{-1}=H^\infty,r\in\dC\}
\cong\dC
\end{eqnarray*}
is compactified by $F^\bullet_\infty$ to $\dP^1$.

\subsection{Applications}
\label{subsecApplications}

In the remainder of this section, we use the results of sections \ref{secLimitTwistor} to \ref{secLimitProofs}
and the construction of the space $\MBLp$ to
prove some applications which are analogues of results for variations of Hodge structures.
They are concerned with extending variations of regular singular, pure polarized TERP-structures over subvarieties.
We first show that  such a variation defined outside a subset of codimension
at least two can be extended to the whole space. This uses the curvature computation
of \cite{HS2} as well as the construction of the compact classifying space. A second application
concerns extensions of variations of TERP-structures over codimension one
subvarieties, here we also use the extension results from the first part, namely theorem
\ref{theoLimitTwistor}.

We associated in \cite[lemma 4.4]{HS2} to any variation of regular singular
TERP-structures with constant spectral pairs a period map to a
classifying space $\cbl$. Here is the analogue if we do not suppose
that the spectral pairs are constant. We use the notion of ``regular
singular mixed TERP-structures'', introduced in definition \ref{defHodgeSpectrum}.
Recall also that $\cL$ denotes the universal locally free sheaf on the classifying space
$\MBL$.
\begin{lemma}
\label{lemPeriodMaps}
Let $\TERP$ be a variation of regular singular TERP-structures on a complex manifold $M$.
Let $H^\infty_\dR, S, M_z, w$ be its topological data
and $\alpha_1\in\dC$ such that $\cH\subset \cV^{\alpha_1}$ (see lemma \ref{lemRegSingBounded}).
\begin{enumerate}
\item
Then there is a unique period map $\widetilde{\phi}:\widetilde{M}\rightarrow \MBL^{H^\infty_\dR,S,M_z,w,\alpha_1}$, where
$\pi:\widetilde{M}\rightarrow M$ (as before, we write $\MBL$ for the target of $\widetilde{\phi}$).
\item
We have that $d\widetilde{\phi}(\cT_{\widetilde{M}})\subset \Theta_{\MBL}\cap{\cH}\!om_{\cO_{\dC\times\MBL}}(\cL,z^{-1}\cL/\cL)$,
and we say that $\widetilde{\phi}$ is horizontal.
\item
If $\Spp(H_{|\dC\times\{x\}},\nabla_z)=\Spp$ for all $x\in M$
for some fixed spectral pairs $\Spp$, then $\textup{Im}(\widetilde{\phi})\subset U_{\Spp}$
(which is equal to some $\cbl$ iff $\TERP$ is mixed TERP).
\item
The image of $\widetilde{\phi}$ is contained in $\MBLp$ if $\TERP$ is pure polarized.
If $\TERP$ has constant spectral pairs and is moreover mixed and pure polarized, then
$\widetilde{\phi}$ is distance decreasing with respect to the distance $d_h$ on
$\cblp\subset\MBLp$ and the Kobayashi pseudo-distance on $\widetilde{M}$.
\item
If $\TERP$ is pure polarized, and if we suppose moreover that
the monodromy representation $\gamma:\pi_1(\dC^*\times M)\rightarrow \Aut(H^\infty_\dR)$ respects
a lattice $H^\infty_\dZ\subset H^\infty_\dR$, then the period map $\widetilde{\phi}$ descends
to a locally liftable map $\phi:M\rightarrow \MBLp/G_\dZ$ where $G_\dZ:=\Aut(H^\infty_\dZ,S,M_z)$.
For constant spectral pairs, its image is contained in
$(U_{\Spp}\cap\MBLp)/G_\dZ$.
\end{enumerate}
\end{lemma}
\begin{proof}
Consider the variation $\pi^*\TERP$, this is obviously an element
of $\MBLfun^{H^\infty_\dR,S,M,w,\alpha}(\widetilde{M})=:\MBLfun(\widetilde{M})$, hence, it corresponds by theorem \ref{theoIsomorphismFunctors} to a unique morphism of complex
spaces $\widetilde{\phi}:\widetilde{M}\longrightarrow\MBL$, with the property
that $\widetilde{\phi}^*\cL\cong \pi^*\cH$ as families of TERP-structures.
The fact that $\cH$ underlies a variation of TERP-structures translates into the horizontality of $\widetilde{\phi}$.
All other statements are obvious consequences of the results of the last sections, the distance decreasing property
follows from \cite[theorem 4.1 and proposition 4.3]{HS2}.
\end{proof}
\begin{theorem}
Let $X$ be a complex manifold and $Y:=X\backslash Z$, where $Z$ is an analytic subspace
of codimension at least two. Then any variation of pure polarized,
regular singular mixed TERP-structures $\TERP$ on $Y$ with constant spectral pairs can be extended
to $X$.
\end{theorem}
\begin{proof}
We argue as in \cite[beginning of \textsection  4]{Sch} and \cite[chapters VI, VII]{Ko2}.
The extension problem is of local nature, therefore, for any $x\in X$, we choose a simply connected
neighborhood $V$ of $x$ in $X$. Then $U:=V\cap Y$ is also simply connected, as $Z$ has codimension at least
two in $X$. We obtain a period map
$\phi:U\rightarrow \cbl\cap\MBLp$.
Notice that by an induction argument on the dimension of the singular locus
of $Z$ we may in fact assume that $Z$ is smooth, and therefore it suffices
to consider the case where $V$ is a polycylinder $V=
\Delta^N$, and $U=\Delta^k\times(\Delta^*)^{N-k}$ for some $k\geq 2$.
\cite[chapter IV, corollary 4.5]{Ko2} yields that the Kobayashi pseudo-distance
$d_U$ is a true distance in this case.
By lemma \ref{lemPeriodMaps}, $\phi$ is distance-decreasing with
respect to the distance $d_h$ on $\MBLp$ and the distance $d_U$ on $U$.
As $\overline{\check{D}}_{\mathit{BL}}\cap\MBLp$ is complete with respect to $d_h$ by theorem \ref{theoCompleteness},
this implies that $\phi$ extends
continuously to the closure $\overline{U}_{d_U}$ of $U$ with respect to $d_U$.
The assumption that $Z\cap V$ is of codimension
at least two in $V$ implies (see \cite[VI, Proposition 5.1]{Ko2}) that the restriction
$(d_V)_{|U}$ agrees with $d_U$, so that $\overline{U}_{d_U}=V$.
This gives the extension $\overline{\phi}:V\rightarrow \overline{\check{D}}_{\mathit{BL}}\cap\MBLp$
we are looking for, which is necessarily holomorphic.
\end{proof}
\textbf{Remark:} Note that it follows from the construction of $\MBL$ that the extension constructed
in this way has in general jumping spectral numbers over the points
lying in $Z$.

\vspace*{0.3cm}

In applications, the extension over subvarieties of codimension one is an even more
important problem. For this, we can combine the limit results from section \ref{secLimitTwistor}
and the properties of the space $\MBLp$ to obtain the following statements
for the period map defined by a variation of pure polarized regular singular TERP-structures.

Let, as in section \ref{secLimitTwistor}, $1\leq l\leq n$, $X:=\Delta^n$, $Y:=(\Delta^*)^l\times\Delta^{n-l}, X\backslash Y=\coprod_{i\in \underline{l}}D_i$,
and consider a variation of pure polarized regular singular TERP-structures $\TERP$ on $Y$.
Denote by $M_i\in\Aut(H^\infty_\dR)$ the monodromy corresponding
to a loop around $\dC^*\times D_i\subset \dC^*\times X$.
As before, we say that the monodromy respects a lattice if
there is a lattice $H^\infty_\dZ\subset H^\infty_\dR$ such that the image of
$\gamma:\pi_1(\dC^*\times Y)\rightarrow \Aut(H^\infty_\dR)$
is contained in $\Aut(H^\infty_\dZ)$, in that case we put $G_\dZ:=\Aut(H^\infty_\dZ,S,M_z)$.
First we have the following rather
simple consequence of the relation between TERP-structures and twistor
structures.
\begin{corollary}
The eigenvalues of the automorphisms $M_i$ are elements in $S^1$. If the monodromy
respects a lattice, then they are roots of unity.
\end{corollary}
\begin{proof}
The first part has already been shown in the proof of lemma
\ref{lemKMSSforTERP}. The second part is the standard
argument known from the case of variations of Hodge structures: If
$M_i\in\Aut(H^\infty_\dZ)$, then its eigenvalues are algebraic integers,
so if they have absolute value one, they are necessarily roots of unity.
\end{proof}
The extension properties of the period map alluded to above can be stated as follows.
\begin{theorem}\label{theoExtensionTERPDivisor}
Let $\TERP$ be a variation of regular singular, pure polarized TERP-structures on $Y$.
\begin{itemize}
\item
If $M_i=\Id$ for all $i\in\{1,\ldots,k\}$, i.e., if there is a period map
$\phi:Y \rightarrow \MBLp$,
then this map extends to
$$
\overline{\phi}:X \rightarrow \MBLp
$$
and the variation $H$ extends to a variation on $X$.
\item
Suppose that the monodromy respects a lattice, so that we have
a locally liftable period map $\phi:Y \rightarrow \MBLp/G_\dZ$.
If all $M_i$ are semi-simple, then $\phi$
extends holomorphically
(not necessarily locally liftable) to
$$
\overline{\phi}:X \rightarrow \MBLp/G_\dZ.
$$
\end{itemize}
In particular, given a $\mu$-constant family of isolated hypersurface
singularities $F:(\dC^{n+1}\times Y,0)\rightarrow (\dC,0)$, the above statements
apply if the variation $\mathit{TERP}(F)$ on $Y$ described in \cite{He4} and \cite{HS1}
is pure polarized. In this case, the extension of the period map is
contained in $\overline{\check{D}}_{BL}\cap \MBLp$ resp.
$(\overline{\check{D}}_{BL}\cap \MBLp)/G_\dZ$ where $\cbl$ is the
classifying space associated to the variation of mixed TERP-structures $\mathit{TERP}(F)$.
\end{theorem}
\begin{proof}
In both cases, for any $x\in D$ consider the maximal subset $I \subset \underline{l}$ such that
$x\in D_I$. Theorem \ref{theoLimitTERP} yields the limit TERP-structure $\cH(x):=\cH(I)_{|\dC\times\{x\}}\in\VB_{|\dC\times\{x\}}$, which
is pure polarized by theorem \ref{theoLimitTwistor} as all nilpotent parts $N_i\in\End(H^\infty_\dR)$ are
zero. By proposition \ref{propLimitTERPregsing}, the smallest spectral number of $\cH(x)$ is not smaller than $\alpha_1$,
so that $[\cH(x)/V^{>\alpha_\mu-1}]\in\MBLp$.
In the first case, putting $\overline{\phi}(x):=[\cH(x)/V^{>\alpha_\mu-1}]\in\MBLp$
defines a continuous and hence holomorphic extension $\overline{\phi}:X \rightarrow \MBLp$.
In the second case we can choose by the last corollary a sufficiently large positive integer $m$
such that $M_i^m=\Id$ for all $i$. Then the lifted map $\widetilde{\phi}(r_1,\ldots,r_l,r_{l+1},\ldots,r_n):=\phi(r_1^m,\ldots,r_l^m,r_{l+1},\ldots,r_n)$
extends as in the first part to a map $\overline{\widetilde{\phi}}:\Delta^n \longrightarrow \MBLp$,
which yields the extension $\overline{\phi}:X\rightarrow\MBLp/G_\dZ$ we are looking for.
\end{proof}

\textbf{Examples:} The following two examples, borrowed from the classifying spaces
$\MBL$ and their strata $U_{\Spp}$ from the last subsections illustrate what kind of phenomena can
occur when extending families of TERP-structures over boundary divisors.
\begin{enumerate}
\item
In subsection \ref{subsecWPS} a variation of TERP-structures on $\dC^2$ with parameters $(r,t)$
and constant spectrum $(\alpha_1,0,-\alpha_1)$ was considered.
For example, its restriction to $\{0\}\times \Delta$ is pure and polarized.
The further restriction to $\{0\}\times \Delta^*$ extends to $(0,0)$ by theorem \ref{theoExtensionTERPDivisor}, 1.,
and gives there a pure and polarized TERP-structure, which is of course the original one at $(0,0)$.

One can also restrict to the variation on $\dC^*\times \{0\}$ with the parameter $\www r =\frac{1}{r}$.
It is pure and polarized for $|r|\neq 1$. By theorem \ref{theoExtensionTERPDivisor}, 1., it has a
pure and polarized limit TERP-structure for $\www r\to 0$. That had also been calculated in
subsection \ref{subsecWPS}, its spectral numbers are $(-\alpha_1-2,0,\alpha_1+2)$, so here the spectral numbers jump.
\item
Now we show an easy example where the second part of theorem \ref{theoExtensionTERPDivisor}
can be applied. Consider the following topological data: Let $H^\infty=\dC A_1\oplus\dC A_2$, $\overline{A}_1=A_2$,
$M_z(A_i)=A_i$ and $S(A_i,A_j)=\delta_{i+j,3}$. Put $w=0$ and $\alpha:=-1$, then the classifying space
$\MBL$ for these data consists of two components of the space considered in
subsection \ref{subsecRedSpaces}, i.e., $\MBL\cong \dP_r^1\cup\dP_s^1$, with the following universal families
$$
\begin{array}{rcl}
\cH^{(r,1)} & := & \cO_{\dC\times\dC_r} (z^{-1}A_1+rA_2)\oplus\cO_{\dC\times\dC_r}zA_2\quad\textup{ over }\quad\dC_r\subset\dP_r^1, \\ \\
\cH^{(r,2)} & := & \cO_{\dC\times(\dP_r^1\backslash\{0\})} (r^{-1}z^{-1}A_1 + A_2)\oplus\cO_{\dC\times(\dP_r^1\backslash\{0\})}A_1 \quad\textup{ over }\quad\dP_r^1\backslash\{0\}, \\ \\
\cH^{(s,1)} & := & \cO_{\dC\times\dC_s} (z^{-1}A_2+sA_1)\oplus\cO_{\dC\times\dC_s}zA_1\quad\textup{ over }\quad\dC_s\subset\dP_s^1, \\ \\
\cH^{(s,2)} & := & \cO_{\dC\times(\dP_s^1\backslash\{0\})} (s^{-1}z^{-1}A_2 + A_1)\oplus\cO_{\dC\times(\dP_s^1\backslash\{0\})}A_2 \quad\textup{ over }\quad\dP_s^1\backslash\{0\}, \\ \\
\end{array}
$$
where the TERP-structures corresponding to $r=\infty$ and $s=\infty$ are the same, i.e., the common point of the two components of $\MBL$.
In subsection \ref{subsecRedSpaces} it is shown that $\MBLp=\left\{r\in\dC\,|\,|r|<1\right\}\coprod\left\{s\in\dC\,|\,|s|<1\right\}$.
Define $H^\infty_\dZ:=\dZ\cdot \frac{A_1+A_2}{2}\oplus \dZ\cdot i\frac{A_1-A_2}{2}$.
It is easy to see that $G_\dZ:=\Aut(H^\infty_\dZ,M_z,S)=\Aut(H^\infty_\dZ,S)=D_4$, and that
the group action of $G_\dZ$ on $\MBLp$ identifies the two components $\{r\in\dC\,|\, |r|<1\}$ and
$\{s\in\dC\,|\, |s|<1\}$, and quotients once more by $r\mapsto -r$, so that the quotient space
is still an open disc $\Delta$, with coordinate $\widetilde{r}=r^2$ (resp., $\widetilde{r}=s^2$).

Now consider the following variation over $Y:=\dC^*$: $\cH:=\cO_{\dC\times Y} (z^{-1}q^{-3/4}A_1+q^{3/4}A_2) \oplus
\cO_{\dC\times Y} (z q^{3/4}A_2)$.
Here $M_q(\underline{A})=\underline{A}\cdot\diag(-i,i)$, so that $\mathit(Im)(\gamma)
\cong \dZ/4\dZ$ is contained in $G_\dZ$
(remember that $\gamma:\pi_1(\dC^*\times Y)\cong\dZ^2\rightarrow \Aut(H^\infty_\dZ)$ is the monodromy representation).
The restriction of this family to $|q|<1$ is pure polarized.
%We have the period map
%$\phi:\widetilde{Y}\rightarrow \MBL/G_\dZ$ which sends $\widetilde{q}$ to $e^{3\pi i \widetilde{q}}$ (here $\widetilde{Y}\cong\dC$ is the universal cover
%of $Y$ and $\widetilde{q}$ a coordinate on $\widetilde{Y}$, with $q=e^{2\pi i \widetilde{q}}$).
We have the period map
$\phi: Y\rightarrow \MBL/G_\dZ$ given by $q\mapsto \widetilde{r}=q^3$.
According to theorem \ref{theoExtensionTERPDivisor}, 2.,
we obtain a holomorphic extension $\overline{\phi}:\Delta\rightarrow \MBLp/G_\dZ
=\Delta$, still given by $\widetilde{r}=r^3$, which is obviously not
locally liftable. Notice finally that both the members of the family
$\cH$ and the image of $0\in\Delta$ of the
extended period map $\overline{\phi}$ are TERP-structure
with spectral numbers $(-1,1)$, so that $\mathit{Im}(\overline{\phi})$ is actually contained in
$(U_{(-1,1)}\cap\MBLp)/G_\dZ$.
\end{enumerate}

\textbf{Remarks:} Pursuing further the analogy with the theory of period maps for variations
of Hodge structures, one might ask whether the asymptotic behavior of the above defined
map $\widetilde{\phi}$ can be controlled by the so called nilpotent orbits of TERP-structures,
as studied in \cite{HS1}. More precisely, given a variation of regular singular, pure polarized TERP-structures $\TERP$
(on $\Delta^*$, say),
one might consider the family $G:=\pi^*(K)$, where $K:=(\lim_{r\rightarrow 0}H) \in \MBL$
and $\pi:\dC\times\Delta^*\rightarrow \dC,(z,r)\mapsto zr$.
This is also a variation over $\Delta^*$, with $G_{|r=1}=K$.
It seems reasonable to expect
that $G$ is a nilpotent orbit of TERP-structures, i.e.,
it lies in $\MBLp$ for $|r|\ll 1$ (see also \cite[theorem 6.6]{HS1}).
Then we can consider the distance of the two families, and ask whether an estimate
as in \cite[theorem 4.9]{Sch} holds. This is particularly interesting if
$K$ has different spectrum  than the general member of $H$, as in this case
these two families are in different strata of the space $\MBL$.

One might also be interested to work out such an estimate of the asymptotical behavior of
the distance between $H$ and $G$ in the higher dimensional case, as in \cite[theorem 4.12]{Sch}.

\bibliographystyle{amsalpha}
\providecommand{\bysame}{\leavevmode\hbox to3em{\hrulefill}\thinspace}
\providecommand{\MR}{\relax\ifhmode\unskip\space\fi MR }
\providecommand{\MRhref}[2]{%
  \href{http://www.ams.org/mathscinet-getitem?mr=#1}{#2}
}
\providecommand{\href}[2]{#2}

\vspace*{1cm}

\nd
Lehrstuhl f\"ur Mathematik VI \\
Institut f\"ur Mathematik\\
Universit\"at Mannheim,
A 5, 6 \\
68131 Mannheim\\
Germany

\vspace*{1cm}

\nd
hertling@math.uni-mannheim.de\\
sevenheck@math.uni-mannheim.de

\end{document}